\documentclass[a4paper,11 pt]{amsart}
\allowdisplaybreaks

 \usepackage[sc]{mathpazo}
\usepackage{eulervm}
\usepackage{hyperref}

\usepackage[utf8]{inputenc}
\usepackage[T1]{fontenc}

\usepackage[english]{babel}
\usepackage{mathtools}
\usepackage{braket}
\usepackage{amsmath, amssymb, stmaryrd, mathabx}
\usepackage{enumerate, enumitem}
\usepackage{tikz-cd}
\usepackage{tikz,float, hyperref, circuitikz}
\usepackage{titlesec, cite}

\titleformat{\section}
{\normalfont\normalsize\bfseries}{\thesection.}{1em}{}
\titleformat{\subsection}
{\normalfont\normalsize\itshape}{\thesubsection.}{1em}{}
\titleformat{\subsubsection}
{\normalfont\normalsize\itshape}{\thesubsubsection.}{1em}{}

\newcommand{\A}{\mathfrak{A}}
\newcommand{\Hs}{\mathcal{H}}
\newcommand{\Ps}{\mathcal{P}}

\newcommand{\TT}{{\mathbb{T}}}

\newcommand{\Bases}{\mathcal{B}}
\newcommand{\bb}{\mathbf{B}}
\newcommand{\bta}[1]{\beta^{[#1]}}
\newcommand{\iber}{\operatornamewithlimits{iRes}}
\newcommand{\ires}{\operatornamewithlimits{iRes}}

\newcommand{\Jac}{\mathrm{Jac}}
\newcommand{\Sym}{\mathrm{Sym}}
\newcommand{\ch}{\mathrm{ch}}

\newcommand{\res}{\mathrm{Res}}
\newcommand{\ParEnd}{\mathrm{ParEnd}}
\newcommand{\End}{\mathrm{End}}
\newcommand{\SParEnd}{\mathrm{SParEnd}}

\newcommand{\ParHom}{\mathrm{ParHom}}
\newcommand{\SParHom}{\mathrm{SParHom}}
\newcommand{\lhat}{\widehat{\lambda}}
\newcommand{\iio}{i_1}
\newcommand{\iit}{i_2}
\newcommand{\jjo}{j_1}

 \theoremstyle{plain}
 \newtheorem{theorem}{Theorem}[section]
\newtheorem{proposition}[theorem]{Proposition}
 \newtheorem{lemma}[theorem]{Lemma}
 
\newtheorem{corollary}[theorem]{Corollary}
\theoremstyle{definition}
 \newtheorem{definition}[theorem]{Definition}

 \newtheorem{remark}[theorem]{Remark}
  
  \newtheorem{construction}[theorem]{Construction}

\title[Equivariant Euler characteristics on parabolic $\mathrm{SL}_2$- and $\mathrm{SL}_3$-Higgs moduli spaces]{Residue formulas for equivariant Euler characteristics on parabolic $\mathrm{SL}_2$- and $\mathrm{SL}_3$-Higgs moduli spaces}
\author{Olga Trapeznikova}
\address{ Faculty of Mathematics, ISTA}
\email{Olga.Trapeznikova@ista.ac.at}
\begin{document}	

\begin{abstract} 
We obtain explicit residue formulas for the equivariant Euler characteristics of line bundles and their twists by tautological $K$-theory classes on moduli spaces of stable strongly parabolic $\mathrm{SL}_2$- and $\mathrm{SL}_3$-Higgs bundles on a curve. For line bundles, these formulas provide an alternative to the
equivariant Verlinde formula, which takes the form of a finite sum over the solutions of Bethe-type equations.  Our approach is geometric and relies on equivariant localization, which reduces the computation to contributions from the components of the fixed-point locus. These are evaluated using residue formulas on moduli spaces of stable parabolic bundles, the intersection theory of symmetric products of the curve, and, in rank $3$, a study of moduli spaces of stable parabolic triples through variation of stability and wall-crossing.
\end{abstract}
\maketitle

\section{Introduction}\label{sec:intro}
Let $C$ be a compact Riemann surface of genus $g\ge2$ with canonical bundle $K$,
and fix a point $p\in C$. For a regular vector $c=(c_1>\dots>c_r)$ of parabolic
weights (see \S\ref{S:notation}), let $\Ps_0(c)$ denote the moduli space of stable
parabolic $\mathrm{SL}_r$-bundles on $C$ with parabolic structure given by a full
flag at $p$ with weights $c$. This space is a smooth projective variety carrying a
natural family of line bundles $\mathcal L(k;\lambda)$, indexed by a positive
integer $k$ and an integral vector $\lambda=(\lambda_1,\dots,\lambda_r)$ with
$\sum_i\lambda_i=0$. When $c=\lambda/k$, the line bundle $\mathcal L(k;\lambda)$ is
ample, and its Euler characteristic is given by the \textit{parabolic Verlinde
formula} \cite{Ver}, which in its standard form is a finite sum over certain
elements of the maximal torus of $\mathrm{SL}_r$. Alternatively, this Euler
characteristic can be expressed as an iterated residue of an explicit meromorphic
function \cite{SzDu,OTASz}; this residue expression remains valid for every regular
$c$ and all $k$ and $\lambda$, and makes the dependence on $k$ and $\lambda$
explicit.

These invariants have natural analogues for parabolic Higgs bundles.
Let $\Hs_0(c)$ denote the moduli space of stable strongly parabolic
$\mathrm{SL}_r$-Higgs bundles, whose points are pairs $(V,\Phi)$, where $V$ is a
parabolic $\mathrm{SL}_r$-bundle as above and
$\Phi\colon V\to V\otimes K(p)$ is a traceless Higgs field whose residue at $p$ is
nilpotent with respect to the parabolic flag. The space $\Hs_0(c)$ contains
$\Ps_0(c)$ as the closed subvariety of pairs of the form $(V,0)$, and the line
bundles $\mathcal L(k;\lambda)$ extend naturally to $\Hs_0(c)$. Although
$\Hs_0(c)$ is not compact, it carries a $\mathbb C^*$-action rescaling the Higgs
field, $t\cdot(V,\Phi)=(V,t\Phi)$, with compact fixed-point locus. The equivariant
Euler characteristic $\chi_t(\Hs_0(c),\mathcal L(k;\lambda))$ can thus be defined
via localization \cite{Thomason,Okounkov} and expanded as a Laurent
series in $t^{-1}$. 

Motivated by physics, Gukov and Pei \cite{GP} proposed a formula for this series,
which was proved by Andersen, Gukov and Pei \cite{AGP}, building on the index
theorem of Teleman and Woodward \cite{TW}, and independently by Halpern-Leistner
\cite{HL}.  Like its classical counterpart, this \textit{equivariant Verlinde formula}
is valid in arbitrary rank and takes the form of a finite sum.
The summation points, however, are now given implicitly as solutions of
Bethe-type equations depending on $t$, which makes the sum more difficult
to evaluate directly. This motivates the search for an explicit residue
expression analogous to the classical one.

For moduli spaces of rank-$2$ degree-$1$ Higgs bundles without parabolic
structure, such a residue expression can be deduced from the equivariant
intersection numbers computed in \cite{ChHSz}. In this paper, we extend this result to moduli spaces of strongly parabolic Higgs
bundles of ranks $2$ and $3$, and from line bundles to their twists by
tautological $K$-theory classes, for which residue formulas on $\Ps_0(c)$
were obtained in \cite{Unyu}. To define these classes, let $\nu$ be a
dominant weight of $\mathrm{GL}_r$ and let $U_\nu$ denote the bundle on $\Hs_0(c)\times C$ associated
with the normalized universal bundle (see \S\ref{S:notation}) via the irreducible
representation $\rho_\nu$ of $\mathrm{GL}_r$ with highest weight $\nu$. Fixing a
square root $K^{1/2}$ of $K$, we consider the tautological class
$\pi_!(U_\nu\otimes K^{1/2})$, where $\pi\colon\Hs_0(c)\times C\to\Hs_0(c)$ is the
projection. Our main results, Theorems~\ref{Rank2Result} and~\ref{thm:main},
express both
$$
\chi_t(\Hs_0(c),\mathcal L(k;\lambda))
\text{\;\;\; and \;\;\;}
\chi_t\bigl(\Hs_0(c),\mathcal L(k;\lambda)\otimes
\pi_!\bigl(U_\nu\otimes K^{1/2}\bigr)\bigr)
$$
as (iterated) residues of explicit meromorphic functions. To illustrate the form
of these expressions, we state below the rank $2$ result for line bundles.
\begin{theorem}\label{thm:intro}
Let $c\in \bigl(0,\tfrac12\bigr)$. Given  $\lambda\in\mathbb Z$ and $k\in\mathbb N$, define 
$$\Omega:=\frac{2^{g}\,e^{\lambda y}\,\bigl(k+2-\frac{2}{1-e^{u+y}}-\frac{2}{1-e^{u-y}}\bigr)^g\,\mathrm{d}y}{W(y)^{g-1}(1-e^{-y})(1-e^{y-u})\big(1-B(y)\big)}, \text{\; where \;} B(y)=e^{(k+2)y}\Big(\tfrac{1-e^{-y-u}}{1-e^{y-u}}\Big)^2$$ and 
$W(y)=(1-e^{-u})(1-e^{y})(1-e^{-y})(1-e^{y-u})(1-e^{-y-u})$. Then, with $t=e^u$,  
\begin{equation*}\chi_{t}(\Hs_0((c,-c)), \mathcal{L}(k;(\lambda,-\lambda))) =
 \underset{y=0}{\mathrm{Res}}\,\Omega + \underset{y=-u}{\mathrm{Res}}\!\!{}^+\Omega  + \underset{y=u}{\mathrm{Res}}\!{}^-\Omega.
\end{equation*}
Here the superscripts $+$ and $-$ indicate that the factor $(1-B(y))^{-1}$ in $\Omega$ is expanded as a geometric series in nonnegative, respectively negative, powers of
$B(y)$.
\end{theorem}  

To compare Theorem \ref{thm:intro} with the equivariant Verlinde formula of \cite{AGP}, note
that the equation $B(y)=1$ is, up to normalization, the rank $2$ Bethe equation
of \cite[(58)]{AGP}. Rather than summing over its solutions, 
our formula involves residues at the three points $y=0$ and $y=\pm u$,
which are independent of the level $k$.

\medskip
\noindent\textbf{Organization of the paper and our approach.}
In \S\ref{S:NotRank2}, we prove Theorem \ref{thm:intro} together with
its extension to twists of line bundles by tautological $K$-theory
classes (Theorem \ref{Rank2Result}). For these twists, the same residues
are applied to a first-order deformation of $\Omega$ by the character
of $\rho_\nu$. The rank $3$ formulas, stated in Theorem \ref{thm:main}
in \S\ref{S:mainrank3}, have an analogous structure, with iterated residues
in two variables.

To obtain these formulas, we use a geometric approach based on equivariant
localization and basic geometric invariant theory. We express the
equivariant Euler characteristic as a sum over the connected components
of the $\mathbb C^*$-fixed locus of $\Hs_0(c)$ and compute their contributions
separately. One of these components is $\Ps_0(c)$ itself, whose contribution
we compute in \S\ref{S:3} using \cite{Unyu}; the resulting expressions
recover the classical formulas of \cite{OTASz,Unyu} after taking the constant
term in $t^{-1}$. The fixed components consisting of Higgs bundles whose
underlying bundle is a direct sum of line bundles are finite \'etale
covers of products of symmetric powers of $C$; in \S\ref{S:111} we reduce
the computation of their contributions to the classical intersection
theory of symmetric products.

In rank $3$, there are also fixed components consisting of Higgs bundles
whose underlying bundle splits as a direct sum of a rank-$2$ bundle and
a line bundle. These are moduli spaces of \textit{stable parabolic triples}.
Following \cite{GPGM}, we study their geometry in \S\ref{S:partriplesgen}
by varying the stability parameter. In \S\ref{S:strategy}, we describe
in detail our strategy for computing their contributions to the
equivariant Euler characteristics, and in \S\ref{S:21and12} we carry out
the computations using the wall-crossing technique of \cite{OTASz}.

\medskip
\noindent\textbf{Further remarks.} A notable feature of our formulas is that their individual residue terms have a clear geometric interpretation.  In Theorem \ref{thm:intro}, for
instance, the residue at $y=0$ gives the contribution of $\Ps_0(c)$ to
the equivariant Euler characteristic, while each of the residues at
$y=\pm u$ gives the total contribution of one of the two families of
fixed components with nonzero Higgs field, distinguished by their
parabolic structure at $p$.

Our results may also be relevant to the study of mirror symmetry for
Higgs bundles. The mirror relationship between the moduli spaces of
$\mathrm{SL}_3$- and $\mathrm{PGL}_3$-Higgs bundles has been investigated in
\cite{HauselThaddeus,DonagiPantev}, and topological mirror symmetry for their
strongly parabolic counterparts has been established in
\cite{GothenOliveira,Shen}.
It would be interesting to understand how the equivariant Euler characteristics
computed here, particularly those involving tautological classes, are
reflected on the $\mathrm{PGL}_3$ side.

\medskip
\noindent\textbf{Acknowledgments.}
I am grateful to András Szenes and Tamás Hausel for sharing their ideas and for helpful discussions, and to Tamás Hausel for his encouragement and support.
This research was supported by the European Union's Horizon 2020 research and innovation programme under the Marie Sk\l{}odowska-Curie Grant Agreement No. 101034413 and by the Austrian Science Fund (FWF) 10.55776/ESP2557325.

\section{Moduli spaces of parabolic Higgs bundles and residue formulas in rank $2$}\label{S:NotRank2}

\subsection{Moduli spaces of stable parabolic Higgs bundles}\label{S:notation}
We begin by introducing the moduli spaces of stable parabolic bundles and parabolic Higgs bundles that will be used throughout the paper, and recall the basic properties that we will need later.  For further details, we refer to  \cite{MehtaSeshadri, OTASz, Yokogawa, Thaddeus, GPGM, GothenOliveira}. 

Fix a point $p\in C$. A \textit{parabolic bundle} on $C$ is a vector bundle $V$ of rank $r$, together with a full flag $V_p=F_r\supsetneq F_{r-1}\supsetneq \ldots\supsetneq F_1\supsetneq F_0=\{0\}$ and a collection of \textit{parabolic weights} $c=(c_1,c_2,\ldots,c_r)$  assigned to the successive quotients of the flag and satisfying the conditions 
\begin{equation}\label{eq:ccond}
c_1>c_2>\ldots>c_r \text{\; and \:} c_1-c_r<1.
\end{equation}
The \textit{parabolic degree}  and \textit{parabolic slope} of $V$  are defined as $$\mathrm{pardeg}(V)=\mathrm{deg}(V)-\sum_{i=1}^r c_i \text{\; and \;} \mathrm{par}\mu(V)= \frac{\mathrm{pardeg}(V)}{r}. $$  Every subbundle $V'\subset{V}$  inherits a parabolic structure from $V$ in a natural way \cite{MehtaSeshadri}, and  $V$ is called \textit{stable}, if  $\mathrm{par}\mu(V')<\mathrm{par}\mu(V)$ for every proper subbundle $V'\subset{V}$.	
 
In this paper, we restrict to the degree-zero case and denote by $\Delta(r)$ the set of vectors $c=(c_1,\ldots,c_r)\in\mathbb Q^r$ satisfying \eqref{eq:ccond} and $\sum_i c_i=0$. We call a vector $c\in\Delta(r)$ \textit{regular} if  no nontrivial subset of its coordinates sums to an integer.  For every regular $c\in\Delta(r)$, there is a smooth projective variety $\Ps_0(c)$ parameterizing isomorphism classes of stable rank-$r$ parabolic bundles with weights $c$ and trivial determinant.

We next recall the notions needed to pass from parabolic bundles to parabolic Higgs bundles. Let $$(V', F'_\bullet \subset V'_p,  (c_1',\ldots,c_{r'}')) \text{\; and \;} (V, F_\bullet\subset V_p, (c_1,\ldots,c_r)),$$ be parabolic bundles of ranks $r'$ and $r$, respectively. A morphism $f\colon V'\to V$ is called \textit{parabolic} if  $f(F'_i)\subset F_{j-1}$ for $c'_{r'-i+1}< c_{r-j+1}$ and  \textit{strongly parabolic} if the same condition holds for $c'_{r'-i+1}\leq c_{r-j+1}$.  We denote the corresponding sheaves of morphisms by $\ParHom(V',V)$ and $\SParHom(V',V)$, and write $\ParEnd_0(V),$ $\SParEnd_0(V)$ for the sheaves of traceless parabolic and strongly parabolic endomorphisms, respectively.

Now let $K$ denote the canonical bundle on $C$. A \textit{parabolic Higgs bundle} is a pair $(V,\Phi)$, where $V$ is a parabolic bundle and $\Phi\in H^0(\SParEnd(V)\otimes K(p))$ is a \textit{Higgs field}.  A parabolic subbundle $V'\subset V$ is called $\Phi$-invariant if $\Phi(V')\subset V'\otimes K(p),$ and the pair $(V,\Phi)$ is called \textit{stable} if $\mathrm{par}\mu(V')<\mathrm{par}\mu(V)$ for every proper nonzero $\Phi$-invariant parabolic subbundle $V'\subset V$.
	
For regular parabolic weights, stable parabolic Higgs bundles admit a well-behaved moduli space. We collect below the properties of this space that will be used later. 
\begin{theorem}\label{thm:modsp} Let $c\in\Delta{(r)}$ be a regular weight vector. Then:
\begin{enumerate}[label=(\roman*), leftmargin=1.5pc]
\item There exists a smooth quasi-projective variety $\Hs_0(c)$ parameterizing isomorphism classes of stable rank-$r$ parabolic Higgs bundles of weight $c$, with trivial determinant and trace-free Higgs field. Moreover, on $\Hs_0(c)\times C$ there exists a universal bundle $\mathcal{U}$, together with a universal flag $$\mathcal{U}\big|_{\Hs_0(c)\times\{p\}} = \mathcal{F}_r\supset \mathcal{F}_{r-1}\supset \ldots \supset \mathcal{F}_1.$$ The universal bundle is unique up to tensoring by the pullback of a line bundle from $\Hs_0(c)$.

\item On $\Hs_0(c)$ there is a natural action of $\TT=\mathbb C^*$ given by rescaling the Higgs field, $t\cdot(V,\Phi)=(V,t\Phi)$, and $\Hs_0(c)$ is semi-projective with respect to this action. 
\item  The $\TT$-equivariant deformation complex of a parabolic Higgs bundle $(V,\Phi)$ is
\[
\ParEnd_0(V)\longrightarrow t\,\SParEnd_0(V)\otimes K(p),
\]
placed in degrees $0$ and $1$, where the differential is given by the commutator with $\Phi$, and $t$ denotes the standard character of $\TT$. Its first hypercohomology is naturally identified with $T_{(V,\Phi)}\Hs_0(c)$, and at a $\TT$-fixed point, this identification is $\TT$-equivariant.
\item  Setting the Higgs field equal to zero defines a natural closed embedding $\Ps_0(c)\hookrightarrow \Hs_0(c)$. 
Under this embedding, the universal bundle and the universal flag on $\Hs_0(c)\times C$ restrict to the corresponding universal objects on $\Ps_0(c)\times C$.
\end{enumerate}
\end{theorem}
\noindent\textbf{Notation:} We say that a universal bundle $\mathcal{U}$ is normalized if $\mathcal{F}_1\subset\mathcal{U}|_{\Hs_0(c)\times\{p\}}$ is equivariantly trivial. We denote such a normalized universal bundle by $U$.

\subsection{Vector bundles and equivariant Euler characteristics on $\Hs_0(c)$}
Let $\pi\colon \Hs_0(c)\times C\to \Hs_0(c)$ denote the natural projection. Given an integer $k\in\mathbb{Z}$ and  a vector $\lambda=(\lambda_1,...,\lambda_r)\in\mathbb{Z}^r$, such that $\sum_{i=1}^r\lambda_i=0$, we define the line bundle 
$$ \mathcal{L}(k; \lambda)=\mathrm{det}\left(\mathcal{U}\big|_{\Hs_0(c)\times\{p\}}\right)^{k(1-g)}\otimes\mathrm{det}(\pi_!\mathcal{U})^{-k} \otimes(\mathcal{F}_r/\mathcal{F}_{r-1})^{\lambda_1}\otimes...\otimes(\mathcal{F}_1)^{\lambda_r} $$ on  $\Hs_0(c)$. A simple calculation shows that $\mathcal{L}(k;\lambda)$ is independent of the choice of the  universal bundle $\mathcal{U}$. 
Later, we will use the following result.
\begin{lemma}\cite[Lemma 2.8]{OTASz}\label{c_2End2}
In the notation introduced above, $c_1(\mathcal{L}(k; 0)) = \frac{k}{2r}\pi_*c_2(\mathrm{End}(\mathcal{U}))$.
\end{lemma}

Now let $\nu = (\nu_1,....,\nu_r)$ be a dominant weight  of $\mathrm{GL}_r$, and let  $\rho_\nu$ denote the irreducible representation with highest weight $\nu$.
We write $U_\nu$ for the vector bundle on $\Hs_0(c)\times C$ associated to the normalized universal bundle $U$ and $\rho_\nu$. It admits the following explicit description: let $f\colon\mathrm{Flag}(U)\to\Hs_0(c)\times C$ be the full flag bundle and denote by $L_1,\dots,L_r$ its tautological quotient line bundles, then 
\begin{equation}\label{constructionUnu}
U_\nu = f_*(L_1^{\nu_1}\otimes L_2^{\nu_2}\otimes...\otimes L_r^{\nu_r}).
\end{equation} 
To $U_\nu$ we associate the $K$-theory class $\pi_!(U_\nu\otimes K^{1/2})$ on $\Hs_0(c)$, where $K^{1/2}$ enters only through its Chern character; in particular, the Grothendieck--Riemann--Roch theorem gives 
\begin{equation}\label{GRRUnyu}
\ch\bigl(\pi_!\bigl(U_\nu\otimes K^{1/2}\bigr)\bigr)=\pi_*\ch(U_\nu).
\end{equation}
 
Our goal is to evaluate explicitly, in ranks \(2\) and \(3\), the \textit{equivariant Euler characteristics}
\[
\chi_t\bigl(\Hs_0(c),\mathcal L(k;\lambda)\bigr) \;\;\; \text{and} \;\;\;
\chi_t\bigl(\Hs_0(c),\mathcal L(k;\lambda)\otimes\pi_!(U_\nu\otimes K^{1/2})\bigr).
\]

Using the equivariant localization formula \cite{AB68, Thomason}, these Euler characteristics can be expressed as sums of contributions from the connected
components of the fixed-point locus. More precisely,   for any $\TT$-equivariant \(K\)-theory class \(\mathcal V\) on \(\Hs_0(c)\), we have 
\begin{equation}\label{Euler}
\chi_t(\Hs_0(c),\mathcal V)
=
\sum_{F} \int_F \frac{\ch_\TT(\mathcal V|_F)}
{e_\TT(N_F^\vee)}
\mathrm{Td}(F),
\end{equation}
where the sum runs over the connected components \(F\) of the fixed-point locus \(\Hs_0(c)^\TT\), and \(e_\TT(N_F^\vee)\) denotes the equivariant Chern character of the $K$-theoretic Euler class of the conormal bundle $N_F^\vee$.  By a result of Simpson \cite[Theorem 8]{Simpson} (see also \cite[Proposition 3.4]{GPGM}), these fixed-point components admit the following description.

\begin{proposition}\label{fixedcomp}
A parabolic Higgs bundle \((V,\Phi)\in\Hs_0(c)\) is fixed under the \(\TT\)-action if and only if  $V$ decomposes as a direct sum $V\simeq V_1\oplus V_2\oplus\ldots\oplus V_m$ of parabolic bundles, such that $\Phi|_{V_i}\in H^0(\SParHom(V_i,V_{i+1})\otimes {K}(p))$ is nonzero for $i=1,\ldots,m-1$, and $\Phi|_{V_m} = 0$.  
\end{proposition}
\noindent\textbf{Notation:}\label{not:type} Given an ordered partition $\mu=(\mu_1,...,\mu_m)$ of the set $\{1,2,...,r\}$ we say that a fixed point $(V,\Phi)\in\Hs_0(c)^\TT$ with a decomposition $V\simeq V_1\oplus\ldots\oplus V_m$ as in Proposition \ref{fixedcomp} has \textit{parabolic type} $\mu$, if the parabolic weights of  $V_i$ are $\{c_{j}\}_{j\in\mu_i}$ for $1\leq i \leq m$. When only the ranks of the summands are relevant, we say that  $(V,\Phi)$ has \textit{type} $(|\mu_1|,|\mu_2|,...,|\mu_m|)$. 
\begin{remark}\label{rem:PtoH}
Note that the fixed-point locus of type $(r)$, or, equivalently, of parabolic type \(\mu=(\{1,\ldots,r\})\), forms the connected component \(\Ps_0(c)\subset\Hs_0(c)^\TT\).
 \end{remark}

\subsection{Residue formulas in rank $2$}
We now briefly discuss the rank-\(2\) case, extending the calculation of \cite{ChHSz}, where the equivariant intersection numbers on the moduli space of rank-\(2\) degree-$1$ Higgs bundles were computed. Although some of the difficulties arising in rank \(3\) are absent here, this case provides a useful illustration of our approach. For those computations that are simpler versions of their rank-$3$ counterparts, we omit the details and refer instead to the more detailed treatment in \S\ref{S:3}--\S\ref{S:111}.

A simple calculation shows that, in rank \(2\), the stability condition for parabolic Higgs bundles is independent of the parabolic weights  \(c=(c_1,c_2)\in\Delta(2)\).  Thus all such weight vectors \(c\) give the same moduli space, which we denote by  \(\Hs_0(2)\).
To describe its \(\TT\)-fixed-point locus,  we set $S^{\mathrm{id}}:=\varnothing$ and $S^{(12)}:=\{p\}$, and define $s^\tau:=|S^\tau|$ for $\tau\in\Sigma_2$.
Then by Proposition \ref{fixedcomp}, every \(\TT\)-fixed point \((V,\Phi)\in\Hs_0(2)\) has  one of the following parabolic types:
\begin{enumerate}[label=(\alph*), leftmargin=1.5pc]
\item \emph{Parabolic type \((\{1,2\})\):} \(V\) is a stable parabolic bundle and \(\Phi=0\). These points form the connected component isomorphic to the moduli space of stable $\mathrm{SL}_2$-parabolic bundles (see Remark \ref{rem:PtoH}), which we denote by $\Ps_0(2)$.
\item \emph{Parabolic type \((\{\tau(1)\},\{\tau(2)\})\):} in this case $V=L\oplus L^\vee,$ where the line bundles  \(L\) and \(L^\vee\) carry parabolic weights \(c_{\tau(1)}\) and \(c_{\tau(2)}=-c_{\tau(1)}\), respectively, and  \(\deg L=l\).  The Higgs field  $\Phi$ vanishes on \(L^\vee\), while its restriction to \(L\) gives a nonzero element $$\Phi|_L\in H^0(\mathrm{SParHom}(L,L^\vee)\otimes K(p)) \simeq H^0\!\left(L^{-2}\otimes K(p-S^\tau)\right).$$ 
Such points form a fixed-point locus $F_l^\tau\subset\Hs_0(c)^\TT$.  An analysis of the stability condition, analogous to that in Proposition \ref{prop:range111}, shows that these loci are labelled by integers $l$ satisfying $c_{\tau(1)}<l\leq g-1$, and that \(F_l^\tau\) is a \(2^{2g}\)-fold \'etale cover of  $\Sym^{2g-2l-1-s^\tau\!}(C).$ \end{enumerate}

Thus for either choice $\mathcal{V}\!=\!\mathcal{L}(k;(\lambda_1,\lambda_2))$ or $\mathcal{V}\!=\!\mathcal{L}(k;(\lambda_1,\lambda_2))\otimes\pi_!\bigl(U_{(\nu_1,\nu_2)}\otimes K^{1/2}\bigr)$, the localization formula \eqref{Euler} reads as 
\begin{equation}\label{explicitrk2}
\chi_t(\Hs_0(2),\mathcal{V})=  \int_{\Ps_0(2)}\frac{ch_\TT\bigl(\mathcal{V}|_{\Ps_0(2)}\bigr)}{e_\TT\bigl(N^\vee_{\Ps_0(2)}\bigr)}\mathrm{Td}(\Ps_0(2)) +  \sum_{\tau\in\Sigma_2} \sum^{g-1}_{c_{\tau(1)}<l}\int_{F^\tau_{l}}\frac{ch_\TT\bigl(\mathcal{V}|_{F^\tau_{l}}\bigr)}{e_\TT\bigl(N^\vee_{F^\tau_{l}}\bigr)}\mathrm{Td}(F^\tau_{l}). 
\end{equation}
We now compute the contributions of $\Ps_0(2)$ and of the loci $F_l^\tau$ to \eqref{explicitrk2}.
\vspace{-5pt}
\subsubsection*{The contribution of the component of parabolic type $(\{1,2\})$}
\vspace{-3pt}
To compute the contribution of $\Ps_0(2)$ to \eqref{explicitrk2}, we proceed exactly as in \S\ref{S:3} for the fixed component of parabolic type $(\{1,2,3\})$ in rank $3$. Since the residue formula of \cite{Unyu} used there holds in arbitrary rank, one obtains that this contribution is given, for either choice of $\mathcal V$ above, by the residue at the origin in the corresponding formula \eqref{PHV2Lkl} or \eqref{PHV2nu}.
\vspace{-5pt}

\subsubsection*{The contribution of the fixed-point locus $F^\tau_l$}
\vspace{-3pt}
We next compute the contribution of $F_l^\tau$ to the equivariant Euler characteristic \eqref{explicitrk2}. Let 
$$\eta\in H^2(\Sym^{2g-2l-1-s^\tau}(C)), \;\;\;  \xi_i\in H^1(\Sym^{2g-2l-1-s^\tau}(C)), \;\;\;  i=1,\ldots,2g,$$
be the cohomology classes arising from the K\"unneth decomposition of $c_1(\mathcal{O}(\Delta))$, where $\Delta\subset \Sym^{2g-2l-1-s^\tau}(C)\times C$ is the universal divisor, and set $\theta:=\sum_{i=1}^g\xi_i\xi_{i+g}.$
We use the same notation for the pullbacks of these classes to \(F_l^\tau\). We now express the characteristic classes entering the integral over
$F_l^\tau$ in \eqref{explicitrk2} in terms of $\eta$ and $\theta$. 

 First, repeating the calculations of Propositions \ref{LklonF}, \ref{ElonF} and Lemma \ref{l:eqT} in the present setting, with $t=e^u$, $l_1=l=-l_2,$ $\eta_1=\eta$, $\xi_i^1= \xi_i$, $\theta_1=\theta,$ and all the remaining cohomology classes equal to zero, we obtain 
\begin{multline}\label{eq:rk2class1}
e_\TT(N^\vee_{F_l^\tau})^{-1} = \frac{\exp\bigl(\theta -\tfrac{\theta}{1-e^{2u-\eta}}\bigr)}{(1-e^{-u}\!)^{g-1}(1-e^{u-\eta}\!)^{g-1+s^{\!\tau}+2l}(1-e^{\eta-u}\!)^{g-s^{\tau}-2l}(1-e^{\eta-2u}\!)^{2l+g-1+s^{\!\tau}}} \\
 \text{\; and \;} \ch_{\TT}\bigl(\mathcal{L}(k;\lambda)\big|_{F_l^\tau}\bigr) = \exp\bigl(-\lambda_{\tau(1)}(\eta-u)+kl(\eta-u)+\tfrac{k}{2}\theta\bigr). \;\;\;\;\;\;\;\;\;\;\;\;\;\;\;\;\;\;
\end{multline}

For the second choice of \(\mathcal V\), set $v:=\left(\tfrac{\nu_1+\nu_2}{2},-\tfrac{\nu_1+\nu_2}{2}\right)$ and let  $$
\phi^{\!\nu}(y):=\frac{e^{(\nu_1-\nu_2+1)y/2}-e^{-(\nu_1-\nu_2+1)y/2}}{e^{y/2}-e^{-y/2}} 
$$
denote the character of the irreducible $\mathrm{GL}_2$-representation of highest weight $\nu=(\nu_1,\nu_2)$ restricted to the maximal torus of $\mathrm{SL}_2$. Repeating the calculation of Proposition \ref{UnyuonF} then yields 
\begin{equation}\label{eq:rk2class2}
\ch_\TT\bigl(\pi_!\bigl(U_\nu\otimes K^{1/2}\bigr)|_{F_l^\tau}\bigr)= \exp(-v_{\tau(1)}(\eta-u))
\bigl(-l \dot{\phi}^\nu(\eta-u)-\tfrac{\theta}{2} \ddot{\phi}^\nu(\eta-u)\bigr), 
\end{equation}
where  we denote $\dot{\phi}^\nu(y):=2\partial_y\phi^{\!\nu}(y)$ and $\ddot{\phi}^\nu(y):=2\partial^2_y\phi^{\!\nu}(y)$.

Finally, since \(F_l^\tau\) is a finite \'etale cover of \(\Sym^{2g-2l-1-s^\tau}(C)\), Lemma \ref{l:intsym}(i) gives
\begin{equation}\label{eq:rk2class3}
\mathrm{Td}(F_l^\tau) = \frac{\eta^{g-2l-s^\tau} \exp\bigl(-\tfrac{\theta}{\eta} -\tfrac{\theta}{1-e^\eta} \bigr)}{(1-e^{-\eta})^{g-2l-s^\tau}}.
\end{equation}

Now, to evaluate the contributions of $F_l^\tau$ to the equivariant Euler characteristics \eqref{explicitrk2}, we introduce a formal parameter \(\delta\), and for  $\lhat=(\lhat_1,\lhat_2)\in\mathbb{Q}^2$ with $\lhat_2=-\lhat_1$, consider the integral
\begin{multline}\label{eq:rk2int}
\!\!\underset{\Sym^{2g\!-\!2l\!-\!1\!-\!s^\tau\!}(C)}{\int}\!\! \frac{ 2^{2g}\, \eta^{g-2l-s^{\!\tau}}  e^{-\lhat_{\tau(1)}(\eta-u)+l(k+2)(\eta-u)-\delta l \dot{\phi}^{\!\nu}(\eta-u)}}{W(\eta\!-\!u)^{\!g-1}(1-e^{u-\eta})^{s^{\!\tau}}(1-e^{\eta-u})^{1-s^{\tau}}} \\
\frac{   \exp\bigl( -\tfrac{\theta}{\eta}+\theta \bigl(\tfrac{k}{2}+1 -\tfrac{1}{1-e^{2u-\eta}}-\tfrac{1}{1-e^{\eta}} - \delta\tfrac{1}{2}\ddot{\phi}^{\nu\!}(\eta-u)\bigr) \bigr)}{(1-e^{\eta-2u})^{2l+s^{\!\tau}}(1-e^{-\eta})^{-2l+1-s^{\!\tau}}   }, 
\end{multline}
where $W(y):=(1-e^{-u})(1-e^y)(1-e^{-y})(1-e^{y-u})(1-e^{-y-u})$. 
Combining \eqref{eq:rk2class1} and \eqref{eq:rk2class3}, we note that the integral \eqref{eq:rk2int}, with \(\widehat\lambda=\lambda\) and \(\delta=0\), is precisely the contribution of the fixed-point locus \(F_l^\tau\) to \eqref{explicitrk2} for $\mathcal V=\mathcal L(k;\lambda). $ 
Using in addition \eqref{eq:rk2class2}, we obtain that the coefficient of \(\delta\) in \eqref{eq:rk2int} with \(\widehat\lambda=\lambda+v\) equals the contribution of \(F_l^\tau\) to \eqref{explicitrk2}  for
$\mathcal V=\mathcal L(k;\lambda)\otimes\pi_!(U_\nu\otimes K^{1/2}).$
It therefore remains to evaluate \eqref{eq:rk2int}. 

Applying Lemma \ref{l:intsym}(ii) and making the change of variable \(y=z-u\), we obtain that \eqref{eq:rk2int} is equal to
\begin{equation}\label{eq:rk2res}
 \underset{y=-u}{\res}\frac{ 2^ge^{\lhat_{\tau(2)} y}\bigl(k+2-\frac{2}{1-e^{u-y}}-\frac{2}{1-e^{u+y}}-\delta\ddot{\phi}^{\!\nu\!}(y)\bigr)^{\!g} \!\left(\!\frac{1-e^{-y-u}}{1-e^{y-u}}\!\right)^{\!2l}\!e^{(k+2)ly-\delta l\dot{\phi}^{\!\nu\!}(y)}  }{W(y)^{\!g-1}(1-e^{-y})^{s^{\!\tau}}(1-e^{y})^{1-s^{\tau}}(1-e^{y-u})^{s^{\!\tau}}(1-e^{-y-u})^{1-s^{\!\tau}}} \,\mathrm{d}y.
\end{equation}

\subsubsection*{The total contribution of the components of parabolic type $(\{\tau(1)\},\{\tau(2)\})$}
For fixed permutation \(\tau\), the total contribution of the fixed-point locus of parabolic type $(\{\tau(1)\},\{\tau(2)\})$ to the equivariant Euler characteristic \eqref{explicitrk2} is obtained by summing \eqref{eq:rk2res} over the integers $l$ satisfying $c_{\tau(1)}<l\leq g-1$. 

In fact, the upper bound $l\leq g-1$ can be omitted. A simple calculation shows that the differential form in $\eqref{eq:rk2res}$ has a pole at $y=-u$ of order at most $2g-2l-s^\tau$. Hence, its residue vanishes for $l\geq g$, so the sum may be extended over all integers $l>c_{\tau(1)}$.

As a result, substituting $s^{\mathrm{id}}=0,$ $s^{(12)}=1$ into \eqref{eq:rk2res} and using the fact that $c_1>0>c_2$ for any $c=(c_1,c_2)\in\Delta(2)$, we find that the total contribution of the fixed-point locus of parabolic type $(\{2\},\{1\})$ to the equivariant Euler characteristics \eqref{explicitrk2} is given by the relevant coefficients in the \(\delta\)-expansion of
$$ \sum_{l\geq 0}\underset{y=-u}{\res}\frac{2^ge^{\lhat_1 y}\bigl(k+2-\frac{2}{1-e^{u-y}}-\frac{2}{1-e^{u+y}}-\delta\ddot{\phi}^{\!\nu\!}(y)\bigr)^{\!g} }{W(y)^{g-1}(1-e^{-y})(1-e^{y-u})}\left(\!\left(\!\frac{1-e^{-y-u}}{1-e^{y-u}}\!\right)^{\!2}\!e^{(k+2)y-\delta\dot{\phi}^{\!\nu\!}(y)}\!\right)^{\!l} \mathrm{d}y,$$
evaluated at the appropriate value of $\lhat$, 
while the total contribution of the fixed-point locus of parabolic type $(\{1\},\{2\})$ to \eqref{explicitrk2} is given by the $\delta$-coefficients of 
$$ \sum_{l\geq 1}\underset{y=-u}{\res}\frac{2^ge^{-\lhat_1 y}\bigl(k+2-\frac{2}{1-e^{u-y}}-\frac{2}{1-e^{u+y}}-\delta\ddot{\phi}^{\!\nu\!}(y)\bigr)^{\!g} }{W(y)^{g-1}(1-e^{y})(1-e^{-y-u})}\left(\!\left(\!\frac{1-e^{-y-u}}{1-e^{y-u}}\!\right)^{\!2}\!e^{(k+2)y-\delta\dot{\phi}^{\!\nu\!}(y)}\!\right)^{\!l} \mathrm{d}y$$ with  the appropriate value of $\lhat$. 
Changing the variable $y$ to $-y$ in the second expression and adding the contribution of $\Ps_0(2)$, we arrive at the following result.
\begin{theorem}\label{Rank2Result}
With the notation introduced  above, for a formal parameter $\delta$ and $\widehat\lambda\in\mathbb Q$, define
\begin{equation*}
\Omega_\delta(\lhat):= \frac{2^{g}\, e^{\lhat y}\, \bigl(k+2-\frac{2}{1-e^{u+y}}-\frac{2}{1-e^{u-y}}-\delta\ddot{\phi}^{\!\nu\!}(y)\bigr)^{g} \, \mathrm{d}y}{W(y)^{g-1}(1-e^{-y})(1-e^{y-u})\big(1-e^{(k+2)y-\delta\dot{\phi}^{\!\nu\!}(y)}\big(\frac{1-e^{-y-u}}{1-e^{y-u}}\big)^2\big)}. 
\end{equation*}
\begin{enumerate}[label=(\roman*), leftmargin=1.5pc]
\item Then, with $t=e^u$, we have 
\begin{equation}\label{PHV2Lkl}
\chi_t(\Hs_0(2), \mathcal{L}(k;(\lambda,-\lambda))) = \underset{y=0}{\mathrm{Res}}\,\Omega_0(\lambda)+ \underset{y=u}{\mathrm{Res}}{}^-\Omega_0(\lambda) + \underset{y=-u}{\mathrm{Res}}\!\!{}^+\Omega_0(\lambda); 
\end{equation}
\vspace{-17pt}
\begin{multline}\label{PHV2nu}
\chi_t(\Hs_0(2), \mathcal{L}(k;(\lambda,-\lambda))\otimes\pi_!(U_{(\nu_1,\nu_2)}\otimes K^{1/2})) =  \\ \frac{\mathrm{d}}{\mathrm{d} \delta}\left.\left(\underset{y=0}{\mathrm{Res}}\,\Omega_\delta\bigl(\lambda+\tfrac{\nu_1+\nu_2}{2}\bigr) +\underset{y=u}{\mathrm{Res}}{}^-\Omega_\delta\bigl(\lambda+\tfrac{\nu_1+\nu_2}{2}\bigr) + \underset{y=-u}{\mathrm{Res}}\!\!{}^+\Omega_\delta\bigl(\lambda+\tfrac{\nu_1+\nu_2}{2}\bigr)\right)\right|_{\delta=0}. \end{multline}
Here the superscripts \(+\) and \(-\) on the residue symbols specify the choice of geometric-series expansion of the last factor in the denominator of \(\Omega_\delta(\lhat)\): \(+\) denotes the expansion in nonnegative powers, while \(-\) denotes the expansion in negative powers.

\item Moreover, the first, second and third summands in \eqref{PHV2Lkl} and \eqref{PHV2nu} give the total contributions to the corresponding equivariant Euler characteristic from the fixed-point loci of parabolic types $(\{1,2\})$, $(\{1\}, \{2\})$ and $(\{2\},\{1\})$, respectively.
\end{enumerate}
\end{theorem}

\section{Residue formulas in rank $3$}\label{S:mainrank3}
We now turn to the moduli spaces $\Hs_0(c)$ of rank-$3$ stable parabolic Higgs bundles with trivial determinant and parabolic weights $c=(c_1,c_2,c_3)\in\Delta:=\Delta(3)$. A simple calculation shows that $\Delta$ is the open triangle with vertices $(0,0,0)$, $(1/3,1/3,-2/3)$ and $(2/3,-1/3,-1/3)$, and that the stability condition depends only on the sign of $c_2$: all weights with $c_2>0$ are regular and define the same stability condition, and likewise for $c_2<0$. Thus there are only two essentially different moduli spaces; cf. \cite{OTASz}.

By Proposition \ref{fixedcomp}, every \(\TT\)-fixed point \((V,\Phi)\in\Hs_0(c)^\TT\) is of one of the following types:
\begin{enumerate}[label=(\alph*), leftmargin=1.5pc]
\item \emph{type \((3)\):} \(V\) is a rank-$3$ stable parabolic bundle and \(\Phi=0\);
\item \emph{type \((1,1,1)\):} \(V\simeq V_1\oplus V_2\oplus V_3\), where \(V_1, V_2\) and $V_3$ are parabolic line bundles and \(\Phi\) restricts to nonzero strongly parabolic morphisms $V_1\to V_2\otimes K(p)$ and $V_2\to V_3\otimes K(p)$;
\item \emph{type \((2,1)\):} \(V\simeq V_1\oplus V_2\), where \(V_1\) is a rank-$2$ parabolic bundle, \(V_2\) is a parabolic line bundle, and $\Phi$ restricts to a nonzero strongly parabolic morphism $V_1\to V_2\otimes K(p)$;
\item \emph{type \((1,2)\):} \(V\simeq V_1\oplus V_2\), where \(V_1\) is a parabolic line bundle, \(V_2\) is a rank-$2$ parabolic bundle, and
$\Phi$ restricts to a nonzero strongly parabolic morphism $V_1\to V_2\otimes K(p)$.
\end{enumerate}

The fixed points of type $(3)$ form the connected component $\Ps_0(c)$ (see Remark \ref{rem:PtoH}).  For the remaining types, we will show that the components of type $(1,1,1)$ are finite covers of products of symmetric powers of \(C\), while those of types \((2,1)\) and \((1,2)\) are isomorphic to moduli spaces of parabolic triples, introduced in \S\ref{S:partriplesgen}. We compute the contributions of these components to the equivariant Euler characteristics \eqref{Euler} in \S\ref{S:3}, \S\ref{S:111}, and \S\ref{S:21and12}, respectively.

Before stating our main result, Theorem \ref{thm:main}, which gives explicit formulas for \eqref{Euler}  for every regular $c\in\Delta$, we introduce some additional notation in the next subsection.

\subsection{Hyperplane arrangements and $\iber$ operator}\label{SS:iberdef}
We denote by  $\mathfrak{t}=\mathbb{R}^3/\mathbb{R}(1,1,1) $ the Cartan subalgebra of the Lie algebra of $SU_3$ and identify its dual with
$$\mathfrak{t}^*=\{a=(a_1,a_2,a_3)\in\mathbb{R}^3|\; a_1+a_2+a_3=0\}.$$ Let $ x_1, x_2, x_3 $ be the coordinates on $ \mathbb{R}^3$; for $ a\in \mathfrak{t}^* $, we write $\langle a,x \rangle = \sum_ia_ix_i $ for the induced linear function on $\mathfrak t$, and sometimes we denote this function simply by $a$.
Given a dominant weight $\nu=(\nu_1, \nu_2, \nu_3)\in\mathbb{Z}^3$ of $\mathrm{GL}_3$, let $m_\mu$ denote the multiplicity of the weight $\mu\in\mathbb{Z}^3$ in the
representation $\rho_\nu$, and put $|\nu|=\nu_1+\nu_2+\nu_3$. The weights of the restriction
$\bar\rho_\nu$ of $\rho_\nu$ to $SU_3$ are $\bar\mu=\mu-\frac{|\nu|}{3}(1,1,1)\in\mathfrak{t}^*$, and the character of $\bar\rho_\nu$ is the holomorphic function on $\mathfrak t\otimes\mathbb C$ given by  
\begin{equation}\label{phinumu}
  \phi^\nu(x):=\mathrm{trace}\,(\bar\rho_\nu\circ\exp(x))   =\sum_\mu m_\mu\, e^{\langle\bar\mu,x\rangle}.
\end{equation}

We now fix a positive real parameter $u$ and consider  the \textit{affine hyperplane arrangement} 
\begin{equation}\label{HPAu}
\A := \bigcup_{1\le i<j\le 3}\{ x\in \mathfrak{t}: \langle \alpha_{i j},x\rangle := x_i-x_j \in \{ 0, u, -u \}  \}.
\end{equation}
We define the set of \textit{poles} in $\A$, denoted by $\mathrm{Poles}(\A)$, to be the set of zero-dimensional intersections of hyperplanes in $\A$. 
\begin{proposition}\label{poles-components}
There is a one-to-one correspondence between the poles of \(\A\) and the ordered partitions of $\{1,2,3\}$. 
\end{proposition}
\begin{proof}
Let $p=(p_1,p_2,p_3)\in\mathrm{Poles}(\A)$. Since $p$ is the intersection of two hyperplanes of $\A$, its distinct coordinate values can be ordered uniquely as $q_1,\ldots,q_m,$ $q_a-q_{a+1}=u.$ Setting $\mu_a=\{ i:p_i=q_a\}$ then defines an ordered partition $\mu[p]=(\mu_1,\ldots,\mu_m)$ of $\{1,2,3\}$.

Conversely, an ordered partition determines a unique pole by assigning equal coordinates to indices in the same group, requiring successive  values to differ by $u$, and imposing $\sum_i p_i=0$. Thus $p\mapsto\mu[p]$ is a bijection.
\end{proof}
Putting together Propositions \ref{poles-components} and \ref{fixedcomp}, we obtain the following statement.
\begin{corollary}\label{cor:identpmup}
There is a one-to-one correspondence between the poles of the hyperplane arrangement  $\A$ and the parabolic types of  connected components of the fixed-point locus $\Hs_0(c)^\TT$.
\end{corollary}

For $p\in\mathrm{Poles}(\A)$, let $\A_p$ denote the \textit{central arrangement} obtained by translating by $-p$ the hyperplanes of $\A$ passing through $p$. Write $\A_p = \{H_i, \, i\in I\} $ for its distinct hyperplanes, indexed by a finite set $I$, and let $$\Lambda=\{\lambda=(\lambda_1,\lambda_2,\lambda_3)\in\mathbb{Z}^3|\; \lambda_1+\lambda_2+\lambda_3=0\}\subset \mathfrak{t}^*$$  be the integer lattice. For each $i\in I$, choose a primitive form $\bta{i}\in \Lambda$, unique up to sign, such that $H_i = \{x\in \mathfrak{t}: \langle \bta{i}, x \rangle =0\}$ and set $\widetilde{\A}_p:=\{\pm\bta{i}\mid i\in I\}.$ We denote by  \[ \Bases(\A_p)= \left\{\mathbf{B}=\left(\bta1,\bta{2}\right) \in\widetilde{\A}_p^{2}|\; \bb \text{ -- basis of }  \mathfrak{t}^*\right\} \] the set of \textit{ordered bases}   of $\mathfrak{t}^*$ consisting of elements of $\widetilde{\A}_p$.

Finally, using the standard Euclidean product to identify $\mathfrak{t}^* \simeq \mathfrak{t}$,  we denote by ${\check{\beta}}$ the image of $\beta\in \mathfrak{t}^*$ under this isomorphism. For a Weyl-invariant function   $Q$ on $\mathfrak{t}$,  we write $Q_{\check{\beta}}$ for the directional derivative of \(Q\) in the direction \(\check\beta\).

With this notation in place, we can now define the iterated residue operator that will be used to state our main result.

\begin{definition}\label{def:ires}
Given a pole $p\in\mathrm{Poles}(\A)$, a basis $\bb=(\bta1,\bta2)\in\Bases(\A_p)$, a meromorphic function $f(x)$ defined in a neighborhood of $p$ in $\mathfrak{t}\otimes\mathbb{C}$, with poles contained in the hyperplanes of $\A$, an element $a\in \Lambda$, and $\pmb{\varepsilon}=(\varepsilon_1,\varepsilon_2)\in\{+,-\}^2$, let $Q(x)$ be a Weyl-invariant function, locally holomorphic away from the hyperplanes of $\A$, such
that $\exp(Q_{\check\beta})$ is meromorphic near $p$ for every $\beta\in\Lambda$ and only finitely many terms in the expansions prescribed by
$\pmb{\varepsilon}$ below have nonzero iterated residue.
Introduce coordinates $y_j=\langle\bta{j},x\rangle,$  $j=1,2,$
and write
\[
\widehat f(y):=f(x),\;\;\; \widehat Q_{\check\beta}(y):=Q_{\check\beta}(x), \text{ \; for \;}  \beta\in\Lambda.\]
Set $y_j(p):=\langle\bta{j},p\rangle$, and define
\[
\iber_{p,\bb,Q}\!{}^{\pmb{\varepsilon}}[f](a) := \underset{y_1=y_1(p)}{\mathrm{Res}^{\varepsilon_1}}\, \underset{y_2=y_2(p)}{\mathrm{Res}^{\varepsilon_2}}\,
\frac{ \widehat f(y)\exp\bigl(\widehat Q_{\check a}(y)\bigr)\, \mathrm d\widehat Q_{\check\beta^{[1]}} \wedge\mathrm d\widehat Q_{\check\beta^{[2]}}}{\bigl(1-\exp(\widehat Q_{\check\beta^{[1]}})\bigr)\bigl(1-\exp(\widehat Q_{\check\beta^{[2]}})\bigr)}.
\]
The residues are taken successively in the order $y_2$, $y_1$, so that $y_1$ is treated as a parameter when taking the residue in $y_2$. The superscript $+$ on the $y_j$-residue means that the corresponding factor in the denominator is expanded in the positive direction, $(1-\exp(\widehat Q_{\check\beta^{[j]}}))^{-1}
=\sum_{n\ge0}\exp(n\widehat Q_{\check\beta^{[j]}}),$
whereas the superscript $-$ means that it is expanded in the negative direction, $(1-\exp(\widehat Q_{\check\beta^{[j]}}))^{-1} =-\sum_{n\ge1}\exp(-n\widehat Q_{\check\beta^{[j]}}).$

When no expansion direction is prescribed, we instead regard the denominator factors as meromorphic functions at $p$ and write
$$\iber_{p,\bb,Q}[f](a):=\underset{y_1=y_1(p)}{\mathrm{Res}}\, \underset{y_2=y_2(p)}{\mathrm{Res}}\, \frac{\widehat f(y)\exp\bigl(\widehat Q_{\check a}(y)\bigr)\mathrm d\widehat Q_{\check\beta^{[1]}}\wedge\mathrm d\widehat Q_{\check\beta^{[2]}}}{\bigl(1-\exp(\widehat Q_{\check\beta^{[1]}})\bigr)\bigl(1-\exp(\widehat Q_{\check\beta^{[2]}})\bigr)},$$
where the residues are taken successively in the order $y_2$, $y_1$.
\end{definition}

We now make a simple observation concerning the choice of expansion direction, which will be used later.  Given $\bb=\left(\bta1,\bta{2}\right) \in \Bases(\A_p)$ and  $a\in \mathfrak{t}^*$, write $$a = [a]_{\bb}+ \{a\}_{\bb}, \;\;  [a]_\bb\in\Lambda, \;\;  \{a\}_{\bb}= s_1\bta1+ s_2\bta2, \; s_1,s_2\in[0,1).$$

\begin{lemma}\label{rem:varChange}
Suppose \(s_j\neq0\). Let \(\bb'\) be obtained from \(\bb\) by replacing \(\beta^{[j]}\) with \(-\beta^{[j]}\), and let  \(\pmb\varepsilon'\) be obtained from \(\pmb\varepsilon\) by reversing its \(j\)-th sign. Then $[a]_{\mathbf B'}=[a]_{\mathbf B}+ \beta^{[j]}$, which implies
$\underset{p,\mathbf B',Q}{\iber}\!{}^{\pmb\varepsilon\!'}[f](-[a]_{\mathbf B'})=\underset{p,\mathbf B,Q}{\iber}\!{}^{\pmb\varepsilon}[f](-[a]_{\mathbf B}).$
\end{lemma}

\subsection{Residue formulas in rank $3$}\label{SS:Main}

Equipped with operator $\iber$,  we can now state the formulas for the equivariant Euler characteristics \eqref{Euler} on the moduli spaces of  stable strongly parabolic $\mathrm{SL}_3$-Higgs bundles.
\begin{theorem}\label{thm:main}
Let $c\in \Delta$ be regular, $k\in \mathbb{Z}_{>0}$, $\lambda\in\Lambda$, and for a dominant weight $\nu\in\mathbb{Z}^3$ of $\mathrm{GL}_3$ set $v^{\det}=\tfrac{\sum_i\nu_i}{3}(1,1,-2)$. Let $\delta$ be a formal parameter and define 
$$Q_\delta(x) := \frac{k+3}{6}\sum_{ i<j}\langle \alpha_{i j},x\rangle^2+\sum_{i\neq j}Li_2(e^{\langle \alpha_{i j},x \rangle-u}) - \delta\phi^\nu(x)$$
and $W_t(x) := (1-t^{-1})^{2}\prod_{i\neq j}((1-t^{-1}e^{\langle\alpha_{i j},x\rangle})(1-e^{\langle\alpha_{i j},x\rangle}))$
\begin{enumerate}[label=(\roman*), leftmargin=1.5pc]
\item Then, with $t=e^u$, we have
\begin{multline}\label{conjectureeq}
 \chi_{t}(\Hs_0(c), \mathcal{L}(k;\lambda))= \\
 \! \sum_{p\in\mathrm{Poles(\A)}}\sum_{\bb\in\Bases[p]} \iber_{p, \mathbf{B},Q_0}\!\!\!{}^{\pmb{\varepsilon}[p]}  \left[ \frac{3^g\cdot \det(\mathrm{Hess}_\bb(Q_0))^{g-1} \cdot e^{\langle \lambda, x\rangle}}{W_{\!e^u\!}(x)^{g-1}\prod_{i<j}(1-e^{\langle \alpha_{j i}, x \rangle})(1-e^{\langle \alpha_{i j}, x \rangle-u})} \right](-[c]_{\mathbf{B}}), \end{multline} 
 \begin{multline}\label{conjectureeqnu}
 \chi_{t}(\Hs_0(c), \mathcal{L}(k;\lambda)\otimes\pi_!(U_\nu\otimes K^{1/2}))=  \\
 \sum_{p\in\mathrm{Poles(\A)}}\sum_{\bb\in\Bases[p]} \frac{\mathrm{d}}{\mathrm{d} \delta} \left.\left(  \iber_{p, \mathbf{B},Q_\delta}\!\!\!{}^{\pmb{\varepsilon}[p]}    \left[ \frac{3^g\cdot \det(\mathrm{Hess}_\bb(Q_\delta))^{g-1} \cdot e^{\langle \lambda+v^{\det}, x\rangle}}{W_{\!e^u\!}(x)^{g-1}\prod_{i<j}(1-e^{\langle \alpha_{j i}, x \rangle})(1-e^{\langle \alpha_{i j}, x \rangle-u})} \right](-[c]_{\mathbf{B}})\right) \right|_{\delta=0}\!\!, \end{multline} 
where, for each pole $p\in\mathrm{Poles}(\A)$,  the set $\Bases[p]\subset \Bases(\A_p)$ is defined in Table \ref{table1}, $\pmb{\varepsilon}[p]=(+,+)$ for $p\neq0$, while for $p=0$ we use the ordinary iterated residue and omit the superscript; 
and, for $\bb=(\bta1,\bta2)$, $\mathrm{Hess}_\bb(Q)$ denotes the matrix $(\partial_{\bta{i}}Q_{\check{\beta}^{[j]}})_{i,j=1,2}$, with $\partial_{\bta{i}}:=\tfrac{\partial}{\partial y_i}$ in the
coordinates $y_i=\langle\beta^{[i]},x\rangle$.
\item 
The left-hand sides of \eqref{conjectureeq} and \eqref{conjectureeqnu}  are sums, over parabolic types, of the total contributions of all fixed-point components of that type, while the right-hand sides are sums over poles of $\A$. By Corollary \ref{cor:identpmup}, these indexing sets are naturally identified. Under this identification, the two sums agree term by term.
\end{enumerate}
\end{theorem}
\vspace{-17pt}
\begin{table}[H] 
\centering
\begin{tabular}{ |c|c|c| }
\hline
$p\in \mathrm{Poles(\A)}$ & $\Bases[p] \subset \Bases(\A_p)$ &  \textit{parabolic type} $\mu[p]$ \\
\hline
$(0,0,0)$  & $(\alpha_{2 3}, \alpha_{1 2}), (\alpha_{3 2}, \alpha_{1 3})$ &   $\{1,2,3\}$ \\
\hline 
$(-u,0,u)$& $(\alpha_{2 3}, \alpha_{1 2})$  & \{3\}, \{2\}, \{1\} \\ 
\hline 
$(u,0,-u)$ & $(\alpha_{2 1}, \alpha_{3 2})$  &  \{1\}, \{2\}, \{3\} \\ 
\hline 
$(u,-u,0)$& $(\alpha_{3 1}, \alpha_{2 3})$ &  \{1\}, \{3\}, \{2\} \\ 
\hline 
$(-u,u,0)$& $(\alpha_{3 2}, \alpha_{1 3})$ & \{2\}, \{3\}, \{1\} \\ 
\hline 
$(0,-u,u)$ & $(\alpha_{1 3}, \alpha_{2 1})$ &  \{3\}, \{1\}, \{2\} \\ 
\hline 
$(0,u,-u)$& $(\alpha_{1 2}, \alpha_{3 1})$ &  \{2\}, \{1\}, \{3\} \\ 
\hline 
$(u/3,u/3,-2u/3)$ & $(\alpha_{2 1}, \alpha_{3 1}), (\alpha_{1 2}, \alpha_{3 2})$   &  \{1, 2\}, \{3\} \\ 
\hline 
$(u/3,-2u/3,u/3)$ & $(\alpha_{3 1}, \alpha_{2 1}), (\alpha_{1 3}, \alpha_{2 3})$  & \{1, 3\}, \{2\} \\ 
\hline 
$(-2u/3,u/3,u/3)$ & $(\alpha_{3 2}, \alpha_{1 2}), (\alpha_{2 3}, \alpha_{1 3})$  & \{2, 3\}, \{1\} \\ 
\hline
$(-u/3,-u/3,2u/3)$ & $(\alpha_{1 2}, \alpha_{1 3}), (\alpha_{2 1}, \alpha_{2 3})$ &  \{3\}, \{1, 2\} \\
\hline
$ (-u/3,2u/3,-u/3)$ & $(\alpha_{1 3}, \alpha_{1 2}), (\alpha_{3 1}, \alpha_{3 2})$ &  \{2\}, \{1, 3\} \\
\hline
$ (2u/3,-u/3,-u/3)$ & $(\alpha_{2 3}, \alpha_{2 1}), (\alpha_{3 2}, \alpha_{3 1})$ &   \{1\}, \{2, 3\}\\
\hline
\end{tabular}
\vspace{2pt}
\setlength{\belowcaptionskip}{-7pt}\caption{\textit{Poles, the corresponding bases, and parabolic types of the fixed-point components.}}\label{table1}
\end{table}
\begin{remark}
The residue formulas \eqref{conjectureeqnu} extend to the equivariant Euler characteristics of bundles of the form
$$\mathcal L(k;\lambda)\otimes \bigotimes_{j=1}^m\pi_!\bigl(U_{\nu[j]}\otimes K^{1/2}\bigr),$$
where \(\nu[1],\ldots,\nu[m]\) are dominant weights of \(\mathrm{GL}_3\). As in \cite[Theorem 6.1]{Unyu}, one replaces the single-parameter deformation of \(Q\) by the multi-parameter deformation 
$$\mathbf{Q}_{\boldsymbol\delta}(x) := \frac{k+3}{6}\sum_{ i<j}\langle \alpha_{i j},x\rangle^2+\sum_{i\neq j}Li_2(e^{\langle \alpha_{i j},x \rangle-u}) - \sum_{j=1}^m\delta_j\phi^{\nu[j]}(x),$$
and replaces the determinant-shift $v^{\det}$ by
$\mathbf{v}^{\det}:=\frac{1}{3}(1,1,-2)\sum_{j=1}^m|\nu[j]|.$
The corresponding formulas are then obtained by taking derivatives with respect to \(\delta_1,\ldots,\delta_m\) and evaluating at \(\delta_1=\cdots=\delta_m=0\).
\end{remark}
We now rewrite the directional-derivative and Hessian  factors appearing in the formulas of Theorem \ref{thm:main} in a form convenient for our later calculations. 
\begin{lemma}\label{lem:phidef}
A simple calculation shows that for any $\beta\in\Lambda$
$$\exp({Q_{\delta}}_{\check \beta}(x)) = e^{(k+3)\langle\beta,x\rangle}\prod_{i<j}\phi(\langle\alpha_{i j},x\rangle)^{\langle\alpha_{i j},\check{\beta}\rangle}\exp(-\delta \phi^\nu_{\check{\beta}}(x)), \text{\; where \;} \phi(z):=\tfrac{1-e^{-z-u}}{1-e^{z-u}}.$$ 
\end{lemma}
To describe the Hessian term, set $z_1=\langle\alpha_{21},x\rangle$, $z_2=\langle\alpha_{32},x\rangle$ and 
write  $\widehat\phi^\nu(z_1,z_2):=\phi^\nu(x)$. Equation \eqref{phinumu} then gives
$$\widehat\phi^\nu(z_1,z_2)=\sum_\mu m_\mu\,e^{-\bar\mu_1z_1+\bar\mu_3z_2}, \text{\; where \;} \bar{\mu}_j = \mu_j - \tfrac13 (\nu_1+\nu_2+\nu_3).$$ 
In fact, since the character $\phi^\nu$ is Weyl invariant, for every $\sigma\in\Sigma_3$ we have  $\phi^\nu(x)=\widehat\phi^\nu(z_1,z_2)$ with $z_1=\langle\alpha_{\sigma(2)\sigma(1)},x\rangle$ and $ z_2=\langle\alpha_{\sigma(3)\sigma(2)},x\rangle.$

We now use these coordinates to rewrite the determinant of the Hessian of $Q_\delta(x)$. 
\begin{lemma}\label{lem:Gdef}
For a permutation $\sigma\in\Sigma_3$, set $\bb_\sigma= (\alpha_{\sigma(2)\sigma(1)}, \alpha_{\sigma(3)\sigma(2)})$.  Using the coordinates $z_1,z_2$ determined by $\bb_\sigma$ as above, define
\begin{multline*}
G(z):= \tfrac{k+3}{3}-\tfrac{1}{1-e^{u+z}}-\tfrac{1}{1-e^{u-z}}, \;\;\; G_{12}(z_1,z_2):= G(z_1+z_2)-\delta\partial_{z_1}\partial_{z_2}\widehat{\phi}^\nu(z_1,z_2) \\
\text{\; and  \;}  G_{j}(z_1,z_2):= G(z_j)-\delta \partial_{z_j}^2\widehat{\phi}^\nu(z_1,z_2)+\delta \partial_{z_1}\partial_{z_2}\widehat{\phi}^\nu(z_1,z_2) \text{\; for \; } j=1,2.
\end{multline*}
A simple calculation then gives
$$\det(\mathrm{Hess}_{\bb_\sigma}(Q_\delta)(x)) = 3(G_1(z_1,z_2)G_2(z_1,z_2) + G_1(z_1,z_2)G_{12}(z_1,z_2)+ G_2(z_1,z_2)G_{12}(z_1,z_2)).$$
\end{lemma}

\section{Contribution of type-$(3)$ components}\label{S:3}
In this section, we compute the contribution of the fixed-point locus of type $(3)$, which is  identified with the moduli space $\Ps_0(c)$,  to the equivariant Euler characteristic \eqref{Euler}. We begin with the following observation.
\begin{proposition}\label{prop:TPs}
The normal bundle to $\Ps_0(c)\subset \Hs_0(c)$ is equivariantly isomorphic to the cotangent bundle $T^{\vee\!}\Ps_0(c)$ with $\TT$-weight $1$.
\end{proposition}
\begin{proof}
Let $V\in\Ps_0(c)$. It follows from  Theorem \ref{thm:modsp}(iii)  that the tangent space to $\Hs_0(c)$ at $(V,0)$ is 
\[ T_{(V,0)}\Hs_0(c)= H^1(\ParEnd_0(V)) \oplus H^0(\SParEnd_0(V)\otimes K(p)),\]
where the first summand has $\TT$-weight zero and hence is identified with the tangent space $T_{V}\Ps_0(c)$, while the second summand has $\TT$-weight $1$. Applying Serre duality, we obtain  \[ H^0(\SParEnd_0(V)\otimes K(p)) \simeq H^1(\ParEnd_0(V))^\vee =T_V^{\vee\!}\Ps_0(c)\]
(see \cite[Proposition 2.3]{GPGM}),  which finishes the proof.
\end{proof}
It follows that, for either choice $\mathcal{V}= \mathcal{L}(k;\lambda)$ or $\mathcal{V}=\mathcal{L}(k;\lambda)\otimes \pi_!(U_\nu\otimes K^{1/2})$, the contribution of the fixed-point locus of type $(3)$ to the equivariant Euler characteristic  \eqref{Euler} is given by 
\begin{equation}\label{eq:type3-first} 
 \sum_{n\geq 0}t^{-n} \chi_t\!\left( \Ps_0(c), \mathcal V\big|_{\Ps_0(c)}\otimes \Sym^n(T\Ps_0(c)) \right).
\end{equation}
Our goal is to evaluate this expression using the residue formulas of \cite{Unyu}. We carry out the computation for $\mathcal V=\mathcal L(k;\lambda)\otimes\pi_!(U_\nu\otimes K^{1/2})$; the case $\mathcal V=\mathcal L(k;\lambda)$ follows by the same argument, with the tautological factor omitted.

We begin by expressing the Chern character of $T\Ps_0(c)$ in terms of cohomology classes associated with the universal bundle. By a slight abuse of notation, throughout this section we use the same symbols $\mathcal U$ and $\pi$ for the restrictions to $\Ps_0(c)\times C$ of the universal bundle on $\Hs_0(c)\times C$ and the projection $ \Hs_0(c)\times C\to\Hs_0(c),$ respectively; see Theorem \ref{thm:modsp}(iv). We set $\mathcal F_0=0$ and, for $1\leq i\neq j\leq3$, write
$$ z_i-z_j:=c_1\!\left(\mathcal F_{4-i}/\mathcal F_{3-i} \otimes \left(\mathcal F_{4-j}/\mathcal F_{3-j}\right)^\vee \right)\big|_{\Ps_0(c)\times\{p\}}.$$
We note that these classes are independent of the choice of the normalization of $\mathcal U$, and, denoting by $\omega\in H^2(C)$ the class Poincar\'e dual to a point, we have 
$$ \ch(\End_0(\mathcal U)) = \ch(\ParEnd_0(\mathcal{U})) + \omega \sum_{i<j}e^{z_i-z_j}  \text{\;\; and \;\;}  \omega\,\ch(\End_0(\mathcal U)) = 2\omega+\omega\sum_{i\neq j}e^{z_i-z_j}.$$
Moreover, since $ T\Ps_0(c)=R^1\pi_*\ParEnd_0(\mathcal U)$ (cf. with the proof of Proposition \ref{prop:TPs}),  the Grothendieck--Riemann--Roch theorem yields
\begin{multline}\label{eq:type3-GRR}
 \ch(T\Ps_0(c)) = -\pi_*\!\left( \ch(\End_0(\mathcal U))(1-(g-1)\omega) \right)+\sum_{i<j}e^{z_i-z_j} = \\
 -\ch\bigl(\pi_!\bigl(\End_0(\mathcal U)\otimes K^{1/2}\bigr)\bigr) +2g-2 +(g-1)\sum_{i\neq j}e^{z_i-z_j} +\sum_{i<j}e^{z_i-z_j}.                
\end{multline}

Next, recall that the standard identity for symmetric powers gives
\begin{equation}\label{eq:type3-GRR-1}
\ch\!\left( \sum_{n\geq0}t^{-n}\Sym^n(T\Ps_0(c))\right)= \exp\!\left(\sum_{n\geq1}\frac{t^{-n}}{n}\ch(\psi^nT\Ps_0(c)) \right), 
\end{equation}
where $\psi^n$ denotes the $n$-th Adams operation: under the splitting principle, it multiplies each Chern root of $T\Ps_0(c)$  by $n$. 
Using \eqref{eq:type3-GRR}, we can rewrite the right-hand side of \eqref{eq:type3-GRR-1} as 
\begin{equation*}\label{eq:type3-sym}
\frac{\exp\!\left( -\sum_{n\geq1}\frac{t^{-n}}{n} \ch\left(\psi^n \pi_!\left(\End_0(\mathcal U)\otimes K^{1/2}\right)\right)\right)}{ (1-t^{-1})^{2g-2} \prod_{i\neq j} (1-t^{-1}e^{z_i-z_j})^{g-1}
 \prod_{i<j} (1-t^{-1}e^{z_i-z_j})}.                                         
\end{equation*}
By \cite[(39)]{Unyu}, we have
$$\ch\bigl(\psi^n \pi_!\bigl(\End_0(\mathcal U)\otimes K^{1/2}\bigr)\bigr)=\frac1n \ch\!\bigl(\pi_!\bigl(\psi^n\End_0(\mathcal U)\otimes K^{1/2}\bigr) \bigr), $$
and thus 
\begin{equation*}\label{eq:type3-sym2}
 \ch\!\left( \sum_{n\geq0}t^{-n}\Sym^n(T\Ps_0(c))\!\right)\! =\!
\frac{ \exp\!\left( -\sum_{n\geq1}\frac{t^{-n}}{n^2}\ch\!\left(\pi_!\left(\psi^n\End_0(\mathcal U)\otimes K^{1/2} \right) \right) \right)}{(1-t^{-1})^{2g-2} \prod_{i\neq j} (1-t^{-1}e^{z_i-z_j})^{g-1}\prod_{i<j} (1-t^{-1}e^{z_i-z_j})}.
\end{equation*}

Thus, applying the Hirzebruch--Riemann--Roch theorem to \eqref{eq:type3-first} with $\mathcal V=\mathcal L(k;\lambda)\otimes\pi_!(U_\nu\otimes K^{1/2})$, we obtain the following expression for the contribution of the fixed-point locus of type $(3)$ to the equivariant Euler characteristic  \eqref{Euler}:
\begin{equation}\label{eq:inttype3}
 \int_{\Ps_0(c)}\!\!\frac{\!\ch\!\left(\!\mathcal L(k;\lambda)\otimes\pi_!(U_\nu\otimes K^{1/2})\big|_{\Ps_0(c)}\right)\exp\!\left(\!-\!\sum_{n\geq1}\frac{t^{-n}}{n^2}\ch( \pi_!(
 \psi^n\End_0(U)\otimes K^{1/2})\!)\!\right)}
 { (1-t^{-1})^{2g-2}\prod_{i\neq j}(1-t^{-1}e^{z_i-z_j})^{g-1}\prod_{i<j} (1-t^{-1}e^{z_i-z_j}) }\mathrm{Td}(\Ps_0(c)).
\end{equation}

We now evaluate this integral by expanding the exponential and applying \cite[Theorem 6.1]{Unyu} coefficientwise in $t^{-1}$. Set $ \phi_{\mathrm{ad}}(x) = 2+\sum_{i\neq j}e^{x_i-x_j}$, and, introducing parameters $\delta_n$, for $n\geq 1$, define
$$\mathbf Q_\delta(x) :=  \frac{k+3}{6}\sum_{i<j}(x_i-x_j)^2 -\delta\phi^\nu(x) - \sum_{n\geq1}\delta_n\phi_{\mathrm{ad}}(nx).$$
Then by \cite[Theorem 6.1 and Lemma 4.6]{Unyu}, \eqref{eq:inttype3} is equal to
\begin{equation}\label{eq:type3-delta}
 \left.
 \frac{\mathrm{d}}{\mathrm{d} \delta}\exp\!\left(\!-\sum_{n\geq1}\!\frac{t^{-n}}{n^2} \! \frac{\partial}{\partial \delta_n}\!\right)\!
 \sum_{\bb\in\mathcal B[0]}  \underset{0,\bb,\mathbf Q_\delta}{\ires}
 \left[ \!\frac{ 3^g \det(\mathrm{Hess}_\bb(\mathbf Q_\delta))^{g-1}  e^{\langle\lambda+v^{\det},x\rangle} }
 {W_{t\!}(x)^{g-1}\prod_{i<j}  (1\!-\!e^{\langle\alpha_{ji},x\rangle})(1\!-\!t^{-1}e^{\langle\alpha_{ij},x\rangle})}\! \right]\! (-[c]_{\bb})\right|_{\substack{\boldsymbol\delta=0 \\ \delta=0}}\!\!,
\end{equation}
where $\boldsymbol\delta=0$ means that $\delta_n=0$ for every $n\geq1$, $\mathcal{B}[0]:=\{(\alpha_{2 3}, \alpha_{ 1 2}), (\alpha_{3 2}, \alpha_{ 1 3}) \},$\footnote{More precisely, the sum in \eqref{eq:type3-delta} may be taken over the elements $\bb$ of any \emph{diagonal basis} of the hyperplane arrangement $\A_0$; see \cite[\S3.2]{OTASz} for the definition. For simplicity, we restrict to the particular choice of $\Bases[0]$.} and $W_t(x)$, $v^{\det}$ are as defined in Theorem \ref{thm:main}. A simple calculation shows that  \eqref{eq:type3-delta}  is equal to
\begin{equation}\label{eq:type3-nodelta}
 \left. \frac{\mathrm{d}}{\mathrm{d} \delta} \,\sum_{\bb\in\mathcal B[0]}\underset{0,\bb,\widetilde Q_\delta}{\ires} \left[\frac{3^g\det(\mathrm{Hess}_\bb(\widetilde Q_\delta))^{g-1} e^{\langle\lambda+v^{\det},x\rangle}}
 {W_{t\!}(x)^{g-1}\prod_{i<j}(1-e^{\langle\alpha_{ji},x\rangle})(1-t^{-1}e^{\langle\alpha_{ij},x\rangle}) }\right] (-[c]_{\bb})\right|_{\delta=0},
\end{equation}
where $ \widetilde Q_\delta(x)$ is obtained  from $\mathbf{Q}_\delta(x)$ by substituting $\delta_n=-\frac{t^{-n}}{n^2}$:
\begin{equation}\label{eq:Qsubst}
\widetilde{Q}_\delta(x) =  \frac{k+3}{6}\sum_{i<j}(x_i-x_j)^2 + \sum_{i\neq j} \mathrm{Li}_2\!\left(t^{-1}e^{x_i-x_j}\right) -\delta \phi^\nu(x)+2\sum_{n\geq1}\frac{t^{-n}}{n^2}.
\end{equation}
Finally, we note that the last term in \eqref{eq:Qsubst} is independent of $x$ and hence does not affect any of the directional derivatives entering the iterated residue \eqref{eq:type3-nodelta};  therefore $\widetilde Q_\delta(x)$ may be replaced by the function $Q_\delta(x)$ defined in Theorem \ref{thm:main}. 

As noted above, the contribution to the equivariant Euler characteristic of $\mathcal{L}(k;\lambda)$ is obtained by the same calculation after omitting the factor $\ch(\pi_!(U_\nu\otimes K^{1/2}))$ in \eqref{eq:inttype3}. This amounts to setting $\nu=0=\delta$ and omitting the operator $\frac{\mathrm{d}}{\mathrm{d}\delta}\big|_{\delta=0}$ in \eqref{eq:type3-delta}.  
We have thus proved the following result.

\begin{theorem}\label{thm:Cont3}
Let $\mathcal{B}[0]=\{(\alpha_{2 3}, \alpha_{ 1 2}), (\alpha_{3 2}, \alpha_{ 1 3}) \}\subset\Bases(\A_{0})$. In the notation of Theorem \ref{thm:main}, the contribution of the fixed-point component of parabolic type $(\{1, 2, 3\})$ to the equivariant Euler characteristic of \(\mathcal L(k;\lambda)\) on \(\Hs_0(c)\) is
$$\sum_{\bb\in\mathcal{B}[0]}\underset{0, \bb,Q_0}{\iber}\left[\frac{3^g \det(\mathrm{Hess}_\bb(Q_0))^{g-1} \cdot e^{\langle \lambda, x\rangle}}{W_{\!e^u\!}(x)^{g-1}\prod_{i<j}(1-e^{\langle \alpha_{j i}, x \rangle})(1-e^{\langle \alpha_{i j}, x \rangle-u})} \right](-[c]_{\bb}),$$ 
while its contribution  to the equivariant Euler characteristic of \(\mathcal L(k;\lambda)\otimes \pi_!(U_\nu\otimes K^{1/2})\) is
$$  \left. \frac{\mathrm{d}}{\mathrm{d} \delta}\left(\sum_{\bb\in\mathcal{B}[0]}\underset{0, \bb,Q_\delta}{\iber}\left[\frac{3^g \det(\mathrm{Hess}_\bb(Q_\delta))^{g-1} \cdot e^{\langle \lambda+v^{\det}, x\rangle}}{W_{\!e^u\!}(x)^{g-1}\prod_{i<j}(1-e^{\langle \alpha_{j i}, x \rangle})(1-e^{\langle \alpha_{i j}, x \rangle-u})} \right](-[c]_{\bb})\right)\right|_{\delta=0}.$$
\end{theorem}

\section{Contribution of type-$(1,1,1)$ components}\label{S:111}
In this section, we study the contribution of the fixed-point components of type-$(1,1,1)$ to the equivariant Euler characteristic \eqref{Euler}. For the result, see Theorem \ref{thm:Cont111}.

\subsection{Fixed-point components of type (1,1,1)}
We begin by describing these components explicitly.  By Proposition \ref{fixedcomp}, each fixed point of this type determines both a permutation $\sigma\in\Sigma_3$, specifying the assignment of parabolic weights to the three line-bundle summands, and the degrees of these summands.  For each $\sigma\in\Sigma_3$, we set
\begin{equation}\label{Sdescr111}
S_1^\sigma = \begin{cases} p, \text{\; if\;} \sigma(1)>\sigma(2); \\ \emptyset,  \text{\; otherwise, \;} \end{cases}  S_2^\sigma = \begin{cases} p, \text{\; if\;} \sigma(2)>\sigma(3); \\ \emptyset,  \text{\; otherwise,\;} \end{cases}  S_3^\sigma = \begin{cases} p, \text{\; if\;} \sigma(1)>\sigma(3); \\ \emptyset,  \text{\; otherwise,\;} \end{cases} 
\end{equation}
and  write $s_i^\sigma = |S_i^\sigma|$, $i=1,2,3$. Given integers $l_1,l_2$, we denote by $F_{l_1,l_2}^\sigma\subset\Hs_0((c_1,c_2,c_3))$ the locus of stable parabolic Higgs bundles of the form $(L_1\oplus L_2\oplus L_3, \Phi)$, where  $L_i$ carries parabolic weight $c_{\sigma(i)}$, and for $j=1,2$   $$\deg(L_j)=l_j, \;\;\; \Phi|_{L_j}\in H^0(\mathrm{SParHom}(L_j,L_{j+1}\otimes K(p)))\cong H^0(L_j^\vee\otimes L_{j+1}\otimes K(p-S_j^\sigma)),$$ while $L_3\simeq (L_1\otimes L_2)^\vee$ and $\Phi|_{L_3} = 0$. The nonvanishing of $\Phi|_{L_1}$ and $\Phi|_{L_2}$, together with the stability conditions for the $\Phi$-invariant subbundles  $L_3$ and $L_2\oplus L_3$, implies 
\begin{equation}\label{111admissibleconditions}
2g-1-s_1^\sigma \geq l_1-l_2, \;\; 2g-1-s_2^\sigma\geq 2l_2+l_1,\;\; -c_{\sigma(3)} < l_1+l_2 \text{\;\; and\;\;} c_{\sigma(1)} <l_1. 
\end{equation}
By passing to the zero divisors of $\Phi|_{L_1}$ and $\Phi|_{L_2}$ and using the relation $L_3\simeq(L_1\otimes L_2)^\vee,$ we obtain the following explicit description of $F^\sigma_{l_1,l_2}$.
\begin{proposition}\label{prop:range111}
Let $\sigma\in\Sigma_3$ and let $l_1, l_2\in\mathbb{Z}$ satisfy \eqref{111admissibleconditions}. Then the fixed-point locus $F_{l_1,l_2}^\sigma\subset \Hs_0(c)^\TT$ is isomorphic to the fiber product
$$F^\sigma_{l_1,l_2} := \Sym^{m_1}(C)\times \Sym^{m_2}(C)\times_{\Jac^{m_1+2m_2}(C)}\Jac^{-l_1-l_2}(C),$$ where $m_1=-l_1+l_2+2g-1-s_1^\sigma$, $m_2=-2l_2-l_1+2g-1-s_2^\sigma$, 
and the maps from $\Sym^{m_1}(C)\times \Sym^{m_2}(C)$ and $\Jac^{-l_1-l_2}(C)$ to ${\Jac^{m_1+2m_2}(C)}$ are given by $(D_1,D_2)\mapsto \mathcal{O}(D_1+2D_2)$ and $\mathcal{L}\mapsto \mathcal{L}^3\otimes K^3\otimes\mathcal{O}(3p-S_1^\sigma-2S_2^\sigma)$, respectively.
\end{proposition}

\begin{corollary}\label{32gcover}
In particular, the induced map $f\colon F^\sigma_{l_1,l_2}\to \Sym^{m_1}(C)\times \Sym^{m_2}(C)$ is a finite \'etale cover of degree $3^{2g}$, and for every class $\gamma\in H^\bullet(\Sym^{m_1}(C)\times \Sym^{m_2}(C))$, the integral of $f^*\gamma$ over $F^\sigma_{l_1,l_2}$  equals  $3^{2g}$ times the integral of $\gamma$ over $\Sym^{m_1}(C)\times\Sym^{m_2}(C)$.
\end{corollary}
Note that $F^\sigma_{l_1,l_2}$ is a union of connected components of $\Hs_0(c)^\TT$ and need not be connected.
From now on, we fix a permutation \(\sigma \in \Sigma_3\) and a pair of integers \(l_1, l_2\) satisfying the admissibility conditions in \eqref{111admissibleconditions}. To compute the contribution of $F^\sigma_{l_1,l_2}\subset \Hs_0(c)^\TT$ to the equivariant Euler characteristic \eqref{Euler}, we need to  evaluate the integrals 
\begin{equation}\label{int111}
\underset{F^\sigma_{l_1,l_2}}{\int} \frac{\ch_\TT\bigl(\mathcal{L}(k;\lambda)|_{F^\sigma_{l_1,l_2}}\bigr)}{e_\TT\bigl(\mathcal{N}^\vee_{F^\sigma_{l_1,l_2}}\bigr)} \mathrm{Td}(F^\sigma_{l_1,l_2}), \;\;\;\;  \underset{F^\sigma_{l_1,l_2}}{\int} \frac{\ch_\TT\bigl(\mathcal{L}(k;\lambda)\otimes\pi_!(U_\nu\otimes K^{1/2})|_{F^\sigma_{l_1,l_2}}\bigr)}{e_\TT\bigl(\mathcal{N}^\vee_{F^\sigma_{l_1,l_2}}\bigr)} \mathrm{Td}(F^\sigma_{l_1,l_2}).
\end{equation}
In the next subsections, we analyze the vector bundles appearing in the integrands and compute \eqref{int111} explicitly; see Proposition \ref{prop:fixedl111cont}  for the result.

\subsection{Vector bundles over the fixed-point locus $F^\sigma_{l_1,l_2}$}\label{S:vbunonF}
We begin by fixing notation for the universal bundles on $F^\sigma_{l_1,l_2}\times C$ and the associated  cohomology classes. 
For $j=1,2$, let $m_j$ be as in Proposition \ref{prop:range111}, and let $\Delta_j\subset\Sym^{m_j}(C)\times C$ denote the universal divisor. We write $\mathcal O(\Delta_j)$ for the associated line bundle and, by abuse of notation, also for its pullback to $F^\sigma_{l_1,l_2}\times C$.

Let $\omega\in H^2(C)$ be the class Poincar\'e dual to a point, and choose a symplectic basis $e_1,\dots,e_{2g}$ of $H^1(C)$, so that $e_i e_{i+g}=-e_{i+g}e_i=\omega$ for $i=1,\dots,g,$ and all other products vanish. We define the classes $\eta_j\in H^2(\Sym^{m_j}(C))$ and $\xi_i^j\in H^1(\Sym^{m_j}(C))$, for $i=1,\dots,2g$, by the K\"unneth decomposition
\begin{equation}\label{eq:ODelta}
c_1(\mathcal{O}(\Delta_j)) = \eta_j \otimes 1 + \sum_{i=1}^{2g} \xi_i^{j} \otimes e_i + m_j \otimes \omega\in H^2(\Sym^{m_j}(C)\times C),
\end{equation}
  regard them as classes on $F^\sigma_{l_1,l_2}$ via pullback, and set   $\theta_j := \sum_{i=1}^g \xi_i^{j} \xi_{i+g}^{j}$ and $\theta_{12} := \sum_{i=1}^g (\xi_i^{1} \xi_{i+g}^{2} - \xi_{i+g}^{1} \xi_i^{2}).$ 
With these conventions, we have
\begin{equation}\label{rel:Sym}
 \left( \sum_{i=1}^{2g} \xi_{i}^j \otimes e_i \right)^2 = -2\theta_j \otimes \omega \text {\; and \;} 
\left( \sum_{i=1}^{2g} \xi_{i}^1 \otimes e_i \right) \left( \sum_{i=1}^{2g} \xi_{i}^2 \otimes e_i \right) = -\theta_{12} \otimes \omega. 
\end{equation}

By construction, the induced $\TT$-linearization of each $\mathcal O(\Delta_j)$ has weight $-1$, and the restriction of the universal bundle $\mathcal U$ to $F^\sigma_{l_1,l_2}\times C$ splits as $\mathcal L_1\oplus\mathcal L_2\oplus\mathcal L_3$, with
\begin{equation}\label{eq:chLLL}
\mathcal{L}_j^\vee\otimes \mathcal{L}_{j+1}\simeq \mathcal{O}(\Delta_j)\otimes K^\vee(S_j^\sigma-p) \text {\; for $j=1,2$,\; and \;} \mathcal{L}_1\otimes \mathcal{L}_2\otimes\mathcal{L}_3\simeq \det\bigl(\mathcal{U}|_{F^\sigma_{l_1,l_2}\times\{p\}}\bigl). 
\end{equation}
Here and below, we use the same notation for bundles on $C$ and for their pullbacks to $F^\sigma_{l_1,l_2}\times C$.
Writing \(t=e^u\), we then have 
\begin{multline}\label{eq:gamma111} 
\;\;\;\;\;\;\;\;\; c_1^\TT(\mathcal L_1^\vee\otimes\mathcal L_{2}) = (\eta_1-u) +\sum_{i=1}^{2g}\xi^1_i\otimes e_i +(l_2-l_1)\omega \text{\;\; and \;\;} \\ 
c_1^\TT(\mathcal L_2^\vee\otimes\mathcal L_{3}) = (\eta_2-u) +\sum_{i=1}^{2g}\xi^2_i\otimes e_i +(-2l_2-l_1)\omega. \;\;\;\;\;\;\;\;\;
\end{multline}

We now describe the restriction of the line bundle $\mathcal{L}(k;\lambda)$ to $F^\sigma_{l_1,l_2}$. 
\begin{proposition}\label{LklonF} In the notation introduced above, 
\begin{multline*}
c_1^\TT(\mathcal{L}(k;(\lambda_1,\lambda_2,\lambda_3))|_{F^\sigma_{l_1,l_2}})= -\lambda_{\sigma(1)}(\eta_1-u) + \lambda_{\sigma(3)}(\eta_2-u) + \\
k(\eta_1-u)l_1 + k(\eta_2-u)(l_1+l_2)+\frac{2k}{3}\theta_1+\frac{2k}{3}\theta_2+\frac{k}{3}\theta_{12}. 
\end{multline*}
\end{proposition}
\begin{proof}
It follows immediately from the definition that
 $$c_1^\TT\bigl(\mathcal{L}(0;(\lambda_1,\lambda_2,\lambda_3))|_{F^\sigma_{l_1,l_2}}\bigr)=-\lambda_{\sigma(1)}(\eta_1-u) + \lambda_{\sigma(3)}(\eta_2-u).$$  
On the other hand, by Lemma \ref{c_2End2}, $c_1(\mathcal{L}(k;0))=\frac{k}{6}\pi_*c_2(\mathrm{End}(\mathcal U))$. 
Using the splitting of \(\mathcal U\), we therefore obtain
\begin{multline*}
c_1^\TT\bigl(\mathcal{L}(k;0)|_{F^\sigma_{l_1,l_2}}\bigr) =  -\frac{k}{6}\sum_{i<j}\pi_*(c^\TT_1(\mathcal{L}_i^\vee\otimes\mathcal{L}_j)^2) = \\ k(\eta_1-u)l_1+k(\eta_2-u)(l_1+l_2)+\frac{2k}{3}\theta_1+\frac{2k}{3}\theta_2+\frac{k}{3}\theta_{12},
\end{multline*}
where for the second equation we used  \eqref{eq:gamma111} and \eqref{rel:Sym}.
\end{proof}

Similarly, we compute the restriction of $\pi_!(U_\nu\otimes K^{1/2})$ to $F^\sigma_{l_1,l_2}$. 
\begin{proposition}\label{UnyuonF} Recall the notation $\widehat\phi^\nu(z_1,z_2)$ for the function $\phi^\nu(x)$ written in coordinates $z_1=\langle\alpha_{\sigma(2)\sigma(1)},x\rangle$ and $ z_2=\langle\alpha_{\sigma(3)\sigma(2)},x\rangle$, set $|\nu|=\nu_1+\nu_2+\nu_3$ and  $v^{\det}=\frac{|\nu|}3(1,1,-2)$.
Then
\begin{multline*}
 \ch_\TT\bigl(\pi_!(U_\nu\otimes K^{1/2})|_{F^\sigma_{l_1,l_2}}\bigr)  =e^{-v^{\det}_{\sigma(1)}(\eta_1-u)+v^{\det}_{\sigma(3)}(\eta_2-u)}\cdot 
 (-(l_1-l_2)\partial_{z_1} - (l_1+2l_2)\partial_{z_2} - \\ \theta_1\,\partial_{z_1}^2-\theta_{12}\,\partial_{z_1}\partial_{z_2} -\theta_2\,\partial_{z_2}^2)\,\widehat\phi^\nu(z_1,z_2)\Big|_{\substack{z_1=\eta_1-u \\ z_2=\eta_2-u }}.
 \end{multline*}
\end{proposition}
\begin{proof}
We note that on  $F^\sigma_{l_1,l_2}\times C$ the normalized universal bundle
$U=\mathcal U\otimes\mathcal F_1^\vee$ splits as
$\overline{\mathcal L}_1\oplus\overline{\mathcal L}_2\oplus\overline{\mathcal L}_3, $ where
$\overline{\mathcal L}_i=\mathcal L_i\otimes\mathcal{F}_1^\vee|_{F^\sigma_{l_1,l_2}}$. Thus, using the notation of \eqref{phinumu}, we obtain
\begin{equation}\label{eq:some1}
U_\nu\big|_{F^\sigma_{l_1,l_2}\times C}=\bigoplus_\mu\big(\overline{\mathcal L}_1^{\mu_1}\otimes\overline{\mathcal L}_2^{\mu_2}\otimes \overline{\mathcal L}_3^{\mu_3}\big)^{\oplus m_\mu}.
\end{equation}
Writing $\bar\mu= \mu-\frac{|\nu|}{3}(1,1,1)$  for each weight $\mu$, we have 
\begin{equation}\label{eq:some2}
\sum_{i=1}^3\mu_i\,c_1^\TT(\overline{\mathcal L}_i) =-\bar\mu_1\,c_1^\TT(\mathcal L_1^\vee\otimes\mathcal L_2) +\bar\mu_3\,c_1^\TT(\mathcal L_2^\vee\otimes\mathcal L_3)+\tfrac{|\nu|}3\,c_1^\TT(\det U).
\end{equation}
Since for every $V\in F^\sigma_{l_1,l_2}$ the restriction of $\det U$ to $\{V\}\times C$ is trivial, we have 
$c_1^\TT(\det U)= c_1^\TT(\mathcal L(0;(1,1,-2)))$, and thus by Proposition \ref{LklonF}
$$\tfrac{|\nu|}3\,c_1^\TT(\det U)\big|_{F^\sigma_{l_1,l_2}} =-v^{\det}_{\sigma(1)}(\eta_1-u)+v^{\det}_{\sigma(3)}(\eta_2-u).$$

Combining this with  \eqref{eq:some1}, \eqref{eq:some2}, and using \eqref{eq:gamma111}, we obtain 
\begin{multline*}
  \ch_\TT\bigl(U_\nu|_{F^\sigma_{l_1,l_2}\times C}\bigr)   =e^{-v^{\det}_{\sigma(1)}(\eta_1-u)+v^{\det}_{\sigma(3)}(\eta_2-u)}
  \sum_\mu m_\mu\,e^{-\bar\mu_1(\eta_1-u)+\bar\mu_3(\eta_2-u)}\cdot\\
  \exp\Big(\sum_{i=1}^{2g}\big(-\bar\mu_1\xi_i^1+\bar\mu_3\xi_i^2\big)\otimes e_i
  -\bar\mu_1(l_2-l_1)\,\omega-\bar\mu_3(l_1+2l_2)\,\omega\Big).
\end{multline*}
Expanding the last exponential and using \eqref{rel:Sym}, we arrive at  
\begin{multline*}
  \ch_\TT\big(\pi_!(U_\nu\otimes K^{1/2})|_{F^\sigma_{l_1,l_2}}\big) = \pi_* \ch_\TT\big(U_\nu|_{F^\sigma_{l_1,l_2}}\big)
  =e^{-v^{\det}_{\sigma(1)}(\eta_1-u)+v^{\det}_{\sigma(3)}(\eta_2-u)} \cdot \\
  \sum_\mu m_\mu\,e^{-\bar\mu_1(\eta_1-u)+\bar\mu_3(\eta_2-u)}
  \bigl(-\bar\mu_1(l_2-l_1)-\bar\mu_3(l_1+2l_2)
  -\bar\mu_1^2\,\theta_1+\bar\mu_1\bar\mu_3\,\theta_{12}-\bar\mu_3^2\,\theta_2\bigr).
\end{multline*}
Finally, observing that  $\partial_{z_1}$ and $\partial_{z_2}$ act on
$e^{-\bar\mu_1z_1+\bar\mu_3z_2}$ by multiplication by $-\bar\mu_1$ and
$\bar\mu_3$, respectively, we arrive at the claimed result.
\end{proof}

\subsection{The normal bundle to $F^\sigma_{l_1,l_2}$}\label{S:normF}
Our next goal is to compute the inverse of the equivariant $K$-theoretic Euler class of the conormal bundle to $F^\sigma_{l_1,l_2}$ in $\Hs_0(c)$. We begin by studying the restriction of the tangent bundle $T\Hs_0(c)$ to $F^\sigma_{l_1,l_2}$. By Theorem \ref{thm:modsp}(iii), we have
$$\ch_\TT(T\Hs_0(c)) = -\ch_\TT(\pi_!(\mathrm{ParEnd}_0(\mathcal{U}))) + t\ch_\TT(\pi_!(\mathrm{SParEnd}_0(\mathcal{U})\otimes K(p))).$$ The Grothendieck-Riemann-Roch theorem then yields
\begin{equation}\label{GRRforT}
\ch_\TT(T\Hs_0(c))= -\pi_*(\ch_\TT(\mathrm{ParEnd}_0(\mathcal{U}))(1-(g-1)\omega)) +t\pi_*(\ch_\TT(\mathrm{SParEnd}_0(\mathcal{U}))(1+g\omega)).
\end{equation}
We now evaluate this expression on  \(F^\sigma_{l_1,l_2}\).
\begin{lemma}\label{l:eqT}
The restriction of the equivariant Chern character of the tangent bundle $T\Hs_0(c)$ to $F^\sigma_{l_1,l_2}$ is given by
\begin{multline}\label{eq:bigchT}
2(g-1)+ e^{\eta_1}(l_2-l_1-\theta_1-s_1^\sigma+g)+e^{\eta_2}(-2l_2-l_1-\theta_2-s_2^\sigma+g)+ \\ 2t(g-1) +
t^{-1}e^{\eta_1}(l_1-l_2+\theta_1+s^\sigma_1+g-1) +  t^{-1}e^{\eta_2}(l_1+2l_2+\theta_2+s^\sigma_2+g-1) + \\
te^{-\eta_1}(l_2-l_1+\theta_1+g-s_1^\sigma) + t^2e^{-\eta_1-\eta_2}(-l_2-2l_1+\theta_1+\theta_2+\theta_{12}+g-s_3^\sigma) + \\
 te^{-\eta_2}(-2l_2-l_1+\theta_2+g-s_2^\sigma) + t^{-2}e^{\eta_1+\eta_2}(l_2+2l_1+\theta_1+\theta_2+\theta_{12}+s_3^\sigma+g-1) + \\
t^3e^{-\eta_1-\eta_2}(l_2+2l_1-\theta_1-\theta_2-\theta_{12}+s_3^\sigma+g-1) + t^2e^{-\eta_2}(l_1+2l_2-\theta_2+s_2^\sigma+g-1) + \\
t^{-1}e^{\eta_1+\eta_2}(-2l_1-l_2-\theta_1-\theta_2-\theta_{12}-s_3^\sigma+g) + t^2e^{-\eta_1}(l_1-l_2-\theta_1+s_1^\sigma+g-1).
\end{multline}
\end{lemma}
\begin{proof}
Note that on \(F^\sigma_{l_1,l_2}\times C\) one has
\begin{multline*}
\ch_\TT(\mathrm{ParEnd}_0(\mathcal{U})) =\ch_\TT(\mathrm{SParEnd}_0(\mathcal{U})) + 2\omega = \\ \ch_\TT(\mathrm{End}_0(\mathcal{U}))  - \omega(t^{-1}e^{\eta_1})^{2s_1^\sigma-1} -  \omega(t^{-1}e^{\eta_2})^{2s_2^\sigma-1}  - 
 \omega (t^{-2}e^{\eta_1+\eta_2})^{2s_3^\sigma-1}.
 \end{multline*}  
 Substituting this into \eqref{GRRforT} yields
\begin{multline}\label{eq:impbhab}
\ch_\TT\bigl(T\Hs_0(c)|_{F^\sigma_{l_1,l_2}}\bigr) = (t-1)\pi_*\ch_\TT\bigl(\mathrm{End}_0(\mathcal{U})|_{F^\sigma_{l_1,l_2}\times C}\bigr)  - 2t + \\  
(1-s_1^\sigma)(te^{-\eta_1}-t^2e^{-\eta_1}) +s_1^\sigma(t^{-1}e^{\eta_1}-e^{\eta_1}) + (1-s_2^\sigma)(te^{-\eta_2}-t^2e^{-\eta_2}) +  \\
s_2^\sigma(t^{-1}e^{\eta_2}-e^{\eta_2}) + s_3^\sigma(t^{-2}e^{\eta_1+\eta_2}-t^{-1}e^{\eta_1+\eta_2})+(1-s_3^\sigma)(t^2e^{-\eta_1-\eta_2}-t^3e^{-\eta_1-\eta_2}) +  \\ 
(g-1+tg)(2+t^{-1}e^{\eta_1}+t^{-1}e^{\eta_2}+t^{-2}e^{\eta_1+\eta_2}+te^{-\eta_1}+te^{-\eta_2}+t^2e^{-\eta_1-\eta_2}).
\end{multline}
On the other hand, \eqref{eq:chLLL} together with \eqref{rel:Sym} give
\begin{multline}\label{eq:bfin}
\pi_*\ch_\TT\bigl(\mathrm{End}(\mathcal{U})|_{F^\sigma_{l_1,l_2}\times C}\bigr)= t^{-1}e^{\eta_1}(l_2-l_1-\theta_1)+ te^{-\eta_1}(l_1-l_2-\theta_1)+  \\ t^{-1}e^{\eta_2}(-2l_2-l_1-\theta_2)+ te^{-\eta_2}(2l_2+l_1-\theta_2) + \\ t^{-2}e^{\eta_1+\eta_2}(-l_2-2l_1-\theta_1-\theta_2-\theta_{12})+ t^2e^{-\eta_1-\eta_2}(l_2+2l_1-\theta_1-\theta_2-\theta_{12}).
\end{multline}
Finally, substituting \eqref{eq:bfin} into \eqref{eq:impbhab} gives \eqref{eq:bigchT}.
\end{proof}

The expression for the inverse Euler class of the conormal bundle now follows directly.
\begin{proposition}\label{ElonF} In the notation introduced above, 
\begin{multline}\label{eq:ElonF}
e_{\TT}(\mathcal{N}^\vee_{F^\sigma_{l_1,l_2}})^{-1} = \frac{\exp(-\frac{\theta_1}{1-t^2e^{-\eta_1}}-\frac{\theta_2}{1-t^2e^{-\eta_2}} - \frac{\theta_1+\theta_2+\theta_{12}}{1-t^{-1}e^{\eta_1+\eta_2}}-\frac{\theta_1+\theta_2+\theta_{12}}{1-t^3e^{-\eta_1-\eta_2}})}{(1-t^{-1})^{2g-2}(1-t^{-2}e^{\eta_1})^{l_1-l_2+s_1^\sigma+g-1}(1-t^{-2}e^{\eta_2})^{2l_2+l_1+g-1+s_2^\sigma}   } \cdot \\
\frac{\exp(2\theta_1+2\theta_2+\theta_{12})}{(1-t^{-3}e^{\eta_1+\eta_2})^{l_2+2l_1+s_3^\sigma+g-1}(1-te^{-\eta_1})^{g-1+s_1^\sigma}(1-t^{-1}e^{\eta_1})^{g-s_1^\sigma} (1-te^{-\eta_2})^{g-1+s_2^\sigma}   }\cdot \\
\frac{(t^{-1}e^{\eta_1})^{3l_1}(t^{-1}e^{\eta_2})^{3l_1+3l_2}}{(1-te^{-\eta_1-\eta_2})^{-l_2-2l_1-s_3^\sigma+g}(1-t^{-1}e^{\eta_2})^{g-s_2^\sigma}(1-t^2e^{-\eta_1-\eta_2})^{g-1+s_3^\sigma}(1-t^{-2}e^{\eta_1+\eta_2})^{g-s_3^\sigma}  }.
\end{multline}
\end{proposition}
\begin{proof}
The weight-zero part of \(T\Hs_0(c)|_{F^\sigma_{l_1,l_2}}\) is \(TF^\sigma_{l_1,l_2}\). Removing it from \eqref{eq:bigchT} therefore gives the equivariant Chern character of the normal bundle $\mathcal{N}_{F^\sigma_{l_1,l_2}}$.
To pass from the Chern character to the inverse \(K\)-theoretic Euler class, we use the following statement. 
\begin{lemma}\label{thetato}
Let $f_\tau(x)=\frac{1}{1-\tau^{-1}e^{-x}}$. If $\theta$ is a cohomology class of degree $2$, then the contribution of a Chern-character term
$\theta e^x$ to the multiplicative class associated with $f_\tau$ is $\exp\left(\frac{\theta}{1-\tau e^{x}}\right).$
\end{lemma}
\begin{proof}
Perturbing a Chern-character term by $e^x\mapsto e^{x+\varepsilon\theta}$, we note that the corresponding first-order variation of the logarithm of the multiplicative class is $ \frac{\theta}{1-\tau e^{x}}.$  Exponentiating, we obtain the result.
\end{proof}
Applying this identity to the $\theta$-terms in \eqref{eq:bigchT}: 
$\theta_1(t^{-1}e^{\eta_1} + te^{-\eta_1}-t^2e^{-\eta_1})$, 
 $\theta_2(t^{-1}e^{\eta_2} + te^{-\eta_2}-t^2e^{-\eta_2}),$ and  
$(\theta_1+\theta_2+\theta_{12}) (t^2e^{-\eta_1-\eta_2}+t^{-2}e^{\eta_1+\eta_2}-t^3e^{-\eta_1-\eta_2}-t^{-1}e^{\eta_1+\eta_2})$, yields \eqref{eq:ElonF}.
\end{proof}

\subsection{Contribution of the fixed-point locus  $F^\sigma_{l_1,l_2}$}
We now turn to the explicit evaluation of the contribution of  $F^\sigma_{l_1,l_2}$ to the equivariant Euler characteristic \eqref{Euler}. Recall that $\delta$ denotes a formal parameter, and $u$ is the additive equivariant parameter: $e^u=t$.  For $\widehat{\lambda}\in\Lambda\otimes\mathbb{Q}$, consider the integral 
\begin{multline}\label{int111-1}
\int_{F^\sigma_{l_1,l_2}}\frac{\exp(-\frac{\theta_1}{1-e^{2u-\eta_1}}-\frac{\theta_2}{1-e^{2u-\eta_2}} - \frac{\theta_1+\theta_2+\theta_{12}}{1-e^{\eta_1+\eta_2-u}}-\frac{\theta_1+\theta_2+\theta_{12}}{1-e^{3u-\eta_1-\eta_2}})e^{(-\lhat_{\sigma(1)}\!+l_1(k+3))(\eta_1\!-u)} }{(1-e^{-u})^{2g-2}(1-e^{\eta_1-2u})^{l_1-l_2+s_1^\sigma+g-1}(1-e^{\eta_2-2u})^{2l_2+l_1+g-1+s_2^\sigma}   } \cdot \\
 \frac{e^{(\lhat_{\sigma(3)}\!+(l_1+l_2)(k+3))(\eta_2\!-u)}e^{{-\delta(l_1-l_2)\partial_{1}\!\widehat\phi^\nu\!(\eta_1\!-\!u,\eta_2\!-\!u)}}e^{{-\delta (l_1+2l_2)\partial_{2}\!\widehat\phi^\nu\!(\eta_1\!-\!u,\eta_2\!-\!u)}}  
 }{(1-e^{\eta_1+\eta_2-3u})^{l_2+2l_1+s_3^\sigma+g-1}(1-e^{u-\eta_1})^{g-1+s_1^\sigma}(1-e^{\eta_1-u})^{g-s_1^\sigma} (1-e^{u-\eta_2})^{g-1+s_2^\sigma}   }\cdot \\
\frac{   \exp\big((\frac{k}{3}+1)(2\theta_1\!+2\theta_2\!+\theta_{12}) -\delta(\theta_1\partial_{1}^2\!+\theta_{12}\partial_{1}\partial_{2}\! +\theta_2 \partial_{2}^2)\widehat\phi^\nu\!(\eta_1\!-\!u,\eta_2\!-\!u)\big)   \,\mathrm{Td}(F^\sigma_{l_1,l_2})}{(1-e^{u-\eta_1-\eta_2})^{-l_2-2l_1-s_3^\sigma+g}(1-e^{\eta_2-u})^{g-s_2^\sigma}(1-e^{2u-\eta_1-\eta_2})^{g-1+s_3^\sigma}(1-e^{\eta_1+\eta_2-2u})^{g-s_3^\sigma} },
\end{multline}
where we use the notation $\partial_j\widehat{\phi}^\nu(\eta_1-u,\eta_2-u):=\partial_{z_j}\widehat{\phi}^\nu(z_1,z_2)|_{z_1=\eta_1-u,\, z_2=\eta_2-u}$. 

By Propositions \ref{LklonF} and \ref{ElonF}, for $\widehat{\lambda}=\lambda$ and $\delta=0$, the integral \eqref{int111-1} is precisely the contribution of the fixed-point locus $F^\sigma_{l_1,l_2}\subset\Hs_0(c)^\TT$ to the equivariant Euler characteristic of $\mathcal L(k;\lambda)$ in \eqref{Euler}. Combining these results with Proposition \ref{UnyuonF}, a simple calculation shows that, for $\widehat{\lambda}=\lambda+v^{\mathrm{det}}$, the coefficient of \(\delta\) in \eqref{int111-1} gives the contribution of $F^\sigma_{l_1,l_2}$ to the equivariant Euler characteristic of $\mathcal L(k;\lambda)\otimes\pi_!(U_\nu\otimes K^{1/2})$  in \eqref{Euler}.

It therefore remains to evaluate \eqref{int111-1} explicitly. Putting together Lemma \ref{l:intsym}(i) and Corollary \ref{32gcover}, we can rewrite this integral as
\begin{multline}\label{int111-2}
\!\!\underset{\Sym^{m_1}\!(C)\!\times\Sym^{m_2}\!(C)\!}{\int}\!\!\!\! \frac{\eta_1^{m_1\!-\!g\!+\!1}\eta_2^{m_2\!-\!g\!+\!1}e^{-\lhat_{\sigma(1)}(\eta_1\!-\!u)\!+\! \lhat_{\sigma(3)}(\eta_2\!-\!u) }   e^{(\eta_1\!-\!u)(k+3)l_1}e^{(\eta_2\!-\!u)(k+3)(l_1+l_2)}}{W_{\!e^u}(\eta_1\!-\!u, \eta_2\!-\!u)^{g\!-\!1} (1-e^{-\!\eta_2})^{1\!-\!s_2^\sigma}   (1-e^{\eta_1\!-\!2u})^{s_1^\sigma}(1-e^{-\!\eta_1})^{1\!-\!s_1^\sigma}} \\
\frac{\exp\big(-\frac{\theta_1}{\eta_1}\!-\!\frac{\theta_2}{\eta_2}\!+\!(\theta_1\!+\!\theta_2\!+\!\theta_{12})G_{12}(\eta_1\!-\!u,\eta_2\!-\!u)\!+\!\theta_1G_{1\!}(\eta_1\!-\!u,\eta_2\!-\!u)\!+\!\theta_2G_2(\eta_1\!-\!u,\eta_2\!-\!u)\big)}{(1-e^{\eta_2 -2u})^{s_2^\sigma}(1-e^{\eta_1+\eta_2-3u})^{s_3^\sigma}(1-e^{-\eta_1-\eta_2+u})^{1-s_3^\sigma}(1-e^{\eta_1-u})^{1\!-\!s_1^\sigma}} \\
 \frac{3^{2g}\phi(\eta_1\!-\!u)^{l_1\!-l_2}\phi(\eta_2\!-\!u)^{l_1\!+\!2l_2}\phi(\eta_1\!+\!\eta_2\!-\!2u)^{2l_1\!+\!l_2} e^{-\delta((l_1\!-l_2)\partial_{1}+(l_1\!+2l_2)\partial_{2})\widehat\phi^\nu\!(\eta_1\!-\!u,\eta_2\!-\!u)}  }{(1-e^{u-\eta_1})^{s_1^\sigma}(1-e^{\eta_2-u})^{1-s_2^\sigma}(1-e^{u-\eta_2})^{s_2^\sigma}(1-e^{\eta_1+\eta_2\!-\!2u})^{1\!-\!s_3^\sigma}(1-e^{2u-\eta_1-\eta_2})^{s_3^\sigma}},
\end{multline}
where $m_j$, $j=1,2$, are as in Proposition \ref{prop:range111},  and  the functions $W_{e^u}(x,y)$, $G_1(z_1,z_2)$, $G_2(z_1,z_2)$, $G_{12}(z_1,z_2)$ and $\phi(z)$ are defined in \S\ref{SS:Main}. 
To perform this integral, we  apply  Lemma \ref{l:intsym}(ii) successively to the two symmetric-product factors in \eqref{int111-2}. First, integrating over $\Sym^{m_2}(C)$, we take 
\begin{itemize}
\item  $B(\eta_2)=-\frac{1}{\eta_2}+G_2(\eta_1-u, \eta_2-u)+G_{12}(\eta_1-u,\eta_2-u)$ and 
\item $b_i(\eta_2) = -G_{12}(\eta_1-u,\eta_2-u)\xi_{g+i}^1$ \\ 
$b_{i+g}(\eta_2) = G_{12}(\eta_1-u,\eta_2-u)\xi_i^{1},$ $\;i=1,2,...,g.$
\end{itemize}
For the remaining integral over \(\Sym^{m_1}(C)\),  we apply the same formula with
\begin{itemize}
\item $B(\eta_1) = -\frac{1}{\eta_1} \!+ \!G_1(\eta_1\!- u,\eta_2\!-u)\! +\!G_{12}(\eta_1\!-u, \eta_2\! - u)\!-  \! \frac{G_{12}(\eta_1-u, \eta_2-u)^2}{G_2(\eta_1\!-\!u,\eta_2\!-\!u) + G_{12}(\eta_1\!-\!u, \eta_2\!-\!u)}$ 
\item and $b_i(\eta_1) = 0$ for all $i$.
\end{itemize}
After the change of variables $x=\eta_1-u,$  $y=\eta_2-u,$ we obtain the following result.

\begin{proposition}\label{prop:fixedl111cont} Let \(\sigma\in\Sigma_3\), and let \(l_1, l_2\in\mathbb{Z}\) satisfy \eqref{111admissibleconditions}. Set $H(x,y)= G_{1\!}(x,y)G_{2\!}(x,y)+G_{1\!}(x,y)G_{12\!}(x,y)+G_{2\!}(x,y)G_{12\!}(x,y)$ and 
\begin{multline*}
\!\!\Omega_{l_1,l_2}^\sigma\!= \!\frac{3^{2g} e^{-\lhat_{\sigma(1)}x+\lhat_{\sigma(3)}y }e^{(k+3)l_1x+(k+3)(l_1+l_2)y}\phi(x)^{l_1-l_2+s_1^\sigma}\phi(y)^{l_1+2l_2+s_2^\sigma}}{W_{\!e^u\!}(x,y)^{g-1}(1-e^{x})^{1-s_1^\sigma\!}(1-e^{-x})^{s_1^\sigma\!}(1-e^{y})^{1-s_2^\sigma\!}(1-e^{-y})^{s_2^\sigma\!}} \\
\frac{\phi(x+y)^{2l_1+l_2+s_3^\sigma} e^{-\delta((l_1\!-l_2)\partial_{x}+(l_1\!+2l_2)\partial_{y})\widehat\phi^\nu\!(x,y)} \,H(x,y)^g\, \mathrm{d}x\, \mathrm{d}y}{(1-e^{x+y})^{1-s_3^\sigma\!}(1-e^{-x-y})^{s_3^\sigma\!}(1-e^{-x-u})(1-e^{-y-u})(1-e^{-x-y-u})}.
\end{multline*}
With $t=e^u\!$, the contribution of the fixed-point locus \(F^\sigma_{l_1,l_2}\) to $\chi_t(\Hs_0(c), \mathcal{L}(k;\lambda))$ 
  is given by $$\underset{x=-u}{\res}\,\underset{y=-u}{\res}\Omega_{l_1,l_2}^\sigma\big|_{\lhat=\lambda, \, \delta=0},$$ while its contribution to $\chi_t(\Hs_0(c), \mathcal{L}(k;\lambda)\otimes\pi_!(U_\nu\otimes K^{1/2}))$
   is given by the coefficient of $\delta$ in $$\underset{x=-u}{\res}\,\underset{y=-u}{\res}\Omega_{l_1,l_2}^\sigma\big|_{{\lhat=\lambda+v^{\det}}}.$$
\end{proposition}

\subsection{Contribution of type-$(1,1,1)$ fixed-point components}\label{S:111trick}
We now turn to the computation of the total contribution of the fixed-point components of parabolic type $(\{\sigma(1)\},\{\sigma(2)\},\{\sigma(3)\})$ to the equivariant Euler characteristics \eqref{Euler}. By Propositions \ref{prop:fixedl111cont} and \ref{prop:range111}, this amounts to evaluating the sum
\begin{equation}\label{eq:111triv1}
\sum_{l_1,l_2} \underset{x=-u}{\res}\, \underset{y=-u}{\res}\Omega_{l_1,l_2}^\sigma
\end{equation}
over the pairs of integers \(l_1,l_2\) satisfying the conditions in \eqref{111admissibleconditions}.

In fact, we can simplify the summation range using the following observations. Note that the form \(\Omega_{l_1,l_2}^\sigma\) has a pole at \(y=-u\)  of order at most $2g-s_2^\sigma-2l_2-l_1.$
Hence $$\underset{y=-u}{\res}\Omega_{l_1,l_2}^\sigma=0, \text{ \; if \;} 2l_2+l_1>2g-1-s_2^\sigma.$$
Similarly, after taking the \(y\)-residue, the resulting form has a pole at \(x=-u\) of order at most  $2g-s_1^\sigma-l_1+l_2,$ so that $$\underset{x=-u}{\res}\,\underset{y=-u}{\res}\,\Omega_{l_1,l_2}^\sigma=0 \text{ \; if \;} l_1-l_2>2g-1-s_1^\sigma.$$
Hence the upper bounds in \eqref{111admissibleconditions} may be omitted without changing the sum \eqref{eq:111triv1}, and, after the change of variables $a=l_1$, $b=l_1+l_2$, \eqref{eq:111triv1} becomes
\begin{multline*}
\!\!\sum_{b>-c_{\sigma(3)}} \sum_{a>c_{\sigma(1)}}\!\underset{x=-u}{\res}\,\underset{y=-u}{\res}\frac{3^{2g} e^{\lhat_{\sigma(3)}y-\lhat_{\sigma(1)}x}e^{(k+3)(ax+by)} \phi(x)^{2a-b+s_1^\sigma}
}{W_{\!e^u\!}(x,y)^{g-1}(1-e^{x})^{1-s_1^\sigma\!}(1-e^{-x})^{s_1^\sigma\!}(1-e^{y})^{1-s_2^\sigma\!}(1-e^{-y})^{s_2^\sigma\!}} \\
\frac{\phi(y)^{2b-a+s_2^\sigma}\phi(x+y)^{a+b+s_3^\sigma} e^{-\delta(a(2\partial_x-\partial_y)+b(2\partial_y-\partial_x))\widehat{\phi}^\nu(x,y)}H(x,y)^g \, \mathrm{d}x\, \mathrm{d}y}{(1-e^{x+y})^{1-s_3^\sigma\!}(1-e^{-x-y})^{s_3^\sigma\!}(1-e^{-y-u})(1-e^{-x-y-u})(1-e^{-x-u})}. 
\end{multline*}

Substituting the values of \(s_i^\sigma\) for the six permutations \(\sigma\in\Sigma_3\) and using Lemmas \ref{lem:phidef} and \ref{lem:Gdef}, we arrive at the following result.
\begin{theorem}\label{thm:Cont111}
Let \(\sigma\in\Sigma_3\) and let $p\in\mathrm{Poles}(\A)$ be the pole with coordinates $(p_{\sigma(1)},p_{\sigma(2)},p_{\sigma(3)})=(u,0,-u)$. Set $\bb:=(\alpha_{\sigma(2) \sigma(1)}, \alpha_{\sigma(3) \sigma(2)})\in\Bases(\A_{p})$. In the notation of Theorem \ref{thm:main}, the total contribution of the fixed-point components of parabolic type $(\{\sigma(1)\},\{\sigma(2)\},\{\sigma(3)\})$
to  the equivariant Euler characteristic of \(\mathcal L(k;\lambda)\) on $\Hs_0(c)$ is
$$\underset{p, \bb,Q_0}{\iber}\!{}^{(\!+\!, +\!)}\left[\frac{3^g\cdot \det(\mathrm{Hess}_\bb(Q_0))^{g-1} \cdot e^{\langle \lambda, x\rangle}}{W_{\!e^u\!}(x)^{g-1}\prod_{i<j}(1-e^{\langle \alpha_{j i}, x \rangle})(1-e^{\langle \alpha_{i j}, x \rangle-u})} \right](-[c]_{\bb}),$$
while its contribution  to the equivariant Euler characteristic of \(\mathcal L(k;\lambda)\otimes \pi_!(U_\nu\otimes K^{1/2})\) is
$$  \left. \frac{\mathrm{d}}{\mathrm{d} \delta}\left(\underset{p, \bb,Q_\delta}{\iber}\!{}^{(\!+\!, +\!)}\left[\frac{3^g\cdot \det(\mathrm{Hess}_\bb(Q_\delta))^{g-1} \cdot e^{\langle \lambda+v^{\det}, x\rangle}}{W_{\!e^u\!}(x)^{g-1}\prod_{i<j}(1-e^{\langle \alpha_{j i}, x \rangle})(1-e^{\langle \alpha_{i j}, x \rangle-u})} \right](-[c]_{\bb})\right)\right|_{\delta=0}.$$
\end{theorem}

\section{Moduli spaces of stable parabolic triples}\label{S:partriplesgen}
In this section, we study moduli spaces of stable parabolic triples, which will be used to compute the contributions  to the equivariant Euler characteristics \eqref{Euler} of the connected components of the fixed-point locus $\Hs_0(c)^\TT$ of types $(2,1)$ and $(1,2)$.

We fix a regular element \(c=(c_1,c_2,c_3)\in\Delta\) and a subset \(\{\iio, \iit\}\subset\{1,2,3\}\), with $\iio<\iit$, and denote by $\jjo$ the remaining index. Throughout this section, \(E\) denotes a rank-\(2\) parabolic bundle of degree \(l\) with parabolic weights ${\hat c}=(c_{i_1}, c_{i_2})$, while \(L\) denotes a parabolic line bundle of degree \(2g-2-l\) with parabolic weight \(c_{\jjo}\).
\subsection{Moduli spaces of stable parabolic triples}
 We begin by briefly recalling the  definitions and basic  properties of  parabolic triples needed below; for further details, we refer to  \cite{GPGM, GPLM}.

A \textit{parabolic triple} is a triple $T=(L,E,\phi),$  where $L$ and $E$ are the parabolic bundles described above and $\phi:E\to L(p)$ is a strongly parabolic morphism. 
A \textit{subtriple} \(T'=(L',E',\phi')\subset T\) consists of parabolic subbundles \(L'\subset L\) and \(E'\subset E\) such that $\phi$ restricts to $E'$ as \(\phi'\colon E' \to L'(p)\); it is \textit{proper} if \(T'\neq0\) and \(T'\neq T\).

Given  \(\alpha\in\mathbb R\), we define the \(\alpha\)-\emph{degree} and \(\alpha\)-\emph{slope} of \(T'=(L',E',\phi')\) as
\[ \deg_\alpha(T'):= \mathrm{pardeg}(L'\oplus E')+\alpha\,\mathrm{rk}(E') \text{\;\; and \;\;}  \mu_\alpha(T'):= \frac{\deg_\alpha(T')}   {\mathrm{rk}(L'\oplus E')}. \]
The triple \(T\) is called  \emph{\(\alpha\)-stable} (respectively, \emph{\(\alpha\)-semistable}) if \(\mu_\alpha(T')<\mu_\alpha(T)\) (respectively, \(\mu_\alpha(T')\leq\mu_\alpha(T)\)) for every proper subtriple \(T'\subset T\), and \emph{\(\alpha\)-polystable} if it is a direct sum of \(\alpha\)-stable triples of the same \(\alpha\)-slope.

\begin{theorem}\label{thm:moduli21} Set  $\alpha_\mathrm{min} := 2g-2-\tfrac{3l}{2}+\tfrac{3}{2}(c_{\iio}+c_{\iit})$ and  $\alpha_{\mathrm{large}}:= 4\alpha_\mathrm{min}+ 3.$
\begin{enumerate}[label=(\roman*), leftmargin=1.5pc]
\item There exists a projective coarse moduli space $N^\alpha_{l}({\hat c})$  whose points correspond to isomorphism classes of \(\alpha\)-polystable parabolic triples $T=(L,E,\phi)$ of the topological and parabolic type described above, satisfying the fixed-determinant condition $\det(E)\otimes L\simeq K$. The locus $N_{l}^{\alpha,\mathrm{st}}({\hat c}) \subset N_{l}^\alpha({\hat c})$ parametrizing \(\alpha\)-stable triples is open.
\item For every $\alpha>\alpha_{\min}$, any $\alpha$-semistable triple $T=(L,E,\phi)$ has a nonzero morphism $\phi$, and the stable locus $N_l^{\alpha,\mathrm{st}}(\hat{c})$ is smooth.
\item A necessary condition for the existence of a \(\alpha\)-stable triple is $\alpha_\mathrm{min}<\alpha<\alpha_\mathrm{large}.$
\end{enumerate}
\end{theorem}

 \subsection{GIT construction of the moduli space $N_l^\alpha({\hat c})$}\label{S:GIT}
We now recall the GIT construction of the moduli space $N_l^\alpha({\hat c})$ for rational $\alpha>\alpha_{\mathrm{min}}$, following Schmitt \cite{DecoSch}. 
 Choose $n\gg0$ such that, for every $\alpha$-semistable triple $T=(L,E,\phi)$, the bundles $E(np)$ and $L(np)$  are globally generated and have vanishing $H^1$; set $\chi:=\dim H^0(E(np))$.
 Let $Q$ be the open subscheme of the Quot scheme parametrizing locally free quotients $$q:\mathbb{C}^\chi\otimes\mathcal O_C(-np)\twoheadrightarrow E$$  for which the induced map $\mathbb{C}^\chi\to H^0(E(np))$ is an isomorphism, and denote by $\tilde{\mathcal E}$ the universal quotient on $Q\times C$.

 Let $FQ=\operatorname{Flag}(\tilde{\mathcal E}\big|_{Q\times \{p\}})$ be the relative flag variety parametrizing parabolic structures on $E$, and denote the corresponding universal flag  by $$\tilde{\mathcal G}_1\subset \tilde{\mathcal G}_2=\tilde{\mathcal E}\big|_{Q\times \{p\}}.$$
Let $\tilde{\mathcal L}$ be the Poincar\'e bundle on $\operatorname{Jac}^{2g-2-l}(C)\times C$ normalized so that its restriction to  ${\Jac^{2g-2-l}\times\{p\}}$ is trivial. By abuse of notation, we continue to denote by $\tilde{\mathcal E}$, $\tilde{\mathcal L}$ and $\tilde{\mathcal G}_1\subset \tilde{\mathcal G}_2$ their pullbacks to the various spaces occurring below.

To impose the fixed-determinant condition, we set $B\subset FQ\times\operatorname{Jac}^{2g-2-l}(C)$ to be the closed subscheme parametrizing pairs $(E,L)$ satisfying
$\det E\otimes L\simeq K$. Let $\pi_B\colon B\times C\to B$ denote the projection, and define
$$\mathcal{H}:=\mathcal{H}\mathrm{om}\left(\mathcal{O}_B^\chi, \pi_{B*}(\tilde{\mathcal{L}}(p+np))\right).$$
Over \(\mathbb P(\mathcal H)\times C\) there is a tautological morphism
\begin{equation}\label{eq:tauth}
h\colon\mathbb{C}^\chi\otimes\mathcal{O}_{\mathbb{P}\mathcal{H}\times C}(-np)\to \tilde{\mathcal{L}}(p)\otimes\pi_{\mathbb{P}\mathcal{H}}^*\mathcal{O}_{\mathbb{P}\mathcal{H}}(1), \text{\; with \;}  \pi_{\mathbb{P}\mathcal{H}}\colon \mathbb{P}\mathcal{H}\times C\to \mathbb{P}\mathcal{H}.
\end{equation}
We introduce an auxiliary torus $\TT=\mathbb C^*$, which acts trivially on $B$ and by scalar multiplication on the fibers of $\mathcal H$. Thus the induced $\TT$-action on $\mathbb P(\mathcal H)$ is trivial, and the tautological line bundle $\mathcal{O}_{\mathbb P(\mathcal H)}(-1)$ carries $\TT$-weight $1$.  We endow both $\widetilde{\mathcal E}$ and $\widetilde{\mathcal L}$ with the trivial $\TT$-linearization.

Now let \(X\subset\mathbb P(\mathcal H)\) be the closed subscheme on which $h$ factors through the pullback of the universal quotient and the induced morphism is strongly parabolic. Thus, over \(X\times C\), we have the universal morphism
\begin{equation}\label{eq:tautphitil}
\widetilde\varphi\colon \widetilde{\mathcal E}_X \longrightarrow\widetilde{\mathcal L}_X(p)\otimes \pi_X^*\mathcal O_X(1) \text{\; with \;}  \pi_{X}\colon X\times C\to X,
\end{equation}
where \(\widetilde{\mathcal E}_X\) and \(\widetilde{\mathcal L}_X\) denote the restrictions of the  corresponding
universal bundles to \(X\times C\).
The natural $\mathrm{SL}_\chi$-action on $\mathbb P(\mathcal H)$
preserves $X$, and the bundles $\widetilde{\mathcal E}_X$,
$\widetilde{\mathcal L}_X$ and $\mathcal O_X(1)$ carry
the induced  $\mathrm{SL}_\chi$-linearizations.
 Since both the factorization and the strongly parabolic condition are invariant under rescaling, the $\TT$-action preserves $X$ and commutes with the $\mathrm{SL}_\chi$-action. With the conventions described above, the universal bundles $\widetilde{\mathcal E}_X$ and $\widetilde{\mathcal L}_X$ carry the trivial $\TT$-linearization, while the universal morphisms \eqref{eq:tauth} and \eqref{eq:tautphitil} have $\TT$-weight $-1$.

Let $\text{det}\colon X\to \operatorname{Jac}^l(C)$ be the morphism associated with $\det\widetilde{\mathcal{E}}_X$ and let $\mathcal{J}$ be the Poincar\'e bundle on $\Jac^l\times C$ normalized as above. Define a line bundle $\mathcal{D}$ on $X$ by $$\det\widetilde{\mathcal{E}}_X= (\text{det}\times \mathrm{id})^* \mathcal{J}\otimes \pi_X^*\mathcal{D}.$$ 
A simple calculation shows that the line bundle $\mathcal D$ has $\TT$-weight $0$.
Then, up to multiplication of all exponents by a common positive integer, the $\mathrm{SL}_\chi$-linearization on $X$ corresponding to \(\alpha\)-stability is
$$ L^{\!\alpha}\!:=\!\mathcal{D}^{\chi-\overline{\alpha}}\otimes \mathcal{O}_X(1)^{2\overline{\alpha}}\otimes\left(\tilde{\mathcal{G}}_2/\tilde{\mathcal{G}}_1\otimes\tilde{\mathcal{G}}_1^\vee\right)^{c_{\iio}-c_{\iit}},$$ where $  \overline{\alpha}:=\frac{2}{3}(\alpha-\alpha_{\mathrm{min}})=\frac{2\alpha-4g+4}{3}+l-c_{\iio}-c_{\iit}.$
The moduli space of $\alpha$-polystable parabolic triples is therefore obtained as  the GIT quotient\footnote{The center of $\mathrm{SL}_\chi$ acts trivially on $X$, and the chosen
polarization has trivial central character; hence both the action and the linearization descend to $\mathrm{PSL}_\chi$.}  $$N_l^\alpha({\hat c})=X\!\sslash^{L^{\!\alpha}}\!\mathrm{PSL}_\chi.$$
\begin{remark}\label{rem:GIT}
By the deformation theory of parabolic triples (see \cite[\S4.3 and \S12.4]{GPGM}), and since every nonzero morphism \(\phi\colon E\to L(p)\) is generically surjective, the parameter space $X$ is smooth. We also note that, although $X$ is not projective, for a sufficiently large choice of $n$, no semistable points with respect to the polarizations $L^\alpha$,
with $\alpha$ varying in a fixed bounded interval in $(\alpha_{\min},\infty)$, lie over the omitted points of the Quot scheme (cf. \cite[Proposition 6.5]{Beck}).
\end{remark}
We now describe some properties of the  moduli space $N_l^\alpha({\hat c})$ that will be needed later.
\begin{proposition}\label{prop:unive21}
Let $(c_1,c_2,c_3)\in\Delta$ be regular, and assume that \(\alpha\) is generic, so that $N^{\alpha, \mathrm{st}}_l({\hat c})=N^\alpha_l({\hat c})$. Then: 
\begin{enumerate}[label=(\roman*), leftmargin=1.5pc]
\item There exists a universal triple $\mathcal T = (\mathcal L,\mathcal E,\varphi)$, $\varphi\colon\mathcal E\to\mathcal L(p)$  on \(N_l^\alpha({\hat c})\times C\), together with the universal parabolic flag $\mathcal G_1\subset\mathcal G_2 =\mathcal E|_{N_l^\alpha({\hat c})\times\{p\}}.$ 
The universal triple is unique up to simultaneous tensoring of \(\mathcal E\) and \(\mathcal L\) by the pullback of a line bundle on \(N_l^\alpha({\hat c})\).

\item The torus \(\TT\) acts on triples by $t\cdot(L,E,\phi)=(L,E,t\phi)$, which induces the trivial action on \(N_l^\alpha({\hat c})\). The universal bundles admit \(\TT\)-linearizations of weights $1$ on $\mathcal{E}$ and $0$ on $\mathcal{L}$;  in particular, \(\varphi\) has \(\TT\)-weight \(-1\).
\end{enumerate}
\end{proposition}
\begin{proof}
For (i), note that the center of \(\mathrm{SL}_\chi\) acts with opposite characters on \(\widetilde{\mathcal E}_X\) and \(\mathcal O_X(-1)\), while it acts trivially on \(\widetilde{\mathcal L}_X\). Hence \(\widetilde{\mathcal E}_X\otimes\mathcal O_X(-1)\) and \(\widetilde{\mathcal L}_X\) descend, by Kempf's descent lemma (see, e.g., \cite[Theorem 2.3]{DrezetNarasimhan}), to the universal triple above. The same argument applies to  the universal parabolic flag.

Part (ii) follows from the construction. Indeed, \((L,E,t\phi)\simeq(L,E,\phi)\), so the induced action on the moduli space is trivial, while
\(\mathcal O_X(-1)\) has \(\TT\)-weight \(1\) and \(\widetilde{\mathcal E}_X,\widetilde{\mathcal L}_X\) have $\TT$-weight \(0\).
\end{proof}

\subsection{Type-(2,1) components as moduli spaces of stable parabolic triples}
We now describe explicitly the relation between the moduli spaces of stable parabolic triples introduced  above and the connected components of the fixed-point locus $\Hs_0(c)^\TT$.
\begin{proposition}\label{prop:tripletofix}
\label{prop:fixed-point-triples}
Let $c= (c_1,c_2,c_3)\in\Delta$ be regular; let $1\leq \iio<\iit\leq3$, denote by $\jjo$ the remaining element of $\{1,2,3\}$ and set $\hat{c}=(c_{\iio},c_{\iit})$. 
The connected components of the fixed-point locus $\Hs_0(c)^\TT$ of parabolic type $(\{\iio,\iit\},\{\jjo\})$ are isomorphic to the moduli spaces $N^{2g-2}_{l}({\hat c}),$ where $l\in\mathbb{Z}$ satisfies \({c}_{\iio}+{c}_{\iit}<l\leq g-2\) if  \(\jjo=1\) and \(c_2<0\), and 
\begin{equation}\label{eq:bound}
{c}_{\iio}+{c}_{\iit}<l\leq g-1
\end{equation}
in all other cases. 
\end{proposition}
\begin{proof}
Let $(V_1\oplus V_2, \Phi)\in\Hs_0(c)$  be a fixed point of parabolic type $(\{\iio,\iit\},\{\jjo\})$, where  $\Phi|_{V_2}=0$ and $\phi:=\Phi|_{V_1}$ is a nonzero strongly parabolic morphism $V_1\to V_2\otimes K(p)$. Let \(l = \deg(V_1)\), so that \(\deg(V_2)=-l\). Since \(V_2\) is \(\Phi\)-invariant, parabolic stability gives $$-l-c_{\jjo}=-l+{c}_{\iio}+{c}_{\iit}<0.$$

Now let $N:=\ker\phi\subset V_1,$ $Q:=V_1/N $ and write  $\deg(Q)=q.$ 
Let $c_a$, with $a\in\{\iio,\iit\}$ be the parabolic weight of $Q$. Since $N$ is $\Phi$-invariant, stability implies $$l-q-({c}_{\iio}+{c}_{\iit}-c_a)<0.$$
Moreover, $\phi$ induces a nonzero strongly parabolic morphism $Q\to V_2\otimes K(p)$, and hence 
$$q\leq -l+2g-1- \mathfrak{s}(a,\jjo), \text{\; where \;}  \mathfrak{s}(a,\jjo)=
\begin{cases}
0,& c_a>c_{\jjo},\\
1,& c_a<c_{\jjo}.
\end{cases}$$
Combining these inequalities, we obtain
$2l<2g-1+{c}_{\iio}+{c}_{\iit}-c_a- \mathfrak{s}(a,\jjo).$
Taking the weaker of the two bounds, as \(a\) varies in \(\{\iio,\iit\}\), gives \(l\leq g-1\), except when \(\jjo=1\) and \(c_2<0\), in which case \(l\leq g-2\). 

Finally, by  \cite[Proposition 4.6]{GPGM} the parabolic Higgs bundle
\((V_1\oplus V_2,\Phi)\) is stable if and only if the parabolic triple $(V_2\otimes K, V_1,\phi)$ is \((2g-2)\)-stable. Since
\(\deg(V_1)=l\), \(\deg(V_2\otimes K)=2g-2-l\), and
$\det(V_1)\otimes V_2\simeq\mathcal{O}$, this identifies the fixed locus with \(N_l^{2g-2}({\hat c})\), which, by Corollary
\ref{cor:max} below, is nonempty and connected.
\end{proof}

We note that, under the identification of a fixed-point component of parabolic type $(\{\iio,\iit\},\{\jjo\})$ with $N_l^{2g-2}(\hat c)$ given in Proposition \ref{prop:tripletofix}, the nonzero component of the Higgs field corresponds to the morphism $\phi$ of the triple. This identification is therefore $\TT$-equivariant, with $\TT$ acting on $N_l^{2g-2}(\hat c)$
by rescaling $\phi$, as described in Proposition \ref{prop:unive21}. Consequently, the equivariant classes entering the contribution of such a component to the equivariant Euler characteristic \eqref{Euler} can be expressed in terms of the $\TT$-linearized tautological bundles of Proposition \ref{prop:unive21}.

To this end, we make the following straightforward observation.
\begin{lemma}\label{lemma:Lon21}
The line bundle $\mathcal L(k;\lambda)$ restricts to  $N^{2g-2}_l({\hat c})\subset \Hs_0(c)$ as
\begin{equation}\label{eq:Lon21}
\mathrm{det}( \mathcal E_{p}\oplus\mathcal L_{p})^{k(1-g)}\otimes\mathrm{det}(\pi_!( \mathcal L\otimes K^\vee\oplus\mathcal E))^{-k} \otimes (\mathcal{G}_2/\mathcal{G}_1)^{\lambda_{\iio}}\otimes (\mathcal{G}_1)^{\lambda_{\iit}}  \otimes \mathcal L_p^{\lambda_{\jjo}},
\end{equation}
while the normalized universal bundle $U$ (see notation after Theorem \ref{thm:modsp}) restricts as
\begin{equation}\label{eq:Uon21}
(\mathcal L\otimes K^\vee\oplus\mathcal E)\otimes \mathcal{L}_p^{-1} \text{\; if \;} \iit<\jjo \text{\; and as \;} (\mathcal L\otimes K^\vee\oplus\mathcal E)\otimes \mathcal{G}_1^{-1} \text{\; if \;} \iit>\jjo.
\end{equation}
Here and below, we use the same notation $K$ for the canonical bundle on $C$ and for its pullback to $N^{2g-2}_l({\hat c})\times C$.
\end{lemma}

We will also need the following description of the normal bundle on $N_l^{2g-2}({\hat c})$.
\begin{lemma}\label{lemma:Non21}
 Let $\pi\colon N_l^{2g-2}({\hat c})\times C\to N_l^{2g-2}({\hat c})$ denote the natural projection. Then  the equivariant Chern character of the normal bundle to $N_l^{2g-2}({\hat c})$ in $\Hs_0(c)$ is given by 
\begin{multline}\label{eq:Normal21}
 -\ch_{\TT}\pi_! \operatorname{ParHom}\bigl(\mathcal E,\mathcal L \otimes K^\vee\bigr) -\ch_{\TT}\pi_!\operatorname{ParHom} \left(\mathcal L\otimes K^\vee,\mathcal E\right) +\\
t\ch_{\TT}\pi_!\left(\operatorname{SParEnd}(\mathcal E)\otimes K(p) \right) +
t\ch_\TT\pi_! \left(\operatorname{SParHom}\left(\mathcal L\otimes K^\vee,\mathcal E\right) \otimes K(p)\right).
\end{multline}
\end{lemma}
\begin{proof}
By Theorem \ref{thm:modsp}(iii), the restriction of the tangent bundle $T\Hs_0(c)$ to $N_l^{2g-2}(\hat c)$ is given by
\begin{equation}\label{eq:defR1comp}
\mathbb R^1_{\TT}\pi_* \left[\operatorname{ParEnd}_0 (\mathcal L\otimes K^\vee\oplus\mathcal E) \longrightarrow
t\operatorname{SParEnd}_0(\mathcal L\otimes K^\vee\oplus\mathcal E)\otimes K(p)\right].
\end{equation}
Decomposing this complex into $\mathbb T$-weight spaces, the weight-zero summand gives the tangent bundle $TN_l^{2g-2}({\hat c})$. The remaining summands have weights $-1$, $1$ and $2$, and are respectively  
\begin{multline*}
\mathbb R^1_{\TT}\pi_* \left[\operatorname{ParHom}(\mathcal E,\mathcal L\otimes K^\vee)\longrightarrow 0\right], \\
\mathbb R^1_{\TT}\pi_* \left[\operatorname{ParHom}(\mathcal L\otimes K^\vee,\mathcal E) \longrightarrow t\left(\operatorname{SParEnd}(\mathcal L\otimes K^\vee)\oplus\operatorname{SParEnd}(\mathcal E)\right)_0\otimes K(p)\right],\\
\text{\; and \;} \mathbb R^1_{\TT}\pi_* \left[ 0\longrightarrow t\operatorname{SParHom}(\mathcal L\otimes K^\vee,\mathcal E)\otimes K(p)\right].
\end{multline*}
By stability and parabolic Serre duality, the two-term deformation complex in \eqref{eq:defR1comp} has vanishing hypercohomology in degrees $0$ and $2$, and these vanishings hold in each $\TT$-weight. Thus, taking equivariant Chern characters of the three summands above, we obtain \eqref{eq:Normal21}.
\end{proof}

\subsection{Strategy for computing the contributions of  type-$(2,1)$ and type-$(1,2)$ components}\label{S:strategy}
Recall that we fixed a regular element $c=(c_1,c_2,c_3)\in\Delta$, chose indices $1\leq\iio<\iit\leq 3$, and denoted the remaining element of $\{1,2,3\}$ by $\jjo$. We now explain our strategy for computing the total contribution of the fixed-point components of parabolic type $(\{\iio,\iit\},\{\jjo\})$ to the equivariant Euler characteristic \eqref{Euler}.

We first compute the contribution  of a single component $N_l^{2g-2}(\hat c)$, where $\hat c=(c_{\iio},c_{\iit})$ and $l$ is an integer in the range specified in Proposition \ref{prop:tripletofix}.
Note that the formulas \eqref{eq:Lon21}, \eqref{eq:Uon21} and \eqref{eq:Normal21} involve only tautological bundles and therefore extend naturally to $N_l^\alpha({\hat c})$ for every generic value of $\alpha$.  
By a slight abuse of notation, we write $\mathcal N_l$ for the equivariant virtual $K$-theory class on $N_l^\alpha(\hat c)$ whose equivariant Chern character is given by the right-hand side of \eqref{eq:Normal21}, with the universal bundles now taken over $N_l^\alpha(\hat c)\times C$ and $\pi$ denoting the natural projection $\pi\colon N_l^\alpha(\hat c)\times C\to N_l^\alpha(\hat c)$.
Similarly, we denote by  $\mathcal L_{l}(k;\lambda)$ the line bundle on $N_l^\alpha({\hat c})$ defined by the same tautological expression as in \eqref{eq:Lon21} and by $U_{l}$ the normalized universal bundle on $N_l^\alpha({\hat c})\times C$ defined by \eqref{eq:Uon21}.  For a dominant weight \(\nu\) of \(\mathrm{GL}_3\), we denote by  \(U_{l,\nu}\) the bundle associated to \(U_l\) via the irreducible representation of \(\mathrm{GL}_3\) with highest weight \(\nu\).

For a generic value of $\alpha$ such that  \(N_l^{\alpha,\mathrm{st}}(\hat c)=N_l^\alpha(\hat c)\),
and for either choice $\mathcal V=\mathcal L_l(k;\lambda)$ or $ \mathcal V=\mathcal L_l(k;\lambda)\otimes\pi_!(U_{l,\nu}\otimes K^{1/2}),$  we set
\begin{equation}\label{eq:chypoly}
\chi_l^\alpha[\mathcal V] := \int_{N_l^\alpha(\hat c)} \frac{\ch_\TT(\mathcal V)}{e_\TT(\mathcal N_l^\vee)}\mathrm{Td}(N_l^\alpha(\hat c)),
\end{equation}
where $e_\TT$ denotes the equivariant Chern character of the $K$-theoretic Euler class, with the convention that $\chi_l^\alpha[\mathcal V]\!=\!0$ if $N_l^\alpha(\hat c)$ is empty. With this notation, $\chi_l^{2g-2}[\mathcal V]$ is precisely the contribution of $N_l^{2g-2}(\hat c)$ to the equivariant Euler characteristic \eqref{Euler}.

We compute this contribution by varying the stability parameter \(\alpha\). More precisely, in \S\ref{S:walls}, following \cite{GPGM}, we study how the moduli spaces \(N_l^\alpha(\hat c)\) depend on \(\alpha\). We show that the stability condition can change only at finitely many \emph{critical values}, which divide the interval \((\alpha_{\min},\alpha_{\mathrm{large}})\) from Theorem \ref{thm:moduli21}(iii) into \emph{chambers} consisting of generic values. The moduli space remains unchanged as \(\alpha\) varies within a chamber, and therefore
\[
\chi_l^\alpha[\mathcal V]=\chi_l^{\alpha'}[\mathcal V]
\]
whenever no critical value lies between \(\alpha\) and \(\alpha'\).
Now let \(\alpha_0\) be a critical value, and choose generic values $\alpha_-<\alpha_0<\alpha_+$ in the two adjacent chambers. We define the corresponding \emph{wall-crossing term} by
\begin{equation}\label{eq:wcgeneral}
\chi_l^{\alpha_-}[\mathcal{V}]- \chi_l^{\alpha_+}[\mathcal{V}].
\end{equation}
We compute these wall-crossing terms in \S\ref{S:masterspace}-\S\ref{S:compWCR} using the method of \cite{OTASz}.

Next, in Proposition \ref{prop:alphacryt}  we show that the moduli space  \(N_l^\alpha({\hat c})\) is empty for \(\alpha>\alpha_{\mathrm{max}}\), where  \(\alpha_\mathrm{max}\) is itself a critical value. Thus, the wall-crossing term \eqref{eq:wcgeneral} at \(\alpha_{\mathrm{max}}\) gives the value of \(\chi_l^\alpha[\mathcal{V}]\) in the chamber right below $\alpha_{\mathrm{max}}$. Adding the wall-crossing terms at the remaining critical values crossed as \(\alpha\) decreases to \(2g-2\), we recover \(\chi_l^{2g-2}[\mathcal V]\) in \S\ref{S:type21}.

Finally, summing the contributions \(\chi_l^{2g-2}[\mathcal V]\) over the values of $l$ specified in Proposition \ref{prop:tripletofix} gives the total contribution of the fixed-point components of parabolic type $(\{\iio, \iit\},\{\jjo\})$ to the equivariant Euler characteristic \eqref{Euler}. The contribution of the components of parabolic  type $(\{\jjo\}, \{\iio, \iit\})$ is then obtained by the duality described in \S\ref{S:type12}.

\section{Contributions of type-(2,1) and type-(1,2) components}\label{S:21and12}
We now compute the contribution of type-\((2,1)\) components to the equivariant Euler characteristics \eqref{Euler}, following the strategy described in \S\ref{S:strategy}. As in the previous section, we fix $1\leq\iio<\iit\leq 3$, denote the remaining element of $\{1,2,3\}$ by $\jjo$, and take $l$ in the range specified in Proposition \ref{prop:tripletofix}.
\subsection{Walls and chambers}\label{S:walls}
We begin by studying how the moduli space $N_l^\alpha(\hat c)$ varies with the stability parameter $\alpha$. Note that the stability of a triple can change only at values of $\alpha$  for which there exists an \(\alpha\)-semistable triple \(T=(L,E,\phi)\) admitting  a proper subtriple $T'=(L',E',\phi')\subset T$ whose $\alpha$-slope is equal to that of $T$. Writing \(r_1'=\operatorname{rk}(L')\) and \(r_2'=\operatorname{rk}(E')\),  this  slope condition becomes
\begin{equation}\label{sigmacrit}
\frac{\operatorname{pardeg}(L'\oplus E')+\alpha r_2'}{r_1'+r_2'} = \frac{2g-2+2\alpha}{3}.
\end{equation}
We say that a value of \(\alpha\) arising in this way is \emph{critical} and call the connected components of the complement of the set of critical values \emph{chambers}. 
These definitions imply the following statements.

\begin{proposition}\label{prop:alphacryt}
\begin{enumerate}[label=(\roman*), leftmargin=1.5pc]
\item If \(\alpha\) and \(\alpha'\) lie in the same chamber, i.e. between two consecutive critical values, then \(\alpha\)-stability and \(\alpha'\)\!-stability conditions are equivalent. In particular,
$N_l^\alpha({\hat c})\simeq N_l^{\alpha'}({\hat c}).$
\item Let $\alpha_{\max}$ be the maximal value of $\alpha$ satisfying $N_l^\alpha({\hat c})\neq\varnothing.$ Then \(\alpha_{\max}\) is critical. 
\end{enumerate}
\end{proposition}

We now describe the set of critical values of $\alpha$ explicitly. A simple calculation shows that the proper subtriples $T'\subset T$ satisfying \eqref{sigmacrit} can be of one of the following three types:
\begin{itemize}[leftmargin=5mm]
\item If \(r_1'=0\) and \(r_2'=1\), then \(T'\) is of the form \(M\to0\), where \(M\subset E\) is a line subbundle of degree \(m\) with parabolic weight \(c_\epsilon\), for some \(\epsilon\in\{\iio,\iit\}\). The quotient triple then has the form \(F\to L(p)\), where \(F=E/M\) is a parabolic line bundle of degree \(l-m\) and  weight \(c_{\bar\epsilon}\), with \(\{\epsilon,\bar\epsilon\}=\{\iio,\iit\}\):
\[
\begin{tikzcd}
0 \arrow[r] & M \arrow[r] \arrow[d] &
E \arrow[r] \arrow[d,"\phi"] &
F \arrow[r] \arrow[d,"\phi''"] & 0\\
& 0 \arrow[r] & L(p) \arrow[r] & L(p) \arrow[r] & 0.
\end{tikzcd}
\]
Equation \eqref{sigmacrit} then  gives the critical value
\begin{equation}\label{sigmac}
\alpha(m,\epsilon):=2g-2-3m+3c_\epsilon.
\end{equation}

\item If \(r_1'=r_2'=1\), then \(T'\) is of the form \(F\to L(p)\), while the quotient triple is  \(M\to0\), with \(F\) and \(M\) as above. The corresponding critical value is again given by \eqref{sigmac}.

\item Finally, if \(r_1'=1\) and \(r_2'=0\), then \(T\) is an extension of \(E\to0\) by \(0\to L(p)\). The corresponding critical value given by \eqref{sigmacrit} is $\alpha_{\min} = 2g-2-\frac{3l}{2}-\frac{3}{2}c_{\jjo}.$ 
\end{itemize}

Now given \(\epsilon\in\{\iio,\iit\}\), we denote by \(\bar\epsilon\)  the other element of \(\{\iio,\iit\}\). Define
\begin{equation}\label{defS}
S^{\epsilon}_1= 
\begin{cases} p, \text{\; if \; } c_{\jjo}\!>c_{\bar{\epsilon}}, \\
\varnothing,  \text{\; otherwise;}  
\end{cases} 
S^{\epsilon}_2= \begin{cases} p, \text{\; if \; } c_{\jjo}\!>c_{\epsilon}, \\
\varnothing,  \text{\; otherwise; }  
\end{cases}
S_3^{\epsilon}= \begin{cases} p, \text{\; if \; } c_\epsilon\!>c_{\bar{\epsilon}}, \\
\varnothing,  \text{\; otherwise;} 
\end{cases}
\end{equation}
and set $ s_i^\epsilon=|S_i^\epsilon|,$ for $i=1,2,3$. With this notation,  Theorem \ref{thm:moduli21}(iii) implies the following statement.
\begin{lemma}
The critical values satisfying \(\alpha\geq 2g-2\) are parametrized by pairs \((m,\epsilon)\), where \(\epsilon\in\{\iio,\iit\}\) and $-2g+1+2l+s_1^{\!\epsilon} \leq m < c_\epsilon.$ They are given by $\alpha(m,\epsilon)=2g-2-3m+3c_\epsilon.$
\end{lemma}
Combining this with Proposition \ref{prop:alphacryt}(ii) and using the description of the flip loci for parabolic triples with fixed determinant in \cite[\S6.1 and \S12.4]{GPGM}, we obtain the following corollary.
\begin{corollary}\label{cor:max}
The  moduli space $N_l^\alpha({\hat c})$ is nonempty if and only if $\alpha_{\mathrm{min}} \leq \alpha\leq \alpha_\mathrm{max}$, where 
$$\alpha_\mathrm{min} = 2g-2-\frac{3l}{2}-\frac{3{c}_{\jjo}}{2} \text {\; and \;} \alpha_{\max} = \begin{cases} 
8g-8-6l+3c_2,  \text{\ if \ } \jjo=1, \\
8g-5-6l+3c_3,  \text{\ if \ } \jjo=2, \\
8g-5-6l+3c_1,  \text{\ if \ } \jjo=3. 
\end{cases}$$ 
Moreover, $N_l^\alpha(\hat c)$ is connected for every non-critical $\alpha_{\mathrm{min}}<\alpha<\alpha_{\mathrm{max}}$.
\end{corollary}

\subsection{Wall-crossing via the master space}\label{S:masterspace}
 We now fix $\epsilon\in\{\iio,\iit\}$ and an integer $m$ satisfying  $-2g+1+2l+s_1^{\!\epsilon}\leq m<c_\epsilon$, and recall the notation $\alpha(m,\epsilon)$ for the corresponding critical value defined in \eqref{sigmac}.  Our next goal is to compute the change in the Euler characteristics \eqref{eq:wcgeneral} as the stability parameter \(\alpha\) crosses the critical value $\alpha(m,\epsilon)$, following the strategy of \cite{OTASz}.

We choose  $\mathfrak{d}\in\mathbb Q_{>0}$ sufficiently small so that $\alpha_-:=\alpha(m,\epsilon)-\mathfrak{d}$ and $\alpha_+:=\alpha(m,\epsilon)+\mathfrak{d}$ lie in the two adjacent chambers. Recall from \S\ref{S:GIT} that the moduli space $N_l^\alpha({\hat c})$ is obtained as a GIT quotient $X\sslash^{L^\alpha}\!G$ with $G=\mathrm{PSL}_\chi$, and that the dependence of the GIT polarization $L^\alpha$ on
$\alpha$ is carried by the $\TT\times G$-equivariant line bundle 
\begin{equation}\label{eq:Lvec}
\vec L:=\mathcal D^{-1}\otimes\mathcal O_X(1)^{2}. 
\end{equation}
Thus, after replacing $L^{\alpha_-}$ and $L^{\alpha_+}$ by suitable common positive tensor powers, we may choose $n>0$ such that $L^{\alpha_+} \simeq L^{\alpha_-}\otimes\vec L^{n}.$

Let $X^0\subset X$ be the locus of points which are semistable with respect to $L^{\alpha(m,\epsilon)}$, but unstable with respect to $L^{\alpha_-}$ and $L^{\alpha_+}$.

\begin{lemma}\label{lem:masterassumptions}
The locus $X^0$ is smooth. The generic stabilizer of the $G$-action on
$X$ is trivial, while for every $x\in X^0$ the stabilizer $G_x$ acts freely on $\vec L_x\setminus\{0\}$,  yielding an identification $G_x \cong \mathbb{C}^*$.
\end{lemma}
\begin{proof}
We note that point \(x\in X^0\) represents a split triple $T=T'\oplus T"$, where both $T'$ and $T"$ are $\alpha(m,\epsilon)$-stable; since both factors are parametrized by smooth varieties, \(X^0\) is smooth.

We now study the action of the stabilizer $\widetilde{G}_x\subset{\mathrm{SL}_\chi}$ of point $x\in{X}$ on the fiber $\vec L_x\setminus{0}$.               
\begin{itemize}[leftmargin=5mm]
\item For a general point $x\in X$, the corresponding triple is stable and hence has only scalar automorphisms. Thus \(\widetilde G_x\) is the center \(\mathbb Z_\chi\subset \mathrm{SL}_\chi\), which acts trivially on \(\vec L_x\setminus\{0\}\).
\item For \(x\in X^0\), the splitting of the corresponding triple \(T=T'\oplus T''\) and the stability of the two factors give $\mathrm{Aut}(T)\simeq\mathbb C^*\times\mathbb C^*$. 
Writing an element of \(\widetilde G_x\) as \((t_1,t_2)\), its action on \(\vec L_x\) is given by the multiplication by \(t_1^{-1}t_2\). If this action is trivial, then \(t_1=t_2\), and  since \((t_1,t_2)\in \mathrm{SL}_\chi\), we obtain \(t_1^\chi=1\). Hence the stabilizer in \(\widetilde G_x\) of any point of \(\vec L_x\setminus\{0\}\) is precisely \(\mathbb Z_\chi\).
\end{itemize}
Thus the generic stabilizer of the \(G=\mathrm{PSL}_\chi\)-action on \(X\) is trivial, while for \(x\in X^0\) the stabilizer $G_x\simeq \widetilde G_x/\mathbb Z_\chi$
acts freely on \(\vec L_x\setminus\{0\}\), and this action induces an isomorphism $G_x\simeq\mathbb C^*.$
\end{proof}

Together with Remark \ref{rem:GIT}, Lemma \ref{lem:masterassumptions} allows us to apply the wall-crossing argument of \cite{OTASz}.  In preparation for this,
we now briefly recall the master-space construction from \cite[\S5]{OTASz} in the form needed for our computation.

Consider the projective bundle $Y := \mathbb P(\mathcal O_X \oplus \vec L)$  and introduce the \emph{master-space torus} $T$ acting on its fibers by $(1,s^{-1}).$
 The $\TT\times G$-linearization of $\vec L$ induces a $\TT\times G$-action on $Y$ commuting with $T$, and the  \emph{master space} is defined as $$Z:=Y\sslash^{L^{\alpha_+}\otimes\mathcal O_Y(n)}G.$$ 
The torus $T\times \TT$ acts naturally on $Z$, and the $T$-fixed locus decomposes as 
$$Z^T= N_l^{\alpha_-}({\hat c}) \sqcup N_l^{\alpha_+}({\hat c})\sqcup  Z^0.$$  
The components $N_l^{\alpha_-}({\hat c})$ and $N_l^{\alpha_+}({\hat c})$ arise from the zero and infinity sections of $Y$; either of them is absent when $\alpha_-<\alpha_{\min}$ or $\alpha_+>\alpha_{\max}$.
To describe the third component $Z^0$, we recall some additional notation from  \cite[\S5.2]{OTASz}. Let $\vec L^\circ:=\vec L\setminus\{0\}$, write $q_\circ$ for the
restriction of the projective bundle map $q\colon Y\to X$ to $\vec L^\circ$, and denote by $\overline q$ the quotient map to $Z$. Then \cite[Proposition 5.2]{OTASz}
gives  an embedding $\iota_0\colon X^0/G\hookrightarrow Z$ with image $Z^0:=\overline q\bigl(q_\circ^{-1}(X^0)\bigr)$. 
\begin{lemma}\label{lem:Z0descr}
There is an isomorphism
$$ Z^0\simeq\mathrm{Jac}^{2g-2-l}(C)\times\Sym^{2g-1-s_1^{\epsilon}-2l+m}(C).$$
\end{lemma}
 \begin{proof}
Since $\iota_0$ identifies $Z^0$ with $X^0/G$, the points of $Z^0$ correspond to the closed $G$-orbits in $X^0$. These represent triples that split as
$T=T'\oplus T''$, where $T'=(0,M,0)$ and $T''=(L,F,\phi)$. Here $M$, $F$ and $L$ are
parabolic line bundles of degrees $m$, $l-m$ and $2g-2-l$, with parabolic weights $c_\epsilon$, $c_{\bar\epsilon}$ and $c_{\jjo}$, respectively, such that
$\det(M\oplus F\oplus L)\cong K$, and $\phi\colon F\to L(p)$ is a nonzero strongly parabolic morphism (see \S\ref{S:walls}).  Since such a triple is determined by $L$ and the divisor of $\phi$,
we obtain  the claimed isomorphism.
 \end{proof}
 
We note that, since the $\TT$-action on $Z$ commutes with $T$, the $T$-fixed  component $Z^0$ is preserved by $\TT$.  Hence we can use the following construction from \cite{OTASz}.

\begin{construction}\label{descconstruction} Given a $G \times \TT$-equivariant vector bundle $\widetilde V$ on $X$, we construct a $\TT\times T$-equivariant bundle $\zeta(\widetilde V)$ on $Z$ as follows. We endow the pullback $q^*{\widetilde{V}}$ with the trivial action of $T$ and the action of $G \times \TT$ pulled back from $X$, and then descend $q^*\widetilde{V}$ to $Z$.  
Restricted to $X^0$, this functor sends $G\times\TT$-equivariant bundles on $X^0$ to $T\times\TT$-equivariant bundles on $Z^0$. 
\end{construction}
The following lemma describes how equivariant weights change under this functor. 
\begin{lemma}\label{lem:Lvectriv}
 Recall that  for $x\in X^0$, we identified $G_x\simeq\mathbb C^*$ through its action on $\vec L_x\setminus\{0\}$.
\begin{enumerate}[label=(\roman*), leftmargin=1.5pc]
\item If a $G \times \TT$-equivariant vector bundle $\widetilde{V} \to X^0$ has $G_x$-weight $a$ and $\TT$-weight $b$, then $\zeta|_{X^0}(\widetilde{V}) \to Z^0$ has $T\times \TT$-weight $(a, b+2a)$.
\item In particular, $\zeta(\vec{L})|_{Z^0}$ is the trivial line bundle on $Z^0$ carrying $T\times \TT$-weight $(1,0)$.
\end{enumerate}
\end{lemma} 
\begin{proof}
Part (i) follows from the definition of the functor $\zeta|_{X^0}$, the identification $G_x\simeq\mathbb C^*$ through its action on $\vec L_x\setminus\{0\}$, and the fact that $T$ and $\TT$ act on $\vec L$ with weights \(-1\) and \(-2\), respectively.

Applying (i) to \(\vec L\) gives the \(T\times\TT\)-weight \((1,0)\). Since \(q^*_\circ\vec L\) admits a tautological nowhere-vanishing \(G\times\TT\)-equivariant section, the line bundle $\zeta(\vec{L})|_{Z^0}$ is trivial.
\end{proof}
To write down the wall-crossing difference \eqref{eq:wcgeneral}, we introduce some extra notation.

\noindent\textbf{Notation:}  For a bundle $V$ on $N_l^\alpha(\hat{c})$ obtained by descent  from a $G\times \TT$-equivariant bundle $\widetilde V$ on $X$, we write  ${V}|_{Z^0} := \zeta(\widetilde{V})|_{Z^0}$.  
Then, expanding the inverse equivariant $K$-theory Euler class of $\mathcal N_l^\vee$ (see \S\ref{S:strategy}) in powers of $t^{-1}$, we write
$$e_{\!\TT}\!\left(\mathcal N_l^\vee\right)^{-1}=\sum_n \ch(E_n) t^n,$$ where each $E_n$ is a $K$-theory class on $N^\alpha_l({\hat c})$, and set  $$e_{T\times\TT}\left({\mathcal N}_l^\vee\right)\!^{-1}\big|_{Z^0}:=\sum_n \ch_{T\times\TT}(E_n|_{Z^0}) t^n.$$

We are now ready to compute the wall-crossing difference \eqref{eq:wcgeneral}. Let \(\alpha_-<\alpha(m,\epsilon)<\alpha_+\) lie in the two chambers separated by the critical value \(\alpha(m,\epsilon)\). We apply \cite[Theorem 5.7]{OTASz} to each class $\mathcal V\otimes E_n$, using localization with respect to the torus $T$ with equivariant parameter $s$.
Multiplying the resulting identities by $t^n$ and summing over $n$, we obtain, in the notation introduced
above,
\begin{equation}\label{intwc0}
\chi_l^{\alpha_-}[\mathcal{V}] - \chi_l^{\alpha_+}[\mathcal{V}]=\underset{s=0, \infty}{\res}\int_{Z^0}\frac{\ch_{T\times\TT\!}\left(\mathcal{V}\big|_{Z^0}\right)}{e_{T\times\TT\!}\left({\mathcal N}_l^\vee\right)\!\!\big|_{Z^0}
e_{T\times\TT}\bigl(N^\vee_{Z^0/Z}\bigr)}\mathrm{Td}(Z^0)\frac{\mathrm{d}s}{s}.
\end{equation}
\begin{remark}\label{rem:extremalwall}
We note that the same formula applies at an extremal critical value $\alpha_\mathrm{max}$. More precisely, if  $N_l^{\alpha_+}({\hat c})$ is empty, then \eqref{intwc0} computes  the equivariant Euler characteristic  $\chi_l^{\alpha_-}[\mathcal V]$.
\end{remark}

To rewrite \eqref{intwc0} in a form more convenient for the calculations below, we recall from  \cite[Lemma 5.3]{OTASz} the following description of the normal bundle \(N_{Z^0/Z}\).
\begin{lemma}\label{lem:Normal}
For any \(x\in X^0\), the fiber of \(N_{Z^0/Z}\) at \(\iota_0(x)\in Z^0\) is isomorphic, as a \(T\times\TT\)-vector space, to $\vec{L}^\circ_x\!\times_{G_x}\!N_xX^0$.
Here \(N_xX^0\) denotes the normal space to \(X^0\) at \(x\), and, writing \((s,t)\in T\times\TT\) for the standard equivariant parameters, the action on the first factor is given by multiplication by \(s^{-1}t^{-2}\).
\end{lemma}
By Lemma \ref{lem:Normal} (see also Proposition \ref{lem:normalZ0} below), the equivariant weights of $N_{Z^0/Z}$ are $st^2$ and $(st^2)^{-1}$.
Hence, for each $n$, the meromorphic form in \eqref{intwc0} corresponding to the class $\mathcal V\otimes E_n$ has poles only at $s=0$, $\infty$, and $t^{-2}$.
Moreover,  Lemma \ref{lem:Lvectriv}(i) implies  that taking the residue at $s=t^{-2}$ preserves the
original $\TT$-grading and hence commutes with summation over $n$.
Thus, applying the residue theorem for each \(n\) and then summing,  we arrive at the following result.

\begin{proposition}\label{prop:ourint}
Fix \(\epsilon\in\{\iio,\iit\}\), an integer \(m\) satisfying $-2g+1+2l+s_1^{\!\epsilon}\le m<c_\epsilon,$ and let $\alpha(m,\epsilon)$ be the critical value defined in \eqref{sigmac}. Then for \(\alpha_-<\alpha(m,\epsilon)<\alpha_+\) lying  in the two adjacent chambers, for either choice $\mathcal V=\mathcal L_l(k;\lambda)$ or $ \mathcal V=\mathcal L_l(k;\lambda)\otimes\pi_!(U_{l,\nu}\otimes K^{1/2})$, 
\begin{equation}\label{intwc1}
\chi_l^{\alpha_-}[\mathcal{V}] - \chi_l^{\alpha_+}[\mathcal{V}]=-\underset{s=t^{-2}}{\res}\int_{Z^0}\frac{\ch_{T\!\times\!\TT}\left(\mathcal{V}\big|_{Z^0}\right)}{e_{T\!\times\!\TT}\left({\mathcal N}_l^\vee\right)\!\!\big|_{Z^0}
e_{T\!\times\!\TT}\bigl(N^\vee_{Z^0/Z}\bigr)}\mathrm{Td}(Z^0)\frac{\mathrm{d}s}{s}, 
\end{equation}
where the bundles $\mathcal L_l(k;\lambda)$, $U_{l,\nu}$ and $\mathcal{N}_l$ are defined in \S\ref{S:strategy}.
\end{proposition}
In the following subsections, we study the vector bundles appearing in \eqref{intwc1} and compute this integral explicitly.

\subsection{Vector bundles on $Z^0$}
We begin by introducing the notation for the universal bundles on $Z^0\times C$. Set $N:=2g-1-2l+m-s_1^{\!\epsilon}$ and recall that  $Z^0\simeq \Jac^{2g-2-l}(C) \times \Sym^N(C).$  
Let  $\mathcal{J}$ denote the Poincar\'e bundle on $\Jac^{2g-2-l}(C) \times C$, normalized so that its restriction to  $\Jac^{2g-2-l}\times\{p\}$ is trivial, and let $\mathcal{O}(\Delta)$ be the line bundle  associated to the universal divisor $\Delta\subset \Sym^{N}(C) \times C$; by abuse of notation, we write $\mathcal{J}$ and $\mathcal{O}(\Delta)$ also for their pullbacks to $Z^0\times C$. Define
$$\mathcal{F} := \mathcal{O}(-\Delta) \otimes \mathcal{J}\otimes\mathcal{O}(p-S_1^{\!\epsilon}) \text{\; and \;} \mathcal{M}:= \mathcal{J}^{-2}\otimes\mathcal{O}(\Delta)\otimes K(S_1^{\!\epsilon}-p),$$ 
where, as usual, we use the same notation for line bundles on $C$ and their pullbacks to $Z^0\times C$.  

\begin{lemma}\label{MFJweight}
The bundles $\mathcal F^\vee \otimes \mathcal J$, $\mathcal M^\vee \otimes \mathcal F$ and $\mathcal M^\vee\otimes \mathcal J$ have $T\times\TT$ weights $(0,-1)$, $(1,2)$ and $(1,1)$, respectively. 
\end{lemma}
\begin{proof}
First, we note that $\zeta(\vec{L})|_{Z^0} = \mathcal{F}^\vee_p\otimes\mathcal{M}^\vee_p\otimes\mathcal{J}^{2}_p$, where, given a bundle $V$ on $Z^0\times C$, we denote by $V_p$ its restriction to $Z^0\times\{p\}$. Thus, by Lemma \ref{lem:Lvectriv}(ii)  the $T\times\TT$ weight of $(\mathcal{M}^\vee_p\otimes\mathcal{F}_p)\otimes(\mathcal{F}^\vee_p\otimes\mathcal{J}_p)^{2}$ is equal to $(1,0)$. 

On the other hand, Lemma \ref{lem:Lvectriv}(i) implies that the $T\times\TT$-weight of $\mathcal{F}^\vee_p\otimes\mathcal{J}_p$ is equal to $(0,-1)$, and thus the $T\times\TT$-weight of $\mathcal{M}^\vee_p\otimes\mathcal{F}_p$ is $(1,2)$. Finally, since $\mathcal{M}^\vee\otimes\mathcal{J} = (\mathcal{M}^\vee\otimes \mathcal{F})\otimes (\mathcal{F}^\vee\otimes\mathcal{J})$, the result follows.
\end{proof}
The normal bundle \(N_{Z^0/Z}\) then admits the following description.

\begin{proposition}\label{lem:normalZ0}
Let $\pi\colon Z^0\times C\to Z^0$ be the natural projection, let $\mathcal M,$  $\mathcal F$, $\mathcal J$ be the line bundles on $Z^0\times C$ introduced above, and let $\varphi\colon \mathcal{F}\to \mathcal{J}(p)$ be the universal strongly parabolic morphism.  Consider the two-term equivariant complexes
$$\mathcal N^+_\bullet :=[\mathrm{ParHom}(\mathcal M,\mathcal F)\to t\,\mathrm{SParHom}(\mathcal M,\mathcal J(p))], \;\;\;  \mathcal N^-_\bullet:=[\mathrm{ParHom}(\mathcal F,\mathcal M)\to 0],$$  placed in degrees \(0\) and \(1\), with the differential induced by composition with \(\varphi\). 
Then, as a \(T\times\TT\)-equivariant vector bundle, $N_{Z^0/Z}\cong \mathbb R^1\pi_*(\mathcal N^+_\bullet)\;\oplus\;\mathbb R^1\pi_*(\mathcal N^-_\bullet).$
\end{proposition}
\begin{proof}
Let $x\in X^0$. By Lemma \ref{lem:Normal}, the fiber of $N_{Z^0/Z}$ at $\iota_0(x)$ is isomorphic to $\vec L_x^\circ\!\times_{G_x}\!N_xX^0,$ where $\vec L_x^\circ= \vec L_x\setminus\{0\}$, and $T\times\TT$ acts on $\vec L_x^\circ$ with weight $(-1,-2)$.

To describe $N_xX^0$, note that the universal triple on $X^0\times C$ splits as $(0,\widetilde{\mathcal M},0)\oplus (\widetilde{\mathcal J},\widetilde{\mathcal F},\widetilde\varphi)$, where the restrictions of $\widetilde{\mathcal M}$, $\widetilde{\mathcal F}$ and $\widetilde{\mathcal J}$ to $\{x\}\times C$ are the parabolic line bundles $M$, $F$ and $L$  introduced in the proof of Lemma \ref{lem:Z0descr}.  By the deformation theory of parabolic triples (see \cite[\S4.3 and \S12.4]{GPGM}), the normal space to $X^0\subset X$ at $x$ is given by
\begin{multline*}
\;\;\;\;\;\;\;\;\;\;\;\;\;\;\; \;\;\;\;\;\;\;\;\;\; N_xX^0=\mathbb{H}^1\bigl(\widetilde{\mathcal{N}}^+_{\bullet}\big|_{\{x\}\times C}\bigr)\oplus \mathbb{H}^1\bigl(\widetilde{\mathcal{N}}^-_{\bullet}\big|_{\{x\}\times C}\bigr), \text{\; where \;\;\;\;\;\;\;\;\;\;\;} \\ 
\widetilde{\mathcal N}^+_{\bullet}=[\mathrm{ParHom}(\widetilde{\mathcal M},\widetilde{\mathcal F})\to t\,\mathrm{SParHom}(\widetilde{\mathcal M},\widetilde{\mathcal J}(p))] \text{\; and \;}   \widetilde{\mathcal N}^-_{\bullet}=[\mathrm{ParHom}(\widetilde{\mathcal F},\widetilde{\mathcal M})\to 0].
\end{multline*}

Given a bundle \(\widetilde V\) on \(X^0\times C\), we denote its restriction to \(X^0\times\{p\}\) by \(\widetilde V_p\). With this notation, \eqref{eq:Lvec} yields
$$\vec L\big|_{X^0}=\widetilde{\mathcal{M}}_p^{-1}\otimes\widetilde{\mathcal F}_p^{-1}\otimes\widetilde{\mathcal{J}}_p^{2}.$$
Since the stabilizer $G_x$ acts with the same character on $F$ and $L$, while $\widetilde{\mathcal M}$, $\widetilde{\mathcal F}$ and $\widetilde{\mathcal J}$
have $\TT$-weights $1$, $1$ and $0$, respectively (see the proof of Proposition \ref{prop:unive21}), the two summands of $N_xX^0$ have $G_x\times\TT$-weights
$(1,0)$ and $(-1,0)$, respectively. By Lemma \ref{lem:Lvectriv}(i), the corresponding summands of $N_{Z^0/Z}$ therefore have $T\times\TT$-weights $(1,2)$ and $(-1,-2)$.

To identify the underlying bundles, we note that
\begin{multline*}
 \vec L_x^\circ\!\times_{G_x}\!\mathbb H^1\bigl(\widetilde{\mathcal N}^{+}_{\bullet}\big|_{\{x\}\times C}\bigr)\simeq\vec L_x^{-1}\otimes \mathbb H^1\bigl(\widetilde{\mathcal N}^{+}_{\bullet}\big|_{\{x\}\times C}\bigr) \simeq \mathbb H^1 (\vec L_x^{-1}\otimes \widetilde{\mathcal N}^{+}_{\bullet}\big|_{\{x\}\times C})  \text{\; and \;} \\ \vec L_x^\circ\!\times_{G_x}\!\mathbb H^1\bigl(\widetilde{\mathcal N}^{-}_{\bullet}\big|_{\{x\}\times C}\bigr) \simeq\vec L_x\otimes \mathbb H^1\bigl(\widetilde{\mathcal N}^{-}_{\bullet}\big|_{\{x\}\times C}\bigr) \simeq \mathbb H^1\bigl(\vec L_x\otimes \widetilde{\mathcal N}^{-}_{\bullet}\big|_{\{x\}\times C} \bigr).
\end{multline*}
By our choice of normalization, the bundles
\begin{equation}\label{eq:MFLdesc}
\widetilde{\mathcal F}\otimes\widetilde{\mathcal J}_p^{-1},\;\;\;\;
 \widetilde{\mathcal M}\otimes \vec{L}\big|_{X^0}\otimes\widetilde{\mathcal J}_p^{-1} \text{\; and \;}
\widetilde{\mathcal J}\otimes\widetilde{\mathcal J}_p^{-1}
\end{equation}
descend, respectively, to  the underlying line bundles $\mathcal F$, $\mathcal M$, and $\mathcal J$ on $Z^0\times C$. Thus the complexes
$$\vec L^{-1}|_{X^0}\otimes\widetilde{\mathcal N}^+_\bullet \text{\; and \;} \vec L|_{X^0}\otimes\widetilde{\mathcal N}^-_\bullet$$ 
descend, as ordinary complexes, to the underlying complexes of $\mathcal N^+_\bullet$ and $\mathcal N^-_\bullet$, respectively. 
By Lemma \ref{MFJweight}, the equivariant weights of
$\mathcal N^+_\bullet$ and $\mathcal N^-_\bullet$ agree with those determined
above. Taking relative hypercohomology completes
the proof.
\end{proof}

\subsection{Calculation of the characteristic classes of vector bundles on $Z^0$}
We now fix the notation for the cohomology classes  on  $Z^0$  used in our computations below. 
Define $\eta \in H^2(\Sym^N(C))$, $\xi_i \in H^1(\Sym^{N}(C)),$ and $\zeta_i \in H^1(\Jac^{2g-2-l}(C))$, for $i = 1, \ldots, 2g$, via the K\"unneth decompositions
\begin{multline*}
  \;\;\;\;\;\;\;\;\;\;\;\;\;\;\;\;\; c_1(\mathcal{O}(\Delta)) = \eta\otimes 1 +   \sum_{i=1}^{2g} \xi_i\otimes e_i +N\otimes \omega\in H^2(\Sym^N(C)\times C) \text{\; and \;} \\
c_1(\mathcal{J}) =   \sum_{i=1}^{2g} \zeta_i\otimes e_i + (2g-2-l)\otimes \omega \in H^2(\mathrm{Jac}^{2g-2-l}(C)\times C).  \;\;\;\;\;\;\;\;\;\;\;\; \;\;\;\;\;
\end{multline*}
Here, as in \S\ref{S:vbunonF}, $\omega\in H^2(C)$ is the class Poincar\'e dual to a point, and  $e_1,\dots,e_{2g}$ is a symplectic basis of $H^1(C)$.
We use the same notation for the pullbacks of \(\eta\), \(\xi_i\), and \(\zeta_i\) to \(Z^0\). Thus, writing $t=e^u$, $s=e^v$ and using Lemma \ref{MFJweight}, we obtain
\begin{multline}\label{c1MFL} 
c_1^{T\!\times\!\TT}(\mathcal{F}^\vee\otimes \mathcal{J}\otimes K^\vee)= \eta-u +   \sum_{i=1}^{2g} \xi_i\otimes e_i +(m-2l)\otimes\omega \text{\; and \;} \\
c_1^{T\!\times\!\TT}(\mathcal{M}^\vee\otimes \mathcal{F})= -2\eta +2u+v +\sum_{i=1}^{2g}(3\zeta_i-2\xi_i)\otimes e_i +(l-2m)\otimes \omega.
\end{multline} 
We also define $\theta:=\sum_{i=1}^g\xi_i\xi_{i+g}$ and $\xi:=\sum_{i=1}^g\zeta_i\zeta_{i+g}$, so that 
\begin{equation}\label{2theta}
\left(\sum_{i=1}^{2g} \xi_i\otimes e_i\right)^2 = -2\theta\otimes \omega, \;\;\; \left(\sum_{i=1}^{2g} \zeta_i\otimes e_i\right)^2 = -2\xi\otimes \omega, 
\end{equation}
and $\xi^g$ evaluates to $g!$ on $\Jac^{2g-2-l}(C)$. Finally, we set $\kappa:=\sum_{i=1}^g(\xi_i\zeta_{i+g}+\zeta_i\xi_{i+g})$.

In the next three Propositions, we rewrite the characteristic classes entering the integral \eqref{intwc1} in terms of the cohomology classes on $Z^0$ introduced above.

\begin{proposition}\label{prop:chNZ0}
In the notation introduced above, 
\begin{equation*}
e_{T\!\times\!\TT}(N^\vee_{Z^0/Z})^{-1} \!=\!\frac{ (-1)^l (s^{-1}t^{-2}e^{2\eta})^{(l-2m)}\exp(4\theta+9\xi-6\kappa-\frac{9\xi+\theta-3\kappa}{1-st^2e^{-\eta}} )}{(1-st^2e^{-2\eta})^{g-1+s_3^{\!\epsilon}}(1-s^{-1}t^{-2}e^{2\eta})^{g-s_3^{\!\epsilon}}(1-s^{-1}t^{-2}e^{\eta})^{g-l-m-s_2^{\!\epsilon}}}.
\end{equation*}
\end{proposition}
\begin{proof}
We use Lemma \ref{lem:normalZ0} to compute the Chern character of $N_{Z^0/Z}$. First, using \eqref{c1MFL}  and Lemma \ref{MFJweight}, we rewrite 
\begin{multline*}
-\ch_{T\!\times\!\TT}(\pi_!\mathrm{ParHom}(\mathcal{F}, \mathcal{M})) -\ch_{T\!\times\!\TT}(\pi_!\mathrm{ParHom}(\mathcal{M}, \mathcal{F})) = \\
-\ch_{T\!\times\!\TT}(\pi_!\mathrm{Hom}(\mathcal{F}, \mathcal{M})) -\ch_{T\!\times\!\TT}(\pi_! \mathrm{Hom}(\mathcal{M}, \mathcal{F}))+(1-s_3^{\!\epsilon})st^2e^{-2\eta}+s_3^{\!\epsilon}s^{-1}t^{-2}e^{2\eta}.
\end{multline*}
Applying \cite[Lemma 6.10]{OTASz}  to $S:=\mathrm{Hom}(\mathcal{M},\mathcal{F})$, with $T\!\times\!\TT$-character $st^2$, we obtain 
$$e_{T\!\times\!\TT}(-(\pi_!\mathrm{Hom}(\mathcal{M}, \mathcal{F})\oplus \pi_!\mathrm{Hom}(\mathcal{F}, \mathcal{M}))^\vee)^{-1} =  \frac{(-1)^l(st^2)^{2m-l}\exp(-\pi_*\ch_2(\mathcal{M}^\vee\otimes\mathcal{F}))}{(1-st^2e^{-2\eta})^{g-1}(1-s^{-1}t^{-2}e^{2\eta})^{g-1}}.$$
Using the relations \eqref{2theta}, the class $-\pi_*\ch_2(\mathcal{M}^\vee\otimes\mathcal{F})$ evaluates to
\begin{multline*}
-\frac{1}{2}\pi_*\bigl(-2\eta\otimes1 +(l-2m)\otimes\omega +\sum_{i=1}^{2g}(3\zeta_i-2\xi_i)\otimes e_i)\bigr)^2 = 
2\eta(l-2m) + 4\theta +9\xi -6 \kappa.
\end{multline*}

Next, applying the Grothendieck--Riemann--Roch theorem to the remaining term of the equivariant Chern character of $N_{Z^0/Z}$ and using \eqref{c1MFL} and \eqref{2theta}, we obtain 
\begin{multline*}
t\,\ch_{T\!\times\!\TT}(\pi_!\mathrm{SParHom}(\mathcal{M}, \mathcal{J}(p))) = 
 (2-g-s_2^{\!\epsilon})st^2e^{-\eta}+ st^2\pi_*\ch(\mathcal{M}^\vee\!\otimes\mathcal{J})= 
  \\ (g-l-m- s_2^{\!\epsilon})st^2e^{-\eta}\! +(- 9\xi -\theta +3\kappa)st^2e^{-\eta}.
\end{multline*}

Finally, applying Lemma \ref{thetato} to  the class $(- 9\xi -\theta +3\kappa)st^2e^{-\eta}$ and combining the result with the calculations above, we obtain the expression stated in Proposition \ref{prop:chNZ0}.
\end{proof}

\begin{proposition}\label{prop:normalN}
In the notation introduced above, 
\begin{multline}\label{eq:normalN}
e_{T\times\TT}\left({\mathcal N}_l^\vee\right)\!^{-1}\big|_{Z^0}=  \frac{(-1)^l(t^{-1}e^\eta)^{2l-m}(ste^{-\eta})^{l+m}\exp(9\xi+2\theta-3\kappa)}{(1-t^{-1})^{2g-2}(1-te^{-\eta})^{g-1+s_1^{\!\epsilon}}(1-t^{-1}e^{\eta})^{g-s_1^{\!\epsilon}}}\cdot\\
\frac{\exp(-\frac{9\xi+4\theta-6\kappa}{1-st^3e^{-2\eta}}-\frac{9\xi+4\theta-6\kappa}{1-s^{-1}t^{-1}e^{2\eta}})}{(1-s^{-1}t^{-1}e^{\eta})^{g-1+s_2^{\!\epsilon}}(1-ste^{-\eta})^{g-s_2^{\!\epsilon}}(1-s^{-1}t^{-3}e^{2\eta})^{g-1+s_3^{\!\epsilon}+l-2m}}\cdot \\
\frac{\exp(-\frac{9\xi+\theta-3\kappa}{1-s^{-1}e^{\eta}}-\frac{\theta}{1-t^2e^{-\eta}})}{(1-ste^{-2\eta})^{g-s_3^{\!\epsilon}+2m-l}(1-se^{-\eta})^{s_2^{\!\epsilon}+g-1+l+m}(1-t^{-2}e^{\eta})^{s_1^{\!\epsilon}+g-1+2l-m}}. 
\end{multline}
\end{proposition}
\begin{proof}
Recall that the Chern character of \(\mathcal N_l\) is given by \eqref{eq:Normal21}. To compute its restriction to \(Z^0\) (see  Construction \ref{descconstruction}), we note that the expression in \eqref{eq:Normal21} is independent of the choice of the normalization of the universal bundles \(\mathcal E\) and \(\mathcal L\). We may therefore choose the normalization so that the restriction of $\mathcal L$ to $N_l^{\alpha}({\hat c})\times\{p\}$ is trivial. In this case, using \eqref{eq:MFLdesc} and Lemma \ref{lem:Lvectriv}, we have 
\begin{equation}\label{eq:identif}
\mathcal{L}|_{Z^0}=\mathcal{J} \text{\; and \;} \mathcal{E}|_{Z^0}=\mathcal{M}\oplus\mathcal{F}.
\end{equation}
Lemmas \ref{lemma:Non21} and \ref{MFJweight} then give
\begin{multline}\label{eq:chNonZ}
\ch_{T\!\times\!\TT}(\mathcal{N}_l|_{Z^0})\!=\! -\ch_{T\!\times\!\TT}(\pi_!\mathrm{Hom}(\mathcal{M}\oplus\mathcal{F}, \mathcal{J}\otimes K^\vee\!))\! - \!\ch_{T\!\times\!\TT}(\pi_!\mathrm{Hom}(\mathcal{J}\otimes K^\vee\!, \mathcal{M}\oplus\mathcal{F})) + \!\\
t \ch_{T\!\times\!\TT}\pi_!(\mathrm{End}(\mathcal{M}\oplus\mathcal{F})\otimes K(p)) + t \ch_{T\!\times\!\TT}\pi_!(\mathrm{Hom}(\mathcal{J}\otimes K^\vee\!, \mathcal{M}\oplus\mathcal{F})\otimes K(p)) + \\ s_2^{\!\epsilon}ste^{-\eta} +(1-s_2^{\!\epsilon})s^{-1}t^{-1}e^{\eta} + s_1^{\!\epsilon}t^{-1}e^\eta + (1-s_1^{\!\epsilon})te^{-\eta}-(1-s_3^{\!\epsilon})st^3e^{-2\eta} - \\ s_3^{\!\epsilon}s^{-1}t^{-1}e^{2\eta} -(1-s_2^{\!\epsilon})e^\eta s^{-1} -(1-s_1^{\!\epsilon})t^2e^{-\eta}-2t.
\end{multline}

We now repeat the argument used in the proof of Proposition  \ref{prop:chNZ0}. Applying \cite[Lemma 6.10]{OTASz} and \eqref{2theta}, we find that the first two summands in \eqref{eq:chNonZ} contribute
$$\frac{(-1)^l\exp(9\xi+2\theta-3\kappa)  (t^{-1}e^\eta)^{2l-m}(ste^{-\eta})^{l+m}}{(1-ste^{-\eta})^{g-1}(1-s^{-1}t^{-1}e^\eta)^{g-1}(1-te^{-\eta})^{g-1}(1-t^{-1}e^\eta)^{g-1}}$$
to the inverse equivariant Euler class \eqref{eq:normalN}. 
Applying the Grothendieck--Riemann--Roch theorem to the third and fourth summands in \eqref{eq:chNonZ}, respectively, gives
\begin{multline*}
st^3e^{-2\eta}(l-2m+g -4\theta-9\xi+6\kappa) + s^{-1}t^{-1}e^{2\eta}(2m-l+g -4\theta-9\xi+6\kappa) + 2gt \\
\text{ \; and \;} 
s^{-1}e^\eta(m+l+g-\theta-9\xi+3\kappa) +t^2e^{-\eta}(2l-m+g-\theta).
\end{multline*}
Combining these summands with the remaining terms in \eqref{eq:chNonZ} and applying Lemma \ref{thetato} to the terms involving $\theta$, $\xi$ and $\kappa$  yields precisely \eqref{eq:normalN}.
\end{proof}

\begin{proposition} \label{prop:chLZ0}
\begin{enumerate}[label=(\roman*), leftmargin=1.5pc]
\item In the notation introduced above,
 \begin{multline*}
\ch_{T\!\times\!\TT}(\mathcal{L}_l(k;\lambda)|_{Z^0}) = (t^{-1}e^\eta)^{kl+ \lambda_{3-s_1^{\!\epsilon}-s_2^{\!\epsilon}}}(s^{-1}t^{-2}e^{2\eta})^{-km+\lambda_{2+s_2^{\!\epsilon}-s_3^{\!\epsilon}}}\exp(6k\xi+2k\theta-3k\kappa).
\end{multline*}
\item Recall the notation $|\nu|=\nu_1+\nu_2+\nu_3$,  $v^{\det}=\tfrac{|\nu|}3(1,1,-2)$ and $\widehat{\phi}^\nu(z_1,z_2)$ for the function $\phi^\nu(x)$ written in coordinates $z_1=\langle \alpha_{\bar{\epsilon} \epsilon},x\rangle$ and $z_2=\langle \alpha_{\jjo \bar{\epsilon}},x\rangle$. We have
\begin{multline*}
 \ch_{T\times\TT}\bigl(\pi_!(U_{l, \nu}\otimes K^{1/2})|_{Z^0}\bigr)  =e^{-v^{\det}_{2+s_2^{\!\epsilon}-s_3^{\!\epsilon}}(2u+v-2\eta)+v^{\det}_{3-s_1^{\!\epsilon}-s_2^{\!\epsilon}}(\eta-u)}
 (-(2m-l)\partial_{z_1} - (2l-m)\partial_{z_2} - \\ (9\xi+4\theta-6\kappa)\partial_{z_1}^2+(4\theta-3\kappa)\partial_{z_1}\partial_{z_2} -\theta \partial_{z_2}^2)\,\widehat\phi^\nu(z_1,z_2)\Big|_{\substack{z_1=2u+v-2\eta \\ z_2=\eta-u }}.
 \end{multline*}
 \end{enumerate}
\end{proposition}
\begin{proof}
Both statements follow from the same calculations as in the proofs of Propositions \ref{LklonF} and \ref{UnyuonF}, with  \(\mathcal L_1, \mathcal L_2, \mathcal L_3\) replaced by \(\mathcal M, \mathcal F, \mathcal J\otimes K^\vee\), respectively, and the corresponding Chern classes computed using \eqref{c1MFL} and \eqref{2theta}.  
\end{proof}
The last ingredient needed to evaluate the integral \eqref{intwc1} is the Todd class of \(Z^0\), which, using Lemmas \ref{lem:Z0descr} and \ref{l:intsym}(i), can be written as
\begin{equation}\label{toddZ}
\mathrm{Td}(Z^0) = \exp\left(\frac{\theta}{e^\eta-1}-\frac{\theta}{\eta}\right)\left(\frac{\eta}{1-e^{-\eta}}\right)^{g-s_1^{\!\epsilon}-2l+m}.
\end{equation}

\subsection{Computations of the wall-crossing term}\label{S:compWCR}
Armed with these results, we can now compute the wall-crossing difference \eqref{intwc1}. Recall that we denote by \(v\) and \(u\) the additive equivariant parameters of \(T\) and \(\TT\), respectively, so that \(s=e^v\) and \(t=e^u\). For \(\widehat\lambda\in\Lambda\otimes\mathbb Q\), consider the integral
\vspace{-3pt}
\begin{multline}\label{intwc2}
-\underset{v=-2u}{\res}\!\int_{Z^0}\!\frac{\eta^{g\!-\!s_1^{\!\epsilon}\!-\!2l\!+\!m}  e^{(k+3)(m(v+2u-2\eta)+l(\eta-u))}e^{(2\eta-v-2u)\lhat_{2+s_2^{\!\epsilon}-s_3^{\!\epsilon}}}e^{(\eta-u){\lhat_{3-s_1^{\!\epsilon}-s_2^{\!\epsilon}}}}  }
{W_{\!e^u\!}(\!-\!2\eta\!+\!v\!+\!2u,\eta\!-\!u)^{g-1}(1-e^{2u+v-2\eta})^{s_3^{\!\epsilon}}(1-e^{2\eta-v-2u})^{1-s_3^{\!\epsilon}}} \\
 \frac{\exp\!\left(-\tfrac{\theta}{\eta}+(9\xi\!+\!\theta\!-\!3\kappa)G_{12\!} + (9\xi\!+\!4\theta\!-\!6\kappa)G_{1\!} + \theta G_{2\!}\right)\phi(v\!+\!2u\!-\!2\eta)^{2m-l}}{(1-e^{u-\eta})^{s_1^{\!\epsilon}}(1-e^{\eta-u})^{1-s_1^{\!\epsilon}}(1-e^{\eta-u-v})^{s_2^{\!\epsilon}}(1-e^{-\eta+v+u})^{1-s_2^{\!\epsilon}}(1-e^{v-\eta})^{s_2^{\!\epsilon}} } \\
\frac{e^{-\delta((2m-l)\partial_{1}+(2l-m)\partial_{2})\widehat\phi^\nu(2u+v-2\eta,\eta-u)}\phi(\eta\!-\!u)^{2l-m}\phi(v\!+\!u\!-\!\eta)^{m+l} \,\mathrm{d}v}{(1-e^{\eta-2u-v})^{1- s_2^{\!\epsilon}}(1-e^{u+v-2\eta})^{1-s_3^{\!\epsilon}}(1-e^{2\eta-v-3u})^{s_3^{\!\epsilon}}(1-e^{-\eta})^{1-s_1^{\!\epsilon}}(1-e^{\eta-2u})^{s_1^\epsilon}},
\end{multline}
where the functions $W_{\!e^u\!}(z_1,z_2)$, $\phi(z)$, $\widehat{\phi}^\nu(z_1,z_2)$ and $G_{1\!}(z_1,z_2)$, $G_{2\!}(z_1,z_2)$, $G_{12\!}(z_1,z_2)$ are  defined in \S\ref{SS:Main}, and, for brevity, we set
\begin{itemize}
\item $G_{12}:=G_{12\!}(v+2u-2\eta,\eta-u),$ and, for $j=1,2$, 
\item $G_j:=G_{j\!}(v+2u-2\eta,\eta-u),$ \\
$ \partial_j\widehat{\phi}^\nu(2u+v-2\eta,\eta-u):= \partial_{z_j}\widehat{\phi}^\nu(z_1,z_2)|_{z_1=2u+v-2\eta,\, z_2=\eta-u}.$
\end{itemize}
Putting together Propositions  \ref{prop:chNZ0}, \ref{prop:normalN}, \ref{prop:chLZ0} and \eqref{toddZ},  
we obtain the following interpretation of \eqref{intwc2}. For $\delta=0$ and $\lhat=\lambda$, it is precisely the wall-crossing difference \eqref{intwc1} for
 $\mathcal{V}=\mathcal{L}_l(k;\lambda)$. On the other hand, the coefficient of $\delta$ in \eqref{intwc2}, with $\lhat=\lambda+v^{\det}$, gives the wall-crossing difference \eqref{intwc1} for $\mathcal{V}=\mathcal{L}_l(k;\lambda)\otimes\pi_!(U_{l,\nu}\otimes K^{1/2})$. 
 
We now evaluate \eqref{intwc2} explicitly. 
Applying Lemma \ref{l:intsym}(ii)  with 
\begin{itemize}
\item $B(\eta)=-\frac{1}{\eta}\!+\!G_{2\!}(2u+v-2\eta,\eta-u)\!+\!G_{12\!}(2u+v-2\eta, \eta-u)\!+\!4G_{1\!}(2u+v-2\eta,\eta-u)$, 
\item $b_i(\eta)=3\zeta_{i+g}(G_{12\!}(2u+v-2\eta, \eta-u)\!+\!2G_{1\!}(2u+v-2\eta, \eta-u))$, \\ $b_{i+g}(\eta)=-3\zeta_{i}(G_{12\!}(2u+v-2\eta, \eta-u)\!+\!2G_{1\!}(2u+v-2\eta, \eta-u))$,  $i=1,\ldots,g$, 
\end{itemize}
we obtain that \eqref{intwc2} is equal to 
\begin{multline*}
\underset{v=-2u}{\res}\,\underset{\eta=0}{\res}\int_{\mathrm{Jac}^{2g-2-l}(C)} \!\frac{- e^{(k+3)(m(v+2u-2\eta)+l(\eta-u))}e^{(2\eta-v-2u)\lhat_{2+s_2^{\!\epsilon}-s_3^{\!\epsilon}}}e^{(\eta-u){\lhat_{3-s_1^{\!\epsilon}-s_2^{\!\epsilon}}}}  }{W_{\!e^u\!}(\!-\!2\eta\!+\!v\!+\!2u,\eta\!-\!u)^{g-1}(1-e^{2u+v-2\eta})^{s_3^{\!\epsilon}}(1-e^{2\eta-v-2u})^{1-s_3^{\!\epsilon}}} \\
\frac{ \phi(v\!+\!2u\!-\!2\eta)^{2m-l} \exp\left(9\xi\left(G_{12} + G_{1} -\frac{(G_{12} + 2G_{1})^2}{G_{2}+G_{12}+4G_{1}} \right)\right)(G_{2}+ G_{12} + 4G_{1})^g}{(1-e^{u-\eta})^{s_1^{\!\epsilon}}(1-e^{\eta-u})^{1-s_1^{\!\epsilon}}(1-e^{\eta-u-v})^{s_2^{\!\epsilon}}(1-e^{-\eta+v+u})^{1-s_2^{\!\epsilon}}(1-e^{v-\eta})^{s_2^{\!\epsilon}}(1-e^{-\eta})^{1-s_1^{\!\epsilon}} } \\
\frac{e^{-\delta((2m-l)\partial_{1}+(2l-m)\partial_{2})\widehat\phi^\nu(2u+v-2\eta,\eta-u)} \phi(\eta-u)^{2l-m}  \phi(v+u-\eta)^{m+l} \,\mathrm{d}v\,\mathrm{d}\eta}{(1-e^{\eta-2u-v})^{1- s_2^{\!\epsilon}}(1-e^{u+v-2\eta})^{1-s_3^{\!\epsilon}}(1-e^{2\eta-v-3u})^{s_3^{\!\epsilon}}(1-e^{\eta-2u})^{s_1^\epsilon}}.
\end{multline*} 
Changing the variables $y=\eta-u$, $x=2u+v-2\eta$ and evaluating the integral over $\mathrm{Jac}^{2g-2-l}(C)$, we arrive at the following result.
\begin{proposition}\label{prop:intwcr3}
Fix \(\epsilon\in\{\iio,\iit\}\), an integer \(m\) satisfying $-2g+1+2l+s_1^{\!\epsilon}\le m<c_\epsilon,$ and let $\alpha(m,\epsilon)$ be the critical value defined in \eqref{sigmac}. Using the notation of \S\ref{SS:Main}, for $\lhat\in\Lambda\otimes\mathbb{Q}$, define 
\begin{multline*}
\Omega_l{(m;\epsilon)}\!:=\!\frac{- 3^{2g}e^{(k+3)(mx+ly)}e^{y{\lhat_{3-s_1^{\!\epsilon}-s_2^{\!\epsilon}}}-x\lhat_{2+s_2^{\!\epsilon}-s_3^{\!\epsilon}}} \phi(x)^{2m-l-s_3^{\!\epsilon}}\phi(y)^{2l-m+s_1^{\!\epsilon}} }{W_{\!e^u\!}(x,y)^{g-1}(1-e^{-x-y-u})(1-e^{x-u})(1-e^{-y-u})(1-e^{x})^{s_3^{\!\epsilon}}} \\
\frac{\phi(x+y)^{m+l+s_2^{\!\epsilon}} e^{-\delta((2m-l)\partial_x +(2l-m)\partial_y)\widehat{\phi}^\nu(x,y)}H(x,y)^g \, \mathrm{d}x\, \mathrm{d}y}{(1-e^{-x})^{1\!-\!s_3^{\!\epsilon}}  (1-e^{-y})^{s_1^{\!\epsilon}}(1-e^{y})^{1\!-\!s_1^{\!\epsilon}}(1-e^{-x-y})^{s_2^{\!\epsilon}}(1-e^{x+y})^{1\!-\!s_2^{\!\epsilon}}}, 
\end{multline*}
where $H(x,y)=G_{1\!}(x,y)G_{2\!}(x,y)+G_{1\!}(x,y)G_{12\!}(x,y)+G_{2\!}(x,y)G_{12\!}(x,y)$.
Let \(\alpha_-<\alpha(m,\epsilon)<\alpha_+\) lie  in the two adjacent chambers. Then  
\begin{equation*}
 \chi^{\alpha_-}_l[\mathcal{L}_l(k;\lambda)]-\chi^{\alpha_+}_l[\mathcal{L}_l(k;\lambda)] = \underset{x=0}{\res}\,\underset{y=-u}{\res} \Omega_l{(m;\epsilon)}\big|_{\delta=0,\, \lhat=\lambda},
\end{equation*}
while the wall-crossing difference 
\begin{equation*}
\chi^{\alpha_-}_l[\mathcal{L}_l(k;\lambda)\otimes \pi_!(U_{l,\nu}\otimes K^{1/2})]-\chi^{\alpha_+}_l[\mathcal{L}_l(k;\lambda)\otimes \pi_!(U_{l,\nu}\otimes K^{1/2})]
\end{equation*}
is given by the coefficient of $\delta$ in 
$$\underset{x=0}{\res}\,\underset{y=-u}{\res} \Omega_l{(m;\epsilon)}\big|_{\lhat=\lambda+v^{\det}}.$$
\end{proposition}

\subsection{Contribution of the type-$(2,1)$ fixed point components}\label{S:type21}
Following the strategy outlined in \S\ref{S:strategy}, we now compute the contribution of the fixed-point component \(N_l^{2g-2}(\hat c)\) to the equivariant Euler characteristics \eqref{Euler}. 
By  Corollary \ref{cor:max}, the moduli space \(N_l^\alpha(\hat c)\) is empty for \(\alpha>\alpha_{\max}\), while Proposition \ref{prop:alphacryt}(ii) shows that \(\alpha_{\max}\) is itself a critical value. Hence, putting together Propositions \ref{prop:intwcr3}, \ref{prop:tripletofix}, and  Remark \ref{rem:extremalwall}, we conclude that the contribution of \(N_l^{2g-2}(\hat c)\) to \eqref{Euler}, for $\mathcal V=\mathcal L_l(k;\lambda)$ or $\mathcal V=\mathcal L_l(k;\lambda)\otimes\pi_!(U_{l,\nu}\otimes K^{1/2}),$ is given by the relevant coefficients in the \(\delta\)-expansion of
\begin{equation}\label{eq:simple1}
 \sum_{\epsilon\in\{\iio,\iit\}}\sum_{m\geq -\!2g\!+\!1\!+\!2l\!+\!s_1^{\!\epsilon}}^{c_\epsilon} \underset{x=0}{\res}\,\underset{y=-u}{\res} \Omega_l{(m;\epsilon)},\end{equation}
evaluated at the appropriate value of \(\widehat\lambda\).
 
 In fact, the same argument as in \S\ref{S:111trick} allows us to extend the range of summation in \eqref{eq:simple1}. Indeed,  the order of the pole of $\Omega_l{(m;\epsilon)}$ at $y=-u$ is $$\underset{y =-u}{\mathrm{ord}}\Omega_l{(m;\epsilon)}\leq 2g-2l+m-s_1^{\!\epsilon}.$$ Thus, for $m<-2g+1+2l+s_1^{\!\epsilon}$, the form \(\Omega_l(m;\epsilon)\) has no pole at \(y=-u\), and its \(y\)-residue vanishes. We may therefore extend the \(m\)-sum in \eqref{eq:simple1} to obtain
\begin{equation}\label{eq:simple1ext}
\sum_{\epsilon\in\{\iio,\iit\}}\sum_{m< c_\epsilon}\!\!\!\underset{x=0}{\res}\,\underset{y=-u}{\res} \Omega_l{(m;\epsilon)}.
\end{equation}

To compute the total contribution of the fixed-point locus of parabolic type $(\{\iio,\iit\},\{\jjo\})$ to the equivariant Euler characteristics \eqref{Euler}, we sum \eqref{eq:simple1ext} over all values of $l$ indexing the connected components of this locus. By Proposition \ref{prop:tripletofix}, these are precisely  the integers satisfying
$-c_{\jjo}<l \leq g-1-\mathfrak{s} $, where $\mathfrak{s}= 1$ if $c_2<0$, $\jjo=1$ and $\mathfrak{s}= 0$ otherwise. Thus the total contribution is given by the $\delta$-coefficients of 
\begin{equation}\label{eq:simple11}
\sum_{l>-c_{\jjo}}^{g-1-\mathfrak{s}} \sum_{\epsilon\in\{\iio,\iit\}}  \sum_{m< c_\epsilon}\!\!\!\underset{x=0}{\res}\,\underset{y=-u}{\res} \Omega_l{(m;\epsilon)}.
\end{equation} 
evaluated at the appropriate value of \(\widehat\lambda\). A simple calculation shows that for $m<c_\epsilon$ and $l\geq g-\mathfrak{s}$ we again have  $$\underset{y =-u}{\mathrm{ord}}\Omega_l{(m;\epsilon)}\leq0,$$ so the corresponding \(y\)-residue vanishes. Hence the upper bound on \(l\) in \eqref{eq:simple11} may also be removed, and thus the total contribution of  the fixed-point components of parabolic type \((\{\iio,\iit\},\{\jjo\})\) to the equivariant Euler characteristic \eqref{Euler} is given by the coefficients in the \(\delta\)-expansion of
$$\sum_{\epsilon\in\{\iio,\iit\}}\sum_{l>-c_{\jjo}} \sum_{m< c_\epsilon}\!\underset{x=0}{\res}\,\underset{y=-u}{\res} \Omega_l{(m;\epsilon)}.$$ 

Now using Lemmas \ref{lem:phidef}, \ref{lem:Gdef}, and substituting the values of $s_i^{\!\epsilon}$, $i=1,2,3$, for the different choices $\{\iio, \iit\} \subset \{1,2,3\}$ and $\epsilon\in \{\iio,\iit\}$, we arrive at the following result.

\begin{theorem}\label{thm:contr21}
Let $\{\iio,\iit,\jjo\}=\{1,2,3\}$ and let $p\in\mathrm{Poles}(\A)$ be the pole with coordinates $(p_{\iio},p_{\iit},p_{\jjo})=(u/3,u/3,-2u/3)$.
For $\epsilon\in\{\iio,\iit\}$, denote by $\bar\epsilon$ the other element of $\{\iio,\iit\}$ and set $\bb_\epsilon:=(\alpha_{\bar{\epsilon} \epsilon }, \alpha_{\jjo \bar{\epsilon}})\in\Bases(\A_{p})$.
In the notation of Theorem \ref{thm:main}, the total contribution of fixed-point components of parabolic type $(\{\iio,\iit\},\{\jjo\})$ to the equivariant Euler characteristic of $\mathcal{L}(k;\lambda)$ on $\Hs_0(c)$ is 
$$
\sum_{\epsilon\in\{\iio,\iit\}}\! \iber_{p,\bb_\epsilon,Q_0}\!\!\!{}^{\!(-, +)}   \left[  \frac{3^g\cdot \det(\mathrm{Hess}_{\bb_\epsilon}(Q_0))^{g-1} \cdot e^{\langle \lambda, x\rangle}}{W_{\!e^u\!}(x)^{g-1}\prod_{i<j}(1-e^{\langle \alpha_{j i}, x \rangle})(1-e^{\langle \alpha_{i j}, x \rangle-u})} \right](-[c]_{\bb_\epsilon}),$$
while its contribution  to the equivariant Euler characteristic of \(\mathcal L(k;\lambda)\otimes \pi_!(U_\nu\otimes K^{1/2})\) is
$$  \left. \frac{\mathrm{d}}{\mathrm{d} \delta}\left(\sum_{\epsilon\in\{\iio,\iit\}}\! \iber_{p,\bb_\epsilon,Q_\delta}\!\!\!{}^{\!(-, +)}   \left[  \frac{3^g\cdot \det(\mathrm{Hess}_{\bb_\epsilon}(Q_\delta))^{g-1} \cdot e^{\langle \lambda+v^{\det}, x\rangle}}{W_{\!e^u\!}(x)^{g-1}\prod_{i<j}(1-e^{\langle \alpha_{j i}, x \rangle})(1-e^{\langle \alpha_{i j}, x \rangle-u})} \right](-[c]_{\bb_\epsilon})\right)\right|_{\delta=0}. $$
\end{theorem}

Applying Lemma \ref{rem:varChange}, with $\bb=(\alpha_{\bar{\epsilon} \epsilon}, \alpha_{\jjo \bar{\epsilon}})$ and $\bb'=(\alpha_{\epsilon \bar{\epsilon}}, \alpha_{\jjo \bar{\epsilon}}),$  to the formulas of Theorem \ref{thm:contr21}, we obtain the total contributions to \eqref{Euler} of the fixed-point components of parabolic type $(\{\iio, \iit\}, \{\jjo\})$ in the form stated in Theorem \ref{thm:main}. 

\subsection{Contribution of type-\((1,2)\) fixed-point components}\label{S:type12}
It remains to compute the contribution of the fixed-point components of type \((1,2)\) to the equivariant Euler characteristics \eqref{Euler}. Rather than treating these components separately, we reduce the calculation to the type-\((2,1)\) case by dualization. For \(c=(c_1,c_2,c_3)\in\Delta\), set \(c^{\!\vee}=(-c_3,-c_2,-c_1)\in\Delta\), and consider the \(\TT\)-equivariant isomorphism
$$h\colon \Hs_0(c)\xrightarrow{\sim} \Hs_0(c^{\!\vee}), \;\;\; (V,\phi)\longmapsto (V^\vee,-\phi^T)$$
(see \cite[Proposition 4.2 and Lemma 9.1]{GPGM}), where the bundle $V^\vee$ is equipped with the dual parabolic structure. More precisely, if $V_p=F_3\supset F_2\supset F_1\supset F_0=0$ is the parabolic flag on $V_p$ with weights $c$, then the corresponding flag on $V_p^\vee$ is  $$0=\check{F}_0\subset \check{F}_1\subset \check{F}_2\subset \check{F}_3=V_p^\vee, \;\;\;\;\;  \check{F}_i:= \ker\!\left(V_p^\vee\to F_{3-i}^\vee\right),$$ 
with parabolic weights $c^{\!\vee}$. 
Note that its successive quotients satisfy
\begin{equation}\label{eq:dual-flag-quotients}
{\check{F}_{4-i}}/{\check{F}_{3-i}}\simeq \left({F_i}/{F_{i-1}}\right)^\vee, \text{\; for \;}  i=1, 2, 3.
\end{equation}

\begin{proposition}\label{prop:duality12}
\begin{enumerate}[label=(\roman*), leftmargin=1.5pc]
\item
Under the isomorphism $h$ defined above, 
\[
h^*\mathcal L\bigl(k;(-\lambda_3,-\lambda_2,-\lambda_1)\bigr)
\simeq
\mathcal L(k;(\lambda_1,\lambda_2,\lambda_3)), \]
and for a dominant weight $\nu=(\nu_1,\nu_2,\nu_3)$  of $\mathrm{GL}_3$ with $|\nu|:=\nu_1+\nu_2+\nu_3$, we have 
\[h^*\bigl(\mathcal L(0;(|\nu|,0,-|\nu|))\otimes U_{(-\nu_3,-\nu_2,-\nu_1)}\bigr) = U_{(\nu_1,\nu_2,\nu_3)}.\]

\item
Moreover, for every $\{\iio,\iit,\jjo\}=\{1,2,3\}$, dualization induces a bijection between the connected components of $\Hs_0(c)^\TT$ of parabolic type $(\{\jjo\},\{\iio,\iit\})$ and the connected components of $\Hs_0(c^{\!\vee})^\TT$ of parabolic type $(\{4-\iit,4-\iio\},\{4-\jjo\}).$
\end{enumerate}
\end{proposition}

\begin{proof}
The first identity in (i) follows directly from \eqref{eq:dual-flag-quotients} and the definition of the line bundle $\mathcal L(k;\lambda)$. 
For the second identity, recall that, by a slight abuse of notation, we write \(U=\mathcal U\otimes\mathcal F_1^\vee\) for the normalized universal bundle on both \(\Hs_0(c)\) and \(\Hs_0(c^{\!\vee})\). Taking this normalization into account, we obtain $h^*(U\otimes\mathcal L(0;(-1,0,1)))\simeq U^\vee.$ Applying the representation of highest weight $\nu^\vee=(-\nu_3,-\nu_2,-\nu_1)$ to this isomorphism and using $U^\vee_{\nu^\vee}\simeq U_\nu$, we arrive at the second identity in (i).

For part (ii), dualization reverses the direction of the Higgs field, interchanging the rank-one and rank-two summands. Moreover, by \eqref{eq:dual-flag-quotients}, it sends each parabolic index \(i\) to \(4-i\), giving the stated correspondence of parabolic types.
\end{proof}

As a corollary, we obtain the following statement.
\begin{corollary}\label{cor:21to12}
The total contribution of the fixed-point components of parabolic type \((\{\jjo\},\{\iio,\iit\})\)  to the equivariant Euler characteristic of $\mathcal{L}(k;\lambda)$ on $\Hs_0(c)$ is equal to the total contribution of the fixed-point  components of parabolic type \((\{4-\iit,4-\iio\},\{4-\jjo\})\) to the equivariant Euler characteristic of $\mathcal{L}(k;(-\lambda_3,-\lambda_2,-\lambda_1))$ on $\Hs_0(c^{\!\vee})$.  
The same statement holds with $\mathcal L(k;\lambda)$ replaced by  $\mathcal L(k;\lambda)\otimes\pi_!(U_\nu\otimes K^{1/2})$, provided that on $\Hs_0(c^{\!\vee})$, the parameters  $\lambda$ and $\nu$ are replaced by $(-\lambda_3+|\nu|,-\lambda_2,-\lambda_1-|\nu|)$ and $(-\nu_3,-\nu_2,-\nu_1),$ respectively.
\end{corollary}
Corollary \ref{cor:21to12} thus reduces the computation of the total contribution to \eqref{Euler} of the fixed-point components of parabolic type \((\{\jjo\},\{\iio,\iit\})\) to the type-\((2,1)\) case treated in Theorem  \ref{thm:contr21}. To make this reduction explicit, consider the involution $$\iota\colon \mathfrak{t}\to \mathfrak{t}, \;\;\;\;\; \iota(a_1,a_2,a_3)=(-a_3,-a_2,-a_1).$$ 
The induced involution $\iota^*$ on $\mathfrak{t}^*$ satisfies $\iota^*x_i=-x_{4-i}$ for $i=1,2,3$.
In particular, $\iota$ sends the pole $p\in\mathrm{Poles}(\A)$ with coordinates $(p_{4-\iio},p_{4-\iit},p_{4-\jjo})=(u/3,u/3,-2u/3)$ to the pole $p'\in\mathrm{Poles}(\A)$ with coordinates $(p'_{\iio},p'_{\iit},p'_{\jjo})=(-u/3,-u/3,2u/3).$
Moreover, we have
\begin{multline*}
\iota^*(-\lambda_3,-\lambda_2,-\lambda_1)=(\lambda_1,\lambda_2,\lambda_3), \;\;  \iota^*(-\nu_3,-\nu_2,-\nu_1)=(\nu_1,\nu_2,\nu_3), \;\; \iota^*\alpha_{a b}=\alpha_{4-b \,4-a},
\end{multline*}
and $\iota^*c=c^{\!\vee}$.  Thus, applying the change of variables $x\mapsto\iota(x)$ to the formulas of Theorem \ref{thm:contr21}, we obtain the following result.

\begin{theorem}\label{thm:contr12}
Let $\{\iio,\iit,\jjo\}=\{1,2,3\}$ and let $p\in\mathrm{Poles}(\A)$ be the pole with coordinates $(p_{\iio},p_{\iit},p_{\jjo})=(-u/3,-u/3,2u/3)$.
For $\epsilon\in\{\iio,\iit\}$, denote by $\bar\epsilon$ the other element of $\{\iio,\iit\}$ and define $\bb_\epsilon:=(\alpha_{ \epsilon \bar{\epsilon}}, \alpha_{\bar{\epsilon} \jjo})\in\Bases(\A_{p})$. Then, in the notation of Theorem \ref{thm:main}, the total contribution of fixed-point components of parabolic type $(\{\jjo\}, \{\iio,\iit\})$ to the equivariant Euler characteristic of $\mathcal{L}(k;\lambda)$ on $\Hs_0(c)$ is 
$$
\sum_{\epsilon\in\{\iio,\iit\}}\! \iber_{p, \bb_\epsilon,Q_0}\!\!\!\!{}^{(-, +)}   \left[ \frac{3^{g}\cdot \det(\mathrm{Hess}_{\bb_\epsilon}(Q_0))^{g-1} \cdot e^{\langle \lambda, x\rangle}}{W_{\!e^u\!}(x)^{g-1}\prod_{i<j}(1-e^{\langle \alpha_{j i}, x \rangle})(1-e^{\langle \alpha_{i j}, x \rangle-u})} \right](-[c]_{\bb_\epsilon}),$$
while its contribution  to the equivariant Euler characteristic of \(\mathcal L(k;\lambda)\otimes \pi_!(U_\nu\otimes K^{1/2})\) is
$$  \left. \frac{\mathrm{d}}{\mathrm{d} \delta}\left( \sum_{\epsilon\in\{\iio,\iit\}}\! \iber_{p, \bb_\epsilon,Q_\delta}\!\!\!\!{}^{(-, +)}   \left[ \frac{3^{g}\cdot \det(\mathrm{Hess}_{\bb_\epsilon}(Q_\delta))^{g-1} \cdot e^{\langle \lambda+v^{\det}, x\rangle}}{W_{\!e^u\!}(x)^{g-1}\prod_{i<j}(1-e^{\langle \alpha_{j i}, x \rangle})(1-e^{\langle \alpha_{i j}, x \rangle-u})} \right](-[c]_{\bb_\epsilon})\right)\right|_{\delta=0}. $$
\end{theorem}
Finally, applying  Lemma \ref{rem:varChange}, with $\bb=(\alpha_{ \epsilon \bar{\epsilon}}, \alpha_{\bar{\epsilon} \jjo})$ and $\bb'=(\alpha_{\bar{\epsilon} \epsilon}, \alpha_{\bar{\epsilon} \jjo})$, to the formulas of Theorem \ref{thm:contr12}, we obtain the total contributions to \eqref{Euler} of the fixed-point components of parabolic type $(\{\jjo\}, \{\iio,\iit\})$ in the form stated in Theorem \ref{thm:main}. This completes the proof of Theorem \ref{thm:main}.

\appendix
\section{Intersection theory on symmetric products of a curve}\label{App1}
We recall several facts about intersection theory on symmetric products of a curve.
\renewcommand{\thelemma}{\Alph{section}.\arabic{lemma}}
\begin{lemma}\label{l:intsym} 
Let \(\mathcal{O}(\Delta)\) denote the line bundle over \(\Sym^m(C) \times C\) corresponding to the universal divisor \(\Delta\), and define cohomology classes \(\xi_i \in H^1(\Sym^m(C))\) for \(i = 1, \dots, 2g\), and \(\eta \in H^2(\Sym^m(C))\), as the components in the K\"unneth decomposition of \(c_1(\mathcal{O}(\Delta))\) (cf. \eqref{eq:ODelta}). 
\begin{enumerate}[label=(\roman*), leftmargin=1.5pc]
\item Define, in addition, $\theta:=\sum_{i=1}^g\xi_i\xi_{i+g}$. Then 
$$\mathrm{Td}(\Sym^{m}(C)) = \left(\frac{\eta}{1-e^{-\eta}}\right)^{m-g+1}\exp\left(\frac{\theta}{e^\eta-1}-\frac{\theta}{\eta}\right).$$
\item Let \({H}^*\) be a \(\mathbb{Z}\)-graded commutative algebra, and let \(b_i(z) \in {H}^{\mathrm{odd}} \otimes \mathbb{Q}[z]\) for \(i = 1, \dots, 2g\). Substituting \(z \mapsto \eta\), we regard \(b_i(\eta)\) as an element of \({H}^{\mathrm{odd}} \otimes H^*(\Sym^m(C))\), where the tensor product is taken in the graded sense. 
Then for any two power series \(A(z)\), \(B(z)\), we have
\begin{multline*}
\underset{{\Sym^{m}(C)}}{\int}\!\!\!A(\eta)\exp\!\left(\sum_{i=1}^{2g}b_i(\eta)\xi_i+B(\eta)\theta\!\right)=
\underset{z=0}{\mathrm{Res}} \frac{\exp\left(\!\frac{-\sum_{i=1}^gb_i(z)b_{i+g}(z)}{1/z+B(z)}\!\right)A(z)(1+zB(z))^g}{z^{m+1}}dz.
\end{multline*}
\end{enumerate}
\end{lemma}
\begin{proof} Part (i) follows from \cite[(7.3)]{Zagier}. For part (ii), we recall from  \cite[(7.2)]{Zagier} that $$\underset{{\Sym^{m}(C)}}{\int}\!\!A(\eta)\exp\left(B(\eta)\theta\right)= \underset{z=0}{\mathrm{Res}}\,\frac{A(z)({1+zB(z)})^g}{z^{m+1}}dz.$$ Let $0\leq n\leq g$. Substituting $B+\varepsilon$ for $B$, differentiating $n$ times with respect to $\varepsilon$, and setting $\varepsilon=0$, we obtain
\begin{equation}\label{Cor:Zag1}
\underset{{\Sym^{m}(C)}}{\int}\!\!A(\eta)\theta^n\exp\left(B(\eta)\theta\right)= \underset{z=0}{\mathrm{Res}}\,\frac{ g!\,A(z)({1+zB(z)})^{g-n}}{(g-n)!\,z^{m-n+1}}dz.
\end{equation}
Now given a set $S\subset\{1,2,...,g\}$ define $\xi_S:=\prod_{i\in S}\xi_i\xi_{i+g}$ and $b_S(z):= \prod_{i\in S}b_i(z)b_{i+g}(z)$. Then 
\begin{multline*}
\underset{{\Sym^{m}(C)}}{\int}\!\!A(\eta)\exp\left(\sum_{i=1}^{2g}b_i(\eta)\xi_i+B(\eta)\theta\right)= \!\!\underset{{\Sym^{m}(C)}}{\int}\!\!A(\eta)\exp\left(B(\eta)\theta\right) \!\!\underset{S\subset\{1,...,g\}}{\sum}\!(-1)^{|S|}b_S(\eta)\xi_S \!= \\
\underset{{\Sym^{m}(C)}}{\int}\!\!A(\eta)\exp\left(B(\eta)\theta\right) \!\!\underset{S\subset\{1,...,g\}}{\sum}\frac{(-1)^{|S|}b_S(\eta)\theta^{|S|}}{g(g-1)...(g-|S|+1)} \overset{\eqref{Cor:Zag1}}{=} \\
 \underset{z=0}{\mathrm{Res}} \underset{S\subset\{1,...,g\}}{\sum}\frac{(-1)^{|S|}b_S(z)}{(1/z+B(z))^{|S|}}   \frac{A(z)(1+zB(z))^g}{z^{m+1}}dz,
\end{multline*} 
where for the last equation we used \eqref{Cor:Zag1}. Finally, observing that $$\underset{S\subset\{1,...,g\}}{\sum}\frac{(-1)^{|S|}b_S(z)}{(1/z+B(z))^{|S|}}  = \exp\left(-\frac{\sum_{i=1}^gb_i(z)b_{i+g}(z)}{1/z+B(z)}\right),$$ we finish the proof of Lemma \ref{l:intsym}.
\end{proof}

\bibliographystyle{plain}
\bibliography{SL3Higgs.bib}

\begin{thebibliography}{10}

\bibitem{AGP}
J.~E. Andersen, S.~Gukov, and D.~Pei.
\newblock The {Verlinde} formula for {Higgs} bundles.
\newblock {\em arXiv:1608.01761}, 2016.

\bibitem{AB68}
M.~F. Atiyah and R.~Bott.
\newblock A {L}efschetz fixed point formula for elliptic complexes. {II}.
  {A}pplications.
\newblock {\em Ann. of Math. (2)}, 88:451--491, 1968.

\bibitem{Beck}
N.~Beck.
\newblock Moduli of decorated swamps on a smooth projective curve.
\newblock {\em Internat. J. Math.}, 26(10):1550086, 29, 2015.

\bibitem{ChHSz}
S.~Chiarello, T.~Hausel, and A.~Szenes.
\newblock An enumerative approach to {P=W}.
\newblock {\em arXiv:2002.08929}, 2020.

\bibitem{DonagiPantev}
R.~Donagi and T.~Pantev.
\newblock Langlands duality for {H}itchin systems.
\newblock {\em Invent. Math.}, 189(3):653--735, 2012.

\bibitem{DrezetNarasimhan}
J.-M. Dr{\'e}zet and M.~S. Narasimhan.
\newblock Groupe de {P}icard des vari\'et\'es de modules de fibr\'es
  semi-stables sur les courbes alg\'ebriques.
\newblock {\em Invent. Math.}, 97(1):53--94, 1989.

\bibitem{GPGM}
O.~Garc{\'i}a-Prada, P.~B. Gothen, and V.~Mu{\~n}oz.
\newblock Betti numbers of the moduli space of rank 3 parabolic {H}iggs
  bundles.
\newblock {\em Mem. Amer. Math. Soc.}, 187(879):viii+80, 2007.

\bibitem{GPLM}
O.~Garc{\'i}a-Prada, M.~Logares, and V.~Mu{\~n}oz.
\newblock Moduli spaces of parabolic {$U(p,q)$}-{H}iggs bundles.
\newblock {\em Q. J. Math.}, 60(2):183--233, 2009.

\bibitem{GothenOliveira}
P.~B. Gothen and A.~G. Oliveira.
\newblock Topological mirror symmetry for parabolic {H}iggs bundles.
\newblock {\em J. Geom. Phys.}, 137:7--34, 2019.

\bibitem{GP}
S.~Gukov and D.~Pei.
\newblock Equivariant {V}erlinde formula from fivebranes and vortices.
\newblock {\em Comm. Math. Phys.}, 355(1):1--50, 2017.

\bibitem{HL}
D.~Halpern-Leistner.
\newblock The equivariant {Verlinde} formula on the moduli of {Higgs} bundles.
\newblock {\em arXiv:1608.01754}, 2016.

\bibitem{HauselThaddeus}
T.~Hausel and M.~Thaddeus.
\newblock Mirror symmetry, {L}anglands duality, and the {H}itchin system.
\newblock {\em Invent. Math.}, 153(1):197--229, 2003.

\bibitem{MehtaSeshadri}
V.~B. Mehta and C.~S. Seshadri.
\newblock Moduli of vector bundles on curves with parabolic structures.
\newblock {\em Math. Ann.}, 248(3):205--239, 1980.

\bibitem{Okounkov}
A.~Okounkov.
\newblock Lectures on {$K$}-theoretic computations in enumerative geometry.
\newblock In {\em Geometry of moduli spaces and representation theory},
  volume~24 of {\em IAS/Park City Math. Ser.}, pages 251--380. Amer. Math.
  Soc., Providence, RI, 2017.

\bibitem{DecoSch}
A.~Schmitt.
\newblock A universal construction for moduli spaces of decorated vector
  bundles over curves.
\newblock {\em Transform. Groups}, 9(2):167--209, 2004.

\bibitem{Shen}
Sh. Shen.
\newblock Mirror symmetry for parabolic {H}iggs bundles via {$p$}-adic
  integration.
\newblock {\em Adv. Math.}, 443:Paper No. 109616, 40, 2024.

\bibitem{Simpson}
C.~T. Simpson.
\newblock Harmonic bundles on noncompact curves.
\newblock {\em J. Amer. Math. Soc.}, 3(3):713--770, 1990.

\bibitem{SzDu}
A.~Szenes.
\newblock Residue theorem for rational trigonometric sums and {V}erlinde's
  formula.
\newblock {\em Duke Math. J.}, 118(2):189--227, 2003.

\bibitem{OTASz}
A.~Szenes and O.~Trapeznikova.
\newblock The parabolic {V}erlinde formula: iterated residues and
  wall-crossings.
\newblock {\em Geom. Topol.}, 28(5):2259--2311, 2024.

\bibitem{TW}
C.~Teleman and C.~T. Woodward.
\newblock The index formula for the moduli of {$G$}-bundles on a curve.
\newblock {\em Ann. of Math. (2)}, 170(2):495--527, 2009.

\bibitem{Zagier}
M.~Thaddeus.
\newblock Stable pairs, linear systems and the {V}erlinde formula.
\newblock {\em Invent. Math.}, 117(2):317--353, 1994.

\bibitem{Thaddeus}
M.~Thaddeus.
\newblock Variation of moduli of parabolic {H}iggs bundles.
\newblock {\em J. Reine Angew. Math.}, 547:1--14, 2002.

\bibitem{Thomason}
R.~W. Thomason.
\newblock Une formule de {L}efschetz en {$K$}-th\'eorie \'equivariante
  alg\'ebrique.
\newblock {\em Duke Math. J.}, 68(3):447--462, 1992.

\bibitem{Unyu}
O.~Trapeznikova.
\newblock Tautological bundles on parabolic moduli spaces: {Euler}
  characteristics and {Hecke} correspondences.
\newblock {\em Algebraic Geom. (to appear), arXiv:2209.15600}, 2022.

\bibitem{Ver}
E.~Verlinde.
\newblock Fusion rules and modular transformations in {2D} conformal field
  theory.
\newblock {\em Nuclear Phys. B}, 300:360--376, 1988.

\bibitem{Yokogawa}
K.~Yokogawa.
\newblock Compactification of moduli of parabolic sheaves and moduli of
  parabolic {H}iggs sheaves.
\newblock {\em J. Math. Kyoto Univ.}, 33(2):451--504, 1993.

\end{thebibliography}

\end{document}